\documentclass[article,11pt,reqno]{amsart}

\usepackage{amsmath,amssymb,amsthm,bm}
\usepackage{caption,graphicx}
\usepackage[english]{babel}
\usepackage{amscd}
\usepackage{amsgen}
\usepackage[final]{epsfig}
\usepackage{latexsym}
\usepackage{amsfonts}
\usepackage{url}
\usepackage{slashed}
\usepackage[all]{xy}
\usepackage{here}
\usepackage{comment}
\usepackage{fixmath}
\usepackage{scalerel,stackengine}
\stackMath
\newcommand\reallywidehat[1]{%
\savestack{\tmpbox}{\stretchto{%
  \scaleto{%
    \scalerel*[\widthof{\ensuremath{#1}}]{\kern-.6pt\bigwedge\kern-.6pt}%
    {\rule[-\textheight/2]{1ex}{\textheight}}%WIDTH-LIMITED BIG WEDGE
  }{\textheight}%
}{0.5ex}}%
\stackon[1pt]{#1}{\tmpbox}%
}
\newtheorem{thm}{Theorem}
\newtheorem{lemm}{Lemma}
\newtheorem{prop}{Proposition}
\newtheorem{cor}{Corollary}

\newtheorem{defn}{Definition}

\allowdisplaybreaks

\begin{document}

\title{ Riemann surfaces of second kind and effective finiteness theorems }

\author{Burglind J\"oricke}

\address{Max-Planck-Institute for Mathematics\\
Vivatsgasse 7,
53111 Bonn\\
Germany}
%\address{IHES, 35 Route de Chartres\\ 91440 Bures-sur-Yvette\\ France}

\email{joericke@googlemail.com}

%\dedicatory{}

\keywords{finiteness theorems, Riemann surfaces of second kind,
$3$-braids, torus bundles, Gromov's Oka principle.}

\subjclass[2020]{Primary 32G13; Secondary 20F36, 32H35, 32Q56, 57Mxx}

\begin{abstract}
The Geometric Shafarevich Conjecture and the Theorem of de Franchis state the finiteness of the number of certain holomorphic objects on closed or punctured Riemann surfaces. The analog of these kind of theorems for Riemann surfaces of second kind is an estimate of the number of irreducible holomorphic objects up to homotopy (or isotopy, respectively). This analog can be interpreted as a quantitatve statement on the limitation for Gromov's Oka principle.

For any finite open Riemann surface $X$ (maybe, of second kind) we
give an effective upper bound for the number of irreducible holomorphic mappings up to homotopy from $X$ to the twice punctured complex plane, and an effective upper bound for the number of irreducible holomorphic torus bundles up to isotopy on such a Riemann surface.
The bound depends on a conformal invariant of the Riemann surface.

If $X_{\sigma}$ is the $\sigma$-neighbourhood of a skeleton of an open Riemann surface with finitely generated fundamental group,
%for any small enough number $\sigma>0$,
then the number of irreducible holomorphic mappings up to homotopy from $X_{\sigma}$ to the twice punctured complex plane grows exponentially in $\frac{1}{\sigma}$.
%This can be interpreted as a quantitative statement on the restricted validity of Gromov's %Oka principle.
\end{abstract}

\maketitle

\centerline \today
%\tableofcontents

\section{Introduction and statements of results}\label{sec:1}
%\label{Results}

It seems that the oldest finiteness theorem for mappings between complex manifolds is the following theorem, which was published by de Franchis \cite{Fr} in 1913.

\medskip

\noindent {\bf Theorem A (de Franchis).} {\it For closed connected Riemann surfaces $X$ and $Y$ with $Y$ of genus at least $2$ there are at most finitely many non-constant holomorphic mappings from $X$ to $Y$.}

\medskip

There is a more comprehensive Theorem in this spirit.

\medskip

\noindent {\bf Theorem B (de Franchis-Severi).} {\it For a closed connected Riemann surface $X$ there are
(up to isomorphism) only finitely many non-constant holomorphic mappings $f:X\to Y$ where $Y$ ranges over all closed Riemann surfaces of genus at least $2$.}

\medskip

A finiteness theorem which became more famous because of its relation to number theory was conjectured by Shafarevich  \cite{Sh}.

\medskip

\noindent {\bf Theorem C (Geometric Shafarevich Conjecture.)} {\it For a given compact or punctured Riemann surface $X$ and given non-negative numbers $\textsf{g}$ and $\textsf{m}$ such that $2\textsf{g}-2+\textsf{m} >0$ there are only finitely many locally holomorphically non-trivial holomorphic fiber bundles over $X$ with fiber of type $(\textsf{g},\textsf{m})$.}

\smallskip

A connected closed Riemann surface (or a smooth connected closed surface) is called of type $(\textsf{g},\textsf{m})$, if it has genus ${\sf{g}}$ and is equipped with ${\sf{m}}$ distinguished points.
Recall that a closed Riemann surface with a finite number of points removed is called a punctured Riemann surface. The removed points are called punctures. Sometimes it is convenient to associate a punctured Riemann surface to a Riemann surface of type $(\textsf{g},\textsf{m})$ by removing the distinguished points.
A Riemann surface is called finite if its fundamental group is finitely generated, and open if no connected component is compact.
A finite connected Riemann surface is called of first kind, if it is a closed or a punctured Riemann surface, otherwise it is called of second kind.

Each finite connected open Riemann surface $X$ is conformally equivalent
to a domain (denoted again by $X$) on a closed Riemann surface $X^c$ such that each connected component of the complement $X^c \setminus X$ is either a point or a closed topological disc with smooth boundary \cite{Sto}.
The connected components of the complement will be called holes. A finite Riemann surface $X$ is of first kind, if and only if all connected components of $X^c \setminus X$
are points. We will say that a connected finite open Riemann surface has only thick ends if all connected components of $X^c \setminus X$ are closed topological discs.

Each finite Riemann surface whose universal covering is equal to the upper half-plane $\mathbb{C}_+$ (a finite hyperbolic Riemann surface for short) is conformally equivalent to the quotient of $\mathbb{C}_+$ by a Fuchsian group. The Riemann surface is of first kind if and only if the Fuchsian group is of first kind (\cite{Kra}, II, Theorem 3.2). %p.52
We will not make use of Fuchsian groups here.

\medskip

Theorem C was conjectured by Shafarevich \cite{Sh} in the case of compact base and fibers of type $(\textsf{g},0)$. It was
proved by Parshin \cite{Pa} in the case of compact base and fibers of type $(\textsf{g},0),\,\textsf{g}\geq 2,$ and by Arakelov \cite{Ar} for punctured Riemann surfaces as base and fibers of type $(\textsf{g},0)$. Imayoshi and Shiga \cite{IS} gave a proof of the quoted version using Teichm\"uller theory.

The statement of Theorem C ''almost'' contains the so called Finiteness Theorem of Sections which is also called
the Geometric Mordell conjecture (see \cite{Mc}),
%in the function field case (Manin's Theorem)
giving an important conceptional connection between geometry and number theory.
For more details we refer to the surveys by C.McMullen \cite{Mc} and B.Mazur \cite{Ma}.

%There is an older finiteness theorem due to de Franchis \cite{Fr} which can be %considered as a special case of the Geometric Shafarevich conjecture.

\medskip
Theorem A is a consequence of Theorem C, and Theorem A has analogs for the source $X$
%being a punctured Riemann surface,
and the target $Y$ being punctured Riemann surfaces. Indeed,
we may associate to any holomorphic mapping $f:X\to Y$ of Theorem A the bundle over $X$ with fiber over $x\in X$ equal to $Y$ with distinguished point $\{f(x)\}$. Thus, the fibers are of type $(\textsf{g},1)$.
A holomorphic self-isomorphism of a locally holomorphically non-trivial $(\textsf{g},1)$-bundle may lead to a new holomorphic mapping from $X$ to $Y$, but there are only finitely many different holomorphic self-isomorphisms.

We will consider here analogs of Theorems A and C for the case when the base $X$ is a Riemann surface of second kind. Notice that finite hyperbolic Riemann surfaces of second kind are interesting from the point of view of spectral theory of the Laplace operator with respect to the hyperbolic metric (see also \cite{Bo}). There are interesting relations to scattering theory and
(the Hausdoff dimension of) the limit set of the Fuchsian group defining $X$.

The Theorems A and C do not hold literally if the base $X$ is of second kind.
If the base is a Riemann surface of second kind the problem to be considered is the finiteness of the number of irreducible isotopy classes (homotopy classes, respectively) containing holomorphic objects. In case the base is a punctured Riemann surface this is equivalent to the finiteness of the number of holomorphic objects. For more detail see Sections \ref{sec:2} and \ref{sec:3}.

We will prove finiteness theorems with effective estimates for Riemann surfaces of second kind.
The estimates depend on a conformal invariant of the base manifold. To define the invariant we recall Ahlfors' definition of extremal length
(see \cite{A1}).
For an annulus $A=\{0\leq r<|z|<R\leq \infty\}$ (and for any open set that is conformally equivalent to $A$) the extremal length equals $\frac{2\pi}{\log\frac{R}{r}}$.
For an open rectangle
$R= \{z=x+iy:0<x<{\sf b},\,0<y<{\sf a}\,\}$ in the plane
with sides parallel to the axes, and with horizontal side length $\textsf{b}$ and vertical side length $\textsf{a}$ the extremal length equals $\lambda(R)=\frac{\sf a}{\sf b}$.
For a conformal mapping $\omega:R\to U$ of the rectangle $R$ onto a domain
$U\subset \mathbb{C}$ the image $U$ is called a curvilinear rectangle, if $\omega$ extends to a continuous mapping on the closure $\bar R$, and the restriction to each (closed) side of $R$ is a homeomorphism onto its image. The images of the vertical (horizontal, respectively) sides of $R$ are called the vertical (horizontal, respectively) curvilinear sides of the curvilinear rectangle $\omega(R)$. The extremal length of the curvilinear rectangle $U$ equals the extremal length of $R$. (See \cite{A1}).

\begin{figure}[h]
\begin{center}
\includegraphics[width=80mm]{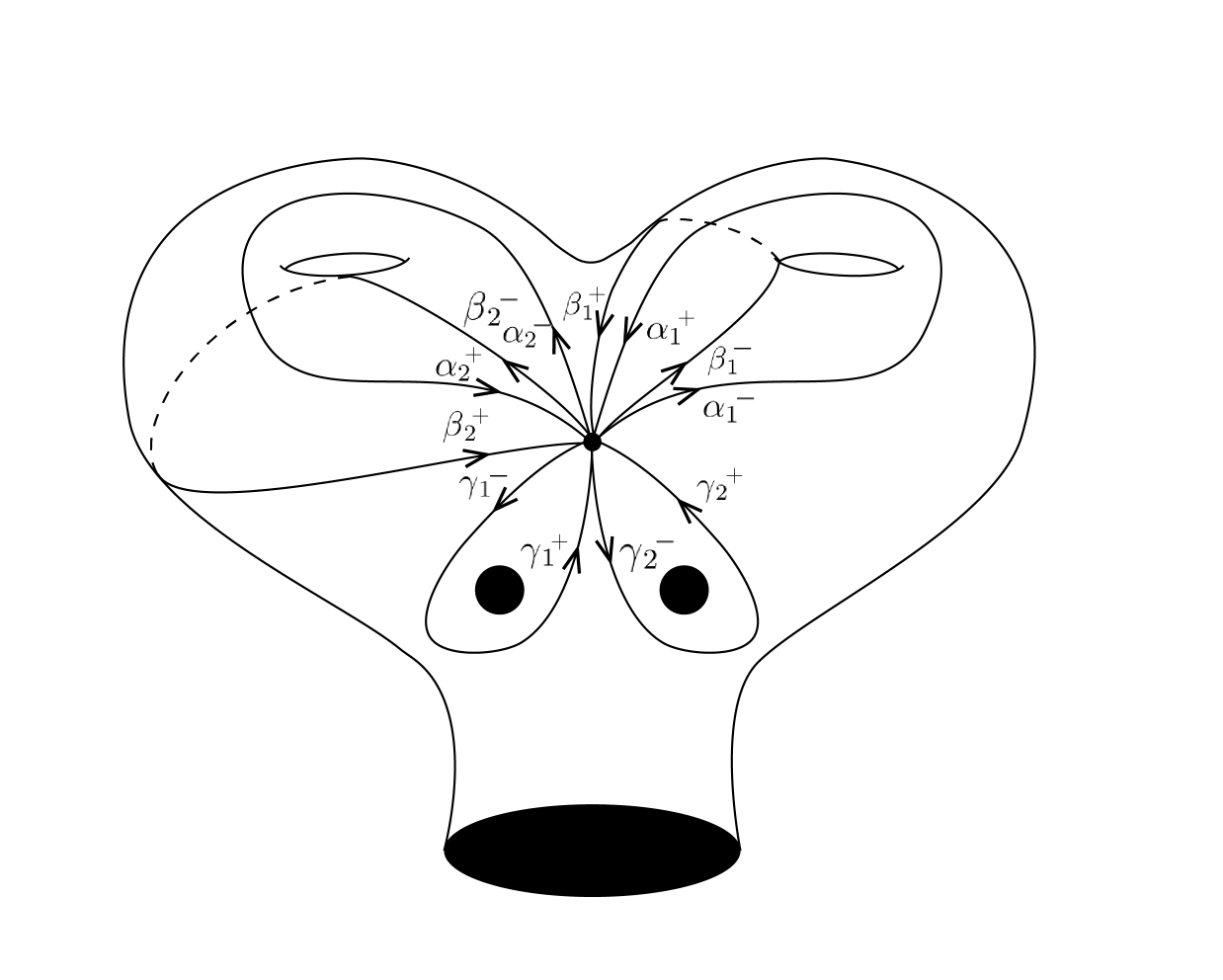}
\end{center}
\caption{A standard bouquet of circles for a connected finite open Riemann surface}\label{fig8.1}
\end{figure}

Let $X$ be a connected open Riemann surface of genus $g\geq 0$ with $m+1$ holes, $m\geq 0$, equipped with a base point $q_0$.
The fundamental group $\pi_1(X,q_0)$ of $X$ is a free group in $2g+m$ generators.
We describe now the conformal invariant of the Riemann surface $X$ that will appear in the mentioned estimate. We take a bouquet of non-contractible circles $S$ in $X$ with base point $q_0$, %which is a deformation retract of $X$, (in particular,
such that $q_0$ is the only common point of any pair of circles in $S$. Moreover, $S$ is the union
of simple closed oriented curves $\alpha_j,\,\beta_j$, $ j=1,\ldots ,g',$  and  $\gamma_{k}$, $k=1,\ldots,m'$, with base point $q_0$ with the following property.
%such that no pair of curves intersects outside $q_0$. Moreover,
Labeling the rays of the loops emerging from the base point $q_0$ by $\alpha_j^-,\,\beta_j^-$  $\gamma_j^-$, and the incoming rays by $\alpha_j^+,\,\beta_j^+$  $\gamma_j^+$, we require that
when moving in counterclockwise direction along a small circle around $q_0$ we meet the rays in the order
\begin{align*}
\ldots, \alpha^-_j,\beta^-_j,\alpha^+_j,\beta^+_j,\ldots, \gamma^-_k,\gamma^+_k,\ldots\;.
\end{align*}
(See Figure \ref{fig8.1}.)
We call such a bouquet of circles a standard bouquet of circles contained in $X$.
If the collection $\mathcal{E}$ of elements of the fundamental group $\pi_1(X,q_0)$ represented by the collection of curves in $S$ is a system of generators of $\pi_1(X,q_0)$ (then in particular, $g'=g$, $m'=m$), we call $S$ a standard bouquet of circles for $X$, and say that the system $\mathcal{E }$ is associated to a standard bouquet of circles for $X$.

The existence of a standard bouquet of circles for a connected finite open Riemann surface can be seen by looking at a fundamental polygon of the compact Riemann surface $X^c$ that contains a lift of each hole of $X$. The pairs of curves $\alpha_j$, $\beta_j$ correspond to the handles of $X^c$. Each curve $\gamma_k,\, k=1,\ldots,m$, surrounds a connected component $\mathcal{C}_k$ of $X^c\setminus X$ counterclockwise. More precisely, $\gamma_{k}$ is contractible in $X\cup \mathcal{C}_k$ and divides $X$ into two connected
components, one of them containing $\mathcal{C}_k$. Moreover, moving along $\gamma_k$ we see $\mathcal{C}_k$ on the left.

Vice versa, if a connected open Riemann surface $X$ contains a standard bouquet of circles consisting of $g$ pairs of curves $\alpha_j$,$\beta_j$, and $m$ curves $\gamma_k$ as above, that represent a system of generators of $\pi_1(X,q_0)$, then $X$ has genus $g$ and $m$ holes. To see this we cut the compact Riemann surface $X^c$ along the $\alpha_j$, $\beta_j$ and obtain a fundamental polygon which corresponds to a closed Riemann surface of genus $g$. The $\gamma_k$ are contractible in $X^c$, hence, each of them surrounds a hole.

Label the generators $\mathcal{E}\subset  \pi_1(X,q_0)$ of a standard bouquet of circles for $X$ as follows.
%associated to a standard bouquet of circles for $X$. It
The elements $e_{2j-1,0}\in \pi_1(X,q_0),$  $j=1,\ldots ,g,$  are represented by $\alpha_j$, the elements $e_{2j,0}\in \pi_1(X,q_0) , \, j=1,\ldots ,g,$  are represented by $\beta_j$, and the elements $e_{2g+k,0}\in \pi_1(X,q_0), \, k=1,\ldots ,m,$  of $\pi_1(X,q_0)$ are represented by $\gamma_k$.
A standard bouquet of circles for a connected finite open Riemann surface is a deformation retract of $X$.
We %denote the described system of generators by $\mathcal{E}$, and
will fix the system of generators $\mathcal{E}$ of $\pi_1(X,q_0)$ throughout the paper.

Let $\tilde X$ be the universal covering of $X$.
For each element $e_0 \in \pi_1(X,q_0)$ we consider the subgroup $\langle e_0\rangle$ of $\pi_1(X,q_0)$ generated by $e_0$.
%Let $\tilde X$ be the universal covering of $X$.
Let $\sigma(e_0)$ be the covering transformation corresponding to $e_0$, and $\langle \sigma( e_0) \rangle$ the group generated by $\sigma(e_0)$.
\begin{defn}\label{defn0}
Denote by $\mathcal{E}_j,\, j=2,\ldots,10,$ the set of primitive elements of $\pi_1(X,q_0)$ which can be written as product of at most $j$ factors with each factor being either an element of $\mathcal{E}$ or an element of $\mathcal{E}^{-1}$, the set of inverses of elements of $\mathcal{E}$. Define $\lambda_j=\lambda_j(X)$ as the maximum over $e_0 \in \mathcal{E}_j$ of the extremal length of the annulus $\tilde{X} \diagup \langle \sigma( e_0) \rangle$.
\end{defn}
The quantity
%$\lambda= \lambda(X)\stackrel{def}=\lambda_5(X)$
$\lambda_7(X)$ (for mappings to the twice punctured complex plane), or $\lambda_{10}(X)$ (for $(1,1)$-bundles) is the mentioned conformal invariant.

\medskip

Let $E$ be a finite subset of the Riemann sphere $\mathbb{P}^1$ which contains at least three points. Let $X$ be a finite open Riemann surface with non-trivial fundamental group.  %non-contractible
A continuous map
$f:X \to  \mathbb{P}^1\setminus E$ is reducible if it is homotopic (as a mapping to $\mathbb{P}^1\setminus E$) to a mapping whose image is contained in $D\setminus E$ for an open topological disc $D\subset \mathbb{P}^1$
with $E\setminus D$
containing at least two points of $E$.
Otherwise the mapping is called irreducible.

In the following theorem we take $E=\{-1,1,\infty\}$. We will often refer to $\mathbb{P}^1\setminus \{-1,1,\infty\}$ as the thrice punctured Riemann sphere or the twice punctured complex plane $\mathbb{C}\setminus \{-1,1\}$. Note that a continuous
%non-contractible
mapping from a Riemann surface to the twice punctured complex plane is reducible, iff it is homotopic to a mapping with image in a once punctured disc contained in $\mathbb{P}^1\setminus E$. (The puncture may be equal to $\infty$.)
There are countably many non-homotopic reducible holomorphic mappings with target being the twice punctured complex plane and source being any finite open Riemann surface with only thick ends and non-trivial fundamental group (see the proof of Lemma 15 in \cite{Jo5}).
On the other hand the following theorem holds.
\begin{thm}\label{thm1}
For each open connected Riemann surface $X$ of genus $g\geq0$ with $m+1\geq 1$ holes there are up to homotopy at most
$3(\frac{3}{2}e^{24\pi \lambda_7(X)})^{2g+m}$
irreducible holomorphic mappings from $X$ into
$Y\stackrel{def}=\mathbb{P}^1\setminus \{-1,1,\infty\}$.
\end{thm}
Notice that the Riemann surface $X$ is allowed to be of second kind.
If $X$ is a torus with a hole, $\lambda_7(X)$ may be replaced by $\lambda_3(X)$.
If $X$ is a planar domain, $\lambda_7(X)$ may be replaced by $\lambda_4(X)$
\smallskip

A holomorphic $(1,0)$-bundle is also called a holomorphic torus bundle. A holomorphic torus bundle equipped with a holomorphic section is also considered as a holomorphic $(1,1)$-bundle. The following lemma holds.

\medskip

\noindent {\bf Lemma D.} {\it A smooth
$(0,1)$-bundle admits a smooth section. A holomorphic torus bundle
is (smoothly) isotopic to a holomorphic torus bundle that admits a holomorphic section.}

\medskip

\noindent For a proof see \cite{Jo5}.

\begin{thm}\label{thm2}
Let $X$ be an open connected Riemann surface of genus $g\geq 0$ with $m+1\geq 1$ holes.
Up to isotopy
there are no more than
%$2^{2g+m}\big(3^3(\exp(6 \pi \lambda_{10}(X)))\big)^{6(2g+m)}=$
$\big(2 \cdot 15^6\cdot\exp(36 \pi \lambda_{10}(X))\big)^{2g+m}$
%$\big(2\cdot 3^6\cdot 5^6\cdot\exp(36 \pi \lambda_8(X))\big)^{2g+m}$
irreducible holomorphic $(1,1)$-bundles over $X$.
\end{thm}

For the definition of irreducible $(\sf g,\sf m)$-bundles
%and a remark treating reducible bundles
see Section \ref{sec:3} below.
%See also the remark in Section \ref{sec:3} concerning
%reducible $(1,1)$-bundles.
Since on each finite open Riemann surface with only thick ends and non-trivial fundamental group there are
countably many non-homotopic reducible holomorphic mappings with target being the twice punctured complex plane, there
are also countably many non-isotopic holomorphic $(1,1)$-bundles over each such Riemann
surface (see Proposition \ref{prop1} below).

We wish to point out that
reducible $(\sf{g,m})$-bundles over finite open Riemann surfaces
can be decomposed into irreducible bundle components, and each reducible bundle is determined by its bundle components up to commuting Dehn twists in the fiber over the base point. (For details see \cite{Jo5}.)

Notice that Caporaso proved the existence of a uniform bound for the number of objects in Theorem C in case $X$ is a closed Riemann surface of genus $g$ with $m$ punctures, and the fibers are closed Riemann surfaces of genus $ \textsf{g} \geq 2$. The bound depends only on the numbers $g$, $\textsf{g}$ and $m$.  Heier gave effective uniform estimates, but the constants are huge and depend in a complicated way on the parameters.

Theorems \ref{thm1} and \ref{thm2} imply effective estimates for the number of locally holomorphically non-trivial holomorphic $(1,1)$-bundles over punctured Riemann surfaces, however, the constants depend also on the conformal type of the base. More precisely, the following corollaries hold.
\begin{cor}\label{cor1a}
There are no more than
$3(\frac{3}{2}e^{24 \pi \lambda_7(X)})^{2g+m}$
non-constant holomorphic mappings from a Riemann surface $X$ of type $(g,m+1)$ to $\mathbb{P}^1\setminus \{-1,1,\infty\}$.
\end{cor}

\begin{cor}\label{cor1b}
There are no more than
$\big(2 \cdot 15^6\cdot\exp(36 \pi \lambda_{10}(X))\big)^{2g+m}$
%$\big(2\cdot 3^6\cdot 5^6\cdot\exp(36 \pi \lambda_8(X))\big)^{2g+m}$
locally holomorphically non-trivial holomorphic $(1,1)$-bundles over a Riemann surface $X$ of type $(g,m+1)$.
\end{cor}

The following examples demonstrate the different nature of the problem in the two cases, the case when the base is a punctured Riemann surface, and when it is a Riemann surface of second kind.

\medskip

\noindent {\bf Example 1.}
There are no non-constant holomorphic mappings from a torus with one puncture to the twice punctured complex plane. Indeed, by Picard's Theorem each such mapping extends to a meromorphic mapping from the closed torus to the Riemann sphere. This implies that the preimage of the set $\{-1,1,\infty\}$ under the extended mapping must contain at least three points, which is impossible.

The situation changes if $X$ is a torus with a large enough hole.
Let $\alpha\geq 1$ and $\sigma \in (0,1)$. Consider the torus with a hole $T^{\alpha,\sigma}$ that is obtained from $\mathbb{\mathbb{C}}\diagup (\mathbb{Z} + i \alpha \mathbb{Z}),$ (with $\alpha\geq 1$ being a real number) by removing a closed geometric rectangle of vertical side length $\alpha-\sigma$ and horizontal side length $1-\sigma$ (i.e. we remove a closed subset that lifts to such a closed rectangle in $\mathbb{C}$). A fundamental domain for
this Riemann surface is ''the golden cross on the Swedish flag'' turned by $\frac{\pi}{2}$ with width of the laths being $\sigma$ and length of the laths being $1$ and $\alpha$.

\begin{prop}\label{prop1a}
Up to homotopy there are at most
$7 e^{3\cdot 2^4 \pi \frac{2\alpha +1}{\sigma}}$
irreducible holomorphic mappings from  $T^{\alpha,\sigma}$
to the twice punctured complex plane.

On the other hand, there are positive constants $c$, $C$, and $\sigma_0$ such that for any positive number $\sigma<\sigma_0$ and any $\alpha >1$ there are at least $ce^{C\frac{\alpha}{\sigma}}$ non-homotopic holomorphic mappings from $T^{\alpha,\sigma}$
to the twice punctured complex plane.

\end{prop}

\noindent {\bf Example 2.}
There are only finitely many holomorphic maps from a thrice punctured Riemann sphere to another thrice punctured Riemann sphere. Indeed, after normalizing both, the source and the target space,
by a M\"obius transformation we may assume that both
are equal to $\mathbb{C} \setminus \{-1,1\}$. Each holomorphic map from
$\mathbb{C} \setminus \{-1,1\}$ to itself extends to a meromorphic map from the Riemann sphere to itself, which maps the set $\{-1,1,\infty\}$ to itself and maps no other point to this set. By the Riemann-Hurwitz formula
the meromorphic map takes each value exactly once. Indeed, suppose it takes each value $l$ times for a natural number $l$. Then each point in $\{-1,1,\infty\}$ has ramification index $l$.
Apply the Riemann Hurwitz formula for the branched covering $ X=\mathbb{P}^1\to Y=\mathbb{P}^1$ of multiplicity $l$
$$
\chi(X)= l \cdot \chi(Y) -\sum_{x\in Y} (e_x-1).
$$
Here $e_x$ is the ramification index at the point $x$. For the Euler characteristic we have $\chi(\mathbb{P}^1)=2$, and $\sum_{x\in Y} (e_x-1)\geq \sum_{x=-1,1,\infty} (e_x-1)=3\,(l-1)$. We obtain $2\leq 2\,l-3\,(l-1)$ which is possible only if $l=1$.
We saw that each non-constant holomorphic mapping from $\mathbb{C}\setminus\{-1,1\}$ to itself extends to a conformal mapping from
the Riemann sphere to the Riemann sphere that maps the set $ \{-1,1,\infty\}$ to itself. There are only finitely may such maps, each a M\"obius transformation commuting the three points.

For Riemann surfaces of second kind
%which are equal to $\mathbb{P}^1$ with three holes
the situation changes, as demonstrated in the following proposition. The proposition does not only concern the
case when the Riemann surface equals $\mathbb{P}^1$ with three holes.
We consider an open Riemann surface $X$ of genus $g$ with $m\geq 1$ holes.

\begin{prop}\label{prop1b} Let $X$ be a connected finite open hyperbolic Riemann surface, that is equipped with a K\"ahler metric. Suppose $S$ is a standard bouquet of piecewise smooth circles in $X$
%with base pointpairs of circles from $S$.
with base point $q_0$.
% which is a deformation retract of $X$. We assume that $q_0$ is the only non-smooth point of the circles and the only intersection point of  circles from $S$.
We assume that $q_0$ is the only non-smooth point of the circles, and
all tangent rays to circles in $S$ at $q_0$ divide a disc in the tangent space into equal sectors.
Let $S_{\sigma}$ be the
%$\frac{1}{2}\sigma$
$\sigma$-neighbourhood of $S$ (in the K\"ahler metric on $X$).

Then there exists a constant $\sigma_0>0$, and positive constants $C'$, $C''$,  $c'$, $c''$,
depending only on $X$, $S$ and the K\"ahler metric, such that for each positive $\sigma<\sigma_0$ the number $N_{S_{\sigma}}^{\mathbb{C}\setminus \{-1,1\}}$ of non-homotopic irreducible holomorphic mappings from $S_{\sigma}$ to the twice punctured complex plane satisfies the
inequalities
\begin{align}\label{eqabc}
c'e^{\frac{c''}{\sigma}} \leq N_{S_{\sigma}}^{\mathbb{C}\setminus \{-1,1\}}\leq C'e^{\frac{C''}{\sigma}}\,.
\end{align}
\end{prop}

The present results may be understood as quantitative statements with regard to limitations for Gromov's Oka principle.
Gromov  \cite{G} formulated his Oka principle as "an expression of an optimistic expectation with regard to the validity of the $h$-principle for holomorphic maps in the situation when the source manifold is Stein". Holomorphic maps $X\to Y$ from a
complex manifold $X$ to a complex manifold $Y$ are said to satisfy the $h$-principle if each continuous map from $X$ to $Y$ is homotopic to a holomorphic map. We call a target manifold $Y$ a Gromov-Oka manifold if the $h$-principle holds for holomorphic maps from any Stein manifold to $Y$.
Gromov \cite{G} gave a sufficient condition on a complex manifold $Y$ to be a Gromov-Oka manifold.

The question of understanding Gromov-Oka manifolds received a lot of attention. It turned out to be fruitful to strengthen the requirement on the target $Y$ by combining the
$h$-principle for holomorphic maps
with a holomorphic approximation property. Manifolds $Y$ satisfying the stronger requirement are called Oka manifolds. For details and an account on modern development of Oka theory based on Oka manifolds see \cite{Forst}.

The twice punctured complex plane %$Y\stackrel {def}=$
$\mathbb{C}\setminus \{-1,1\}$ is not a Gromov-Oka manifold. Then the question becomes, what prevents a continuous map from a Stein manifold $X$ to $\mathbb{C}\setminus \{-1,1\}$  to be homotopic to a holomorphic map, and ''how many'' homotopy classes contain a holomorphic map? As for the first question in case the source manifold is a finite open  Riemann surface $X$, Proposition \ref{prop2} below says that an irreducible map $X\to\mathbb{C}\setminus \{-1,1\}$ can only be homotopic to a holomorphic map, if the ''complexity'' of the monodromies of the map are compatible with conformal invariants of the source manifold.
Theorem \ref{thm1} gives an upper bound related to the second question.
Propositions \ref{prop1a} and \ref{prop1b} can be interpreted as statements related to the following question. Consider a family of Riemann surfaces $Y_{\sigma},\,\sigma\in (0,\sigma_0) $, obtained by continuously changing the conformal structure of a fixed Riemann surface. Determine the growth rate for $\sigma\to 0$ of the number of irreducible holomorphic mappings $X_\sigma\to \mathbb{C}\setminus\{-1,1\}$ up to homotopy. In Proposition
\ref{prop1a} the family of Riemann surfaces depends also on a second parameter $\alpha$, and the growth rate is determined in $\alpha$ and $\sigma$.

The proof of both propositions uses solutions of a $\overline{\partial}$-problem. The solution in the case of Proposition \ref{prop1a} uses a simple explicit formula.

\medskip
The author is grateful to B.Farb who suggested to use the concept of conformal module and extremal length for a proof of finiteness theorems, and to B.Berndtsson for proposing the kernel for solving the $\bar{\partial}$-problem that arises in the proof of Proposition \ref{prop1a}. The work on the paper was started while the author was visiting the Max-Planck-Institute and was finished during a stay at IHES. The author would like to thank these institutions for the support. The author is also indebted to Fanny Dufour for drawing the figures and to an anonymous referee whose critics helped to improve the overall quality of the paper.

\section{preliminaries on mappings, coverings, and extremal length}
%{Holomorphic mappings into the twice punctured plane}
\label{sec:1a}
In this section we will prepare the proofs of the Theorems.
% \ref{thm1}.
%We first need some preparation.

\noindent {\bf The change of the base point.}
Let $\mathcal{X}$ be a connected smooth open surface, and let $\alpha$ be an arc in $\mathcal{X}$ with initial point $x_0$ and terminating point $x$. Change the base point $x_0 \in \mathcal{X}$
along a curve $\alpha$ to the point $x \in \mathcal{X}$. This leads to an isomorphism $\mbox{Is}_{\alpha}: \pi_1(\mathcal{X},x_0) \to \pi_1(\mathcal{X},x)$ of fundamental groups induced by the correspondence $\gamma \to \alpha^{-1} \gamma \alpha$ for any loop $\gamma$ with base point $x_0$ and the arc $\alpha$ with initial point $x_0$ and terminating point $x$. We will denote the correspondence $\gamma \to \alpha^{-1} \gamma \alpha$ between curves also by $\mbox{Is}_{\alpha}$.

We call two homomorphisms $h_j:G_1 \to G_2, \, j=1,2,\,$ from a group $G_1$ to a group $G_2$ conjugate if there is an element $g' \in G_2$ such that for each $g \in G_1$ the equality $h_2(g)= {g'}^{-1} h_1(g) g'$ holds.
For two
arcs $\alpha_1$ and $\alpha_2$ with initial point $x_0$ and terminating point $x$ we have
$\alpha_2^{-1} \gamma \alpha_2= (\alpha_1^{-1}\alpha_2)^{-1} \alpha_1^{-1} \gamma \alpha_1 (\alpha_1^{-1}\alpha_2)$.
Hence, the two isomorphisms $\mbox{Is}_{\alpha_1}$ and $\mbox{Is}_{\alpha_2}$ differ by conjugation with the element of $\pi_1(\mathcal{X},x)$ represented by $\alpha_1^{-1}\alpha_2$.

Free homotopic curves are related by homotopy with fixed base point and an application of a homomorphism $\mbox{Is}_{\alpha}$
that is defined up to conjugation. Hence, free homotopy classes of curves can be identified with conjugacy classes of elements of the fundamental group  $\pi_1(\mathcal{X},x_0)$ of $\mathcal{X}$.

For two smooth manifolds $\mathcal{X}$ and $\mathcal{Y}$ with base points $x_0 \in \mathcal{X}$ and $y_0 \in \mathcal{Y}$ and a continuous mapping $F:\mathcal{X} \to \mathcal{Y}$ with $F(x_0)=y_0$ we denote by $F_*: \pi_1(\mathcal{X},x_0) \to \pi_1(\mathcal{Y},y_0) $ the induced map on fundamental groups. For each element $e_0\in \pi_1(\mathcal{X},x_0)$ the image $F_*(e_0)$ is called the monodromy along $e_0$, and the homomorphism $F_*$ is called the monodromy homomorphism corresponding to $F$.
The homomorphism $F_*$ determines the homotopy class of $F$ with fixed base point in the source and fixed value at the base point.
%in the target.
Consider a free homotopy $F_t, \, t \in (0,1)$, of homeomorphisms from $\mathcal{X}$ to $\mathcal{Y}$ such that the value $F_t(x_0)$ at the base
point $x_0$ of the source space varies along a loop. Then the homomorphisms $(F_0)_*$ and $(F_1)_*$ are related by conjugation with the element of the fundamental group of $\mathcal{Y}$ represented by the loop.

Using deformation retractions we see that each homomorphism $h: \pi_1(\mathcal{X},x_0) \to \pi_1(\mathcal{Y},y_0) $ equals $F_*$ for a continuous mapping $F:X\to Y$. Moreover,
if two homomorphisms $h_j: \pi_1(\mathcal{X},x_0) \to \pi_1(\mathcal{Y},y_0),\, j=0,1, $ are related by conjugation, $h_1=e^{-1}h_2 e$ for an element $e\in  \pi_1(\mathcal{Y},y_0)$, then  there is a free homotopy $F_t$  of mappings $X\to Y$ such that
$F_t(x_0)$ changes along a loop representing $e$ and
$(F_0)_*=h_0$, $(F_1)_*=h_1$.
Further, since the fundamental group $\pi_1(\mathcal{Y},y)$ with base point $y$ is related to the fundamental  group $\pi_1(\mathcal{Y},y_0)$ with base point $y_0$ by an isomorphism determined up to conjugation
we obtain the following theorem (see \cite{Ha},\cite{St}).
\smallskip

\noindent {\bf Theorem E.} {\it The free homotopy classes of continuous mappings from $\mathcal{X}$ to $\mathcal{Y}$
are in one-to-one correspondence to the set of conjugacy classes of homomorphisms between the fundamental groups of $\mathcal{X}$ and $\mathcal{Y}$.}
\smallskip

\noindent {\bf Extremal length.} The fundamental group $\pi_1  \stackrel{def}{=} \pi_1(\mathbb
{C}\setminus \{-1,1\},0)$ is canonically isomorphic to the fundamental group $ \pi_1(\mathbb{C}\setminus \{-1,1\},q')$ for an arbitrary point $q' \in (-1,1)$. For the arc $\alpha$ defining the isomorphism we take the unique arc contained in $(-1,1)$ that joins $0$ and $q'$. The fundamental group $\pi_1(\mathbb{C}\setminus \{-1,1\},0)$ is a free group in two generators. We choose standard generators $a_1$ and $a_2$, where $a_1$ is represented by a simple closed curve with base point $0$ which surrounds $-1$ counterclockwise, and $a_2$ is represented by a simple closed curve with base point $0$ which surrounds $1$ counterclockwise. For $q' \in (-1,1)$
we also denote by $a_j$ the generator of
$\pi_1(\mathbb{C}\setminus \{-1,1\},q')$ which is obtained from the respective standard generator of $\pi_1(\mathbb{C}\setminus\{-1,1\}, 0)$ by the standard isomorphism between fundamental groups with base point on $(-1,1)$.
More detailed, $a_1$ is the generator of $\pi_1(\mathbb{C}\setminus \{-1,1\},q')$ which is represented by a loop with base point $q'$ that surrounds $-1$ counterclockwise, and $a_2$ is the generator of $\pi_1(\mathbb{C}\setminus \{-1,1\},q')$ which is represented by a loop with base point $q'$ that surrounds $1$ counterclockwise. We refer to $a_1$ and $a_2$ as to the standard generators of $\pi_1(\mathbb{C}\setminus \{-1,1\},q')$.

Further, the group $ \pi_1(\mathbb {C}\setminus \{-1,1\},0)$ is canonically
isomorphic to the relative fundamental group $\pi_1^{tr}(\mathbb {C}\setminus \{-1,1\})
\stackrel{def}{=}\pi_1(\mathbb {C}\setminus \{-1,1\},(-1,1))$
whose elements are homotopy classes of (not necessarily closed) curves in $\mathbb{C}\setminus \{-1,1\}$ with end points on the interval $(-1,1)$.
We refer to $\pi_1^{tr}(\mathbb {C}\setminus \{-1,1\}) $ as fundamental group with totally real horizontal
boundary values ($tr$-boundary values for short). For an element
$w \in \pi_1(\mathbb{C}\setminus \{-1,1\}, q')$ with base point $q' \in (-1,1)$ we denote by $w_{tr}$ the element of the relative fundamental group $\pi_1^{tr}(\mathbb{C}\setminus \{-1,1\})$ with totally real boundary values, corresponding to $w$. For more details see \cite{Jo2}.

Each element of a free group can be written uniquely as a reduced word in the generators. (A word is reduced if neighbouring terms are powers of different generators.)
The degree (or word length) $d(w)$ of a reduced word $w$ in the generators of a free group is the sum of the absolute values of the powers of generators in the reduced word. If the word is the identity its degree is defined to be zero. We will identify elements of a free group with reduced words in generators of the group.

For a rectangle $R$  let
$f:R \to \mathbb{C} \setminus \{-1,1\}$ be a mapping which
admits a continuous extension to the closure $\bar R$ (denoted again by $f$) which maps the (open) horizontal sides to $(-1,1)$.
We say that the mapping $f$ represents an element $w_{tr} \in \pi_1^{tr}(\mathbb{C} \setminus \{-1,1\})$
if for each maximal vertical line segment contained in $R$
(i.e. $R$
intersected with a vertical line in $\mathbb{C}$) the
restriction of
$f$ to the closure of the line segment represents $w_{tr}$.

The extremal length $\;\;\Lambda(w_{tr})\;\;$ of an element $\;\;w_{tr}\;\;$ in the relative fundamental group $\;\pi_1^{tr}(\mathbb{C}\setminus \{-1,1\})\;$ is defined as
\begin{align}\label{eq1}
\Lambda(w_{tr})\stackrel{def}{=}& \inf \{\lambda(R): R\, \mbox{ a rectangle
which
admits a holomorphic map to} \;\mathbb{C}\setminus \{-1,1\} \nonumber \\
&  \,\mbox{ that
represents}\; w_{tr}\}\,.
\end{align}
For an element $w\in \pi_1(\mathbb{C}\setminus \{-1,1\},q')$ and the associated element $w_{tr}$ we will also write $\Lambda_{tr}(w)$ instead of $\Lambda(w_{tr})$.

Any reduced word $\;w\;$ in  $\;\pi_1(\mathbb{C}\setminus \{-1,1\},q')\;$ can be uniquely decomposed into syllables. They are defined as follows.
Each term $a_{j_i}^{k_i}$ with $|k_i|\geq 2$ is a syllable, and any maximal sequence of consecutive terms  $a_{j_i}^{k_i}$,  for which $|k_i|=1$ and all $k_i$ have the same sign, is a syllable (see \cite{Jo2}, \cite{Jo3}).
Let $d_k$ be the degree of the $k$-th syllable from the left. (We consider each syllable as a reduced word in the elements of the fundamental group.)
Put
\begin{equation}\label{eq3+}
\mathcal{L}_-(w)\stackrel{def}=   \sum \log(3 d_k), \quad \mathcal{L}_+(w)\stackrel{def}=   \sum \log(4 d_k)\,,
\end{equation}
where the sum runs over the degrees of all syllables of $w_{tr}$.
Notice that  $\mathcal{L}_{\pm}(w^{-1})=\mathcal{L}_{\pm}(w)$.  We define $\mathcal{L}_-({\rm Id})=\mathcal{L}_+({\rm Id})=0$ for the identity ${\rm Id}$.
We need the following theorem which is proved in \cite{Jo2} (see Theorem 1 there).

\medskip

\noindent {\bf Theorem F.} {\it If $w \in \pi_1(\mathbb{C}\setminus\{-1,1\},0)$ is not equal to a (trivial or non-trivial) power of $a_1$ or of $a_2$
then
\begin{equation}\label{eq3}
\frac{1}{ 2 \pi} \mathcal{L}_-(w)\leq \Lambda(w_{tr}) \leq 300 \mathcal{L}_+(w)\,.
 \end{equation}
}

\medskip

\noindent {\bf Regular zero sets.}
We will call a subset of a smooth manifold $\mathcal{X}$ a simple relatively closed curve
if it is the connected component of a regular level set of a smooth real-valued function on $\mathcal{X}$.

Let $\mathcal{X}$ be a connected finite open Riemann surface. Suppose the zero set $L$ of a non-constant smooth real valued function on $\mathcal{X}$ is regular.
Each component of $L$ is either a simple closed curve or
it can be parameterized by a continuous mapping $\ell:(-\infty,\infty)\to \mathcal{X}$. We call a component of the latter kind a simple relatively closed arc in $\mathcal{X}$.

A relatively closed curve $\gamma$ in a connected finite open Riemann surface $\mathcal{X}$
is said to be contractible to a hole of $\mathcal{X}$, if the following holds. Consider $\mathcal{X}$ as domain $\mathcal{X}^c\setminus \cup\mathcal{C}_j$ on a closed Riemann surface $\mathcal{X}^c$. Here the $\mathcal{C}_j$ are the holes, each is either a closed topological disc  with smooth boundary or a point. The condition is the following. For each pair  $U_1$, $U_2$ of open subsets of $ \mathcal{X}^c$, $\cup \mathcal{C}_j\subset U_1\Subset U_2$, there exists a homotopy of $\gamma$ that fixes $\gamma\cap U_1$ and moves $\gamma$ into $U_2$.
Taking for $U_2$ small enough neighbourhoods of  $\cup \mathcal{C}_j$ we see that the homotopy moves $\gamma$ into an annulus adjacent to one of the holes.

For each relatively compact domain $\mathcal{X}'\Subset \mathcal{X}$ in $\mathcal{X}$ there is a finite cover of $L\cap \overline{\mathcal{X}'}$ by open subsets $U_k$ of $\mathcal{X}$ such that each $L\cap U_k$ is connected. Each set $L\cap U_k$ is contained in a component of $L$. Hence, only finitely many connected components of $L$ intersect $\mathcal{X}'$.
Let $L_0$ be a connected component of $L$ which is a simple relatively closed arc parameterized by $\ell_0:\mathbb{R}\to \mathcal{X}$. Since each set $L_0\cap U_k$ is connected it is the image of an interval under $\ell_0$. Take real numbers $t_0^-$ and $t_0^+$ such that all these intervals are contained in $(t_0^-,t_0^+)$. Then the images $\ell\big((-\infty,t_0^-)\big)$ and
$\ell\big((t_0^+ ,+\infty)\big)$ are contained in $\mathcal{X}\setminus \mathcal{X}'$, maybe, in different components. Such parameters $t_0^-$ and $t_0^+$ can be found for each relatively compact deformation retract $\mathcal{X}'$ of $\mathcal{X}$. Hence for each relatively closed arc $L_0\subset L$ the set of limit points $L_0^+$ of $\ell_0(t)$ for $t\to \infty$ is contained in a boundary component of $\mathcal{X}$. Also, the set of limit points $L_0^-$ of $\ell_0(t)$ for $t\to -\infty$ is contained in a boundary component of $\mathcal{X}$. The boundary components may be equal or different.

Moreover, if
$\mathcal{X}'\Subset \mathcal{X}$ is a relatively compact domain in $\mathcal{X}$ which is a deformation retract of $\mathcal{X}$, and a connected component $L_0$ of $L$ does not intersect $\overline{\mathcal{X}'}$ then $L_0$ is contractible to a hole of $\mathcal{X}$. Indeed, $\mathcal{X}\setminus \overline{\mathcal{X}'}$ is the union of disjoint annuli, each of which is adjacent to a boundary component of $\mathcal{X}$, and the connected set $L_0$ must be
%the union of the limit sets $L_0^+ \cup L_0^-$
contained in a single annulus.

Further, denote by $L'$ the union of all connected components of $L$ that are simple relatively closed arcs. Consider those components $L_j$ of $L'$ that intersect $\mathcal{X}'$. There are finitely many such $L_j$. Parameterize each $L_j$ by a mapping $\ell_j:\mathbb{R}\to \mathcal{X}$. For each $j$ we let $[t_j^-,t_j^+]$ be a compact interval for which
\begin{equation}\label{eq2d}
\ell_j(\mathbb{R}\setminus [t_j^-,t_j^+]) \subset \mathcal{X}\setminus \overline{\mathcal{X}'}\,.
\end{equation}
Let $\mathcal{X}''$, $\mathcal{X}'\Subset \mathcal{X}''\Subset \mathcal{X}$, be a domain which is a deformation retract of $\mathcal{X}$ such that  $\ell_j([t_j^-,t_j^+])\subset \mathcal{X}''$ for each $j$. Then all connected components of $L'\cap \mathcal{X}''$, that do not contain a set $\ell_j([t_j^-,t_j^+])$, are contractible to a hole of $\mathcal{X}''$. Indeed, each such component is contained in the union of annuli $\mathcal{X}''\setminus \overline{\mathcal{X}'}$.

\smallskip

\noindent {\bf Some remarks on coverings.}
By a covering $P:\mathcal{Y} \to \mathcal{X}$ we mean a continuous map $P$ from a topological space $\mathcal{Y}$ to  a topological space $\mathcal{X}$ such that for each point $x \in \mathcal{X}$ there is a neighbourhood $V(x)$ of $x$ such that the mapping $P$ maps each connected component of the preimage of $V(x)$ homeomorphically onto $V(x)$. (Note that in function theory sometimes these objects are called unlimited unramified coverings to reserve the notion ''covering'' for more general objects.)

Let $X$ be a connected finite open Riemann surface with base point $q_0$ and let ${\sf P}: \tilde X \to X$
be the universal covering map.
Recall that a homeomorphism $\varphi:\tilde{X} \to \tilde{X}$ for which ${\sf P}\circ \varphi = {\sf P}$ is called a covering transformation (or deck transformation). The covering transformations form a group, denoted by ${\rm Deck}(\tilde{X},X)$.
For each pair of points $\tilde{x}_1, \tilde{x}_2 \in \tilde{X}$ with ${\sf P}(\tilde{x}_1)={\sf P}(\tilde{x}_2)$ there exists exactly one covering transformation that maps $\tilde{x}_1$ to $\tilde{x}_2$. (See e.g.  \cite{Fo}).

Throughout the paper we will fix a base point $q_0\in X$ and a base point $\tilde{q}_0\in {\sf P}^{-1}(q_0)\subset \tilde{X}$. The group of covering transformations of  $\tilde X$ can be identified with the fundamental group $\pi_1(X,q_0)$ of $X$ by the following correspondence. (See e.g. \cite{Fo}).

Take a covering transformation $\sigma \in {\rm Deck}(\tilde{X},X)$.
Let $\tilde{\gamma}_0$ be an arc in $\tilde X$ with initial point $\sigma(\tilde{q}_0)$
and terminating point $\tilde{q}_0$.
Denote by ${\rm Is}^{\tilde{q}_0}(\sigma)$ the element of  $\pi_1(X,q_0)$ represented by the loop ${\sf P}(\tilde{\gamma}_0)$. The mapping ${\rm Deck}(\tilde{X},X)\ni\sigma\to {\rm Is}^{\tilde{q}_0}(\sigma)\in \pi_1(X,q_0)$ is a group homomorphism. The homomorphism ${\rm Is}^{\tilde{q}_0}$ is injective and surjective, hence it is a group isomorphism. The inverse $({\rm Is}^{\tilde{q}_0})^{-1}$ is obtained as follows. Represent an element $e_0\in \pi_1(X,q_0)$ by a loop $\gamma_0$. Consider the lift $\tilde{\gamma}_0$ of $\gamma_0$ to $\tilde X$ that has terminating point $\tilde{q}_0$. Then
$({\rm Is}^{\tilde{q}_0})^{-1}(e_0)$ is the covering transformation that maps $q_0$ to the initial point of of $\tilde{\gamma}_0$.

For another point $\tilde{q}$ of $\tilde X$ and the point $q\stackrel{def}={\sf P}(\tilde{q})\in X$ the isomorphism ${\rm Is}^{\tilde{q}}:{\rm Deck}(\tilde{X},X)\to\pi_1(X,q)$ assigns to each $\sigma \in{\rm Deck}(\tilde{X},X)$ the element of $\pi_1(X,q)$ that is represented by ${\sf P}(\tilde{\gamma})$ for a curve $\tilde{\gamma}$ in $\tilde X$ that joins $\sigma(\tilde{q})$ with $\tilde q$.
%is defined in the same way as ${\rm Is}^{\tilde{q}_0}$
%was defined and
${\rm Is}^{\tilde{q}}$ is related to ${\rm Is}^{\tilde{q}_0}$ as follows.
Let $\tilde{\alpha}$ be an arc in $\tilde X$ with initial point $\tilde{q}_0$ and terminating point $\tilde{q}$.
Put $q={\sf P}(\tilde{q})$ and $\alpha={\sf P}(\tilde{\alpha})$.
Then for the isomorphism ${\rm Is}_{\alpha}:\pi_1(X,q_0)\to \pi_1(X,q)$ the equation
\begin{equation}\label{eq1''}
{\rm Is}^{\tilde{q}}(\sigma)=
{\rm Is}_{\alpha}\circ{\rm Is}^{\tilde{q}_0}(\sigma),\;\; \sigma \in {\rm Deck}(\tilde{X},X),
\end{equation}
holds, i.e. the following diagram Figure \ref{fig2} is commutative.
\begin{figure}[h]
\begin{center}
$$
\xymatrix{
&{\pi_1(X,q_0)} \ar[dd]^{{\rm Is}_{\alpha}    } \\
{{\rm Deck}(\tilde{X},X)  } \ar[ru]^{{\rm Is}^{\tilde{q}_0}} \ar[rd]_{{\rm Is}^{\tilde{q}}  } \\
& \pi_1(X,q)
}
$$
\end{center}
\caption{A commutative diagram related to the change of the base point}
\label{fig2}
\end{figure}
%{for covering transformations and the fundamental group with different base points}
Indeed, let
$\tilde{\alpha}^{-1}$ denote the curve that is obtained from
$\tilde{\alpha}$ by inverting the direction on $\tilde{\alpha}$, i.e. moving from $\tilde{q}$ to $\tilde{q}_0$. The curve $\sigma((\tilde{\alpha})^{-1})$ has initial point $\sigma(\tilde{q})$
and terminating point $\sigma(\tilde{q}_0)$.
Hence, for a curve $\tilde{\gamma}_0$ in $\tilde X$ that joins $\sigma(\tilde{q}_0)$ with $\tilde{q}_0$, the curve
$\sigma((\tilde{\alpha})^{-1}) \; \tilde{\gamma}_0\;\tilde{\alpha}$ in $\tilde X$ has initial point $\sigma(\tilde{q})$ and terminating point
$\tilde{q}$. Therefore ${\sf P}(\sigma((\tilde{\alpha})^{-1}) \; \tilde{\gamma}_0 \;\tilde{\alpha})$ represents ${\rm Is}^{\tilde{q}}(\sigma)$.
On the other hand
\begin{equation}\label{eq1'}
{\sf P}(\sigma(\tilde{\alpha}^{-1}) \; \tilde{\gamma}_0\;\tilde{\alpha})=
{\sf P}(\sigma(\tilde{\alpha}^{-1}))\;{\sf P}(\tilde{\gamma}_0)\; {\sf P}(\tilde{\alpha})=
\alpha^{-1} \gamma_0 \,\alpha
\end{equation}
represents ${\rm Is}_{\alpha}(e_0)$ with $e_0={\rm Is}^{\tilde q_0}(\sigma)$.
In particular, if $\tilde{q}'_0\in {\sf P}^{-1}(q_0)$ is another preimage
of the base point $q_0$ under the projection $\sf P$,
%by another preimage $\tilde{q}'_0\in {\sf P}^{-1}(q_0)$.
then the associated isomorphisms to the fundamental group $\pi_1(X,q_0)$ are conjugate, i.e. ${\rm Is}^{\tilde{q}'_0}(e_0)=(e'_0)^{-1} {\rm Is}^{\tilde{q}_0}(e_0) e'_0$ for each $e_0\in\pi_1(X,q_0)$. The element $e'_0$ is represented by the projection of an arc in $\tilde X$ with initial point $\tilde{q}_0$ and terminating point $\tilde{q}'_0$.

Keeping fixed $\tilde{q}_0$ and $q_0$ we will say that a point $\tilde{q}\in \tilde X$ and a curve $\alpha$ in $X$ are compatible if the diagram Figure \ref{fig2} is commutative, equivalently, if equation \eqref{eq1''} holds.
We may also start with choosing a curve $\alpha$ in $X$ with initial point $q_0$ and terminating point $q$.
Then there is a point $\tilde{q}=\tilde{q}(\alpha)$, such that $\tilde{q}$ and $\alpha$
are compatible.
Indeed,
let $\tilde{\alpha}$ be the lift of $\alpha$, that has initial point $\tilde{q}_0$.
Denote the terminating point of $\tilde{\alpha}$ by $\tilde{q}(\alpha)$, and repeat the previous arguments.

\smallskip

Let $N$ be a subgroup of $\pi_1(X,q_0)$. Denote by $X(N)$ the quotient $\tilde X \diagup ({\rm Is}^{\tilde{q}_0})^{-1}(N)$. We obtain a covering
$\omega^N_{{\rm Id}} :\tilde X \to X(N)$ with group of covering transformations isomorphic to $N$.
The fundamental group of $X(N)$ with base point $(q_0)_N\stackrel{def}= \omega_{\rm Id}^N(\tilde{q}_0)$ can be identified with $N$.

\begin{figure}[h]
\begin{center}
\includegraphics[width=7cm]{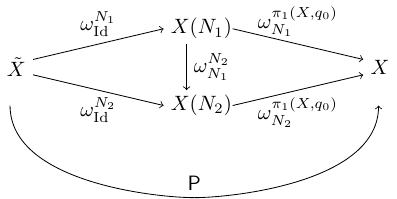}
\end{center}
\caption{A commutative diagram related to subgroups of the group of covering transformations}\label{fig3}
\end{figure}

If $N_1$ and $N_2$ are subgroups of $\pi_1(X,q_0)$ and $N_1$ is a subgroup of $N_2$ (we write $N_1 \leq N_2$), then there is a covering map $\omega^{N_2}_{N_1} :\tilde X \diagup  ({\rm Is}^{\tilde{q}_0})^{-1}(N_1) \to \tilde X \diagup  ({\rm Is}^{\tilde{q}_0})^{-1}(N_2)$, such that $\omega^{N_2}_{N_1} \circ \omega^{N_1}_{\rm Id}=\omega^{N_2}_{\rm Id}$. Moreover,
the diagram Figure \ref{fig3} below is commutative.

\medskip

Indeed, take any point $x_1 \in \tilde X \diagup  ({\rm Is}^{\tilde{q}_0})^{-1}(N_1)$
and a preimage $\tilde x$ of $x_1$ under $\omega^{N_1}_ {\rm Id} $.
There exists
a neighbourhood $V(\tilde x)$ of $\tilde x$ in $\tilde X$ such that $V(\tilde x)\cap \sigma(V(\tilde x))=\emptyset$ for all covering transformations $\sigma\in {\rm Deck}(\tilde{X}, X)$.
Then for $j=1,2$ the mapping $\omega^{N_j,\tilde x} _{\rm Id} \stackrel{def}{=} \omega^{N_j}_{\rm Id}  \mid V(\tilde x)$ is a homeomorphism from $V(\tilde x)$ onto its image denoted by $V_j$. Put $x_2=\omega^{N_2,\tilde x} _{\rm Id}(\tilde{x})$. The set $V_j\subset \tilde{X}\diagup ({\rm Is}^{\tilde{q}_0})^{-1}(N_j)  $ is a neighbourhood of $x_j$ for $j=1,2$.

For each preimage $\tilde{x}' \in (\omega^{N_1}_{\rm Id})^{-1}(x_1)$
there is a covering transformation $\varphi_{\tilde x,\tilde{x}'}$ in $({\rm Is}^{\tilde{q}_0})^{-1}(N_1)$ which maps a neighbourhood $V(\tilde {x}')$ of $\tilde{x}'$ conformally onto the neighbourhood $V(\tilde {x})$ of $\tilde{x}$ so that
on $V(\tilde{x}')$ the equality
$\omega^{N_1,\tilde{x}'} _{\rm Id}  =  \omega^{N_1,\tilde x}_{\rm Id} \circ \varphi _{\tilde x,\tilde{x}'}$ holds. Choose $\tilde{x} \in (\omega^{N_1} _{\rm Id})^{-1}(x_1)$ and define
\begin{equation}\label{eq1a'}
\omega^{N_2}_{N_1}(y) = \omega^{N_2, \tilde x}_{N_1}(y) \stackrel{def}= \omega^{N_2, \tilde x} _{\rm Id} ((\omega^{N_1,\tilde x} _{\rm Id} )^{-1}(y)) {\mbox{ \;for each\;}} y \in V_1\,.
\end{equation}
We get a correctly defined mapping from $V_1$ onto $V_2$. Indeed, since $N_1$ is a subgroup of $N_2$, the covering transformation $\varphi_{\tilde x,\tilde{x}'}$ is contained in $({\rm Is}^{\tilde{q}_0})^{-1}(N_2)$, and we get for another point $\tilde{x}' \in (\omega^{N_1}_ {\rm Id})^{-1}(x_1)$ the equality  $\omega^{N_2, \tilde{ x}'} _{\rm Id}  =  \omega^{N_2, \tilde x} _{\rm Id} \circ \varphi _{\tilde x,\tilde{x}'}$. Hence, for $y \in V_1(x_1)$
\begin{equation}\label{eq1b'}
\omega^{N_2, \tilde{ x}'} _{\rm Id} \circ(\omega^{N_1,\tilde{x}'} _{\rm Id})^{-1} (y)  = (\omega^{N_2, \tilde x} _{\rm Id} \circ \varphi_{\tilde x,\tilde{x}'})\circ (\omega^{N_1,\tilde x} _{\rm Id} \circ \varphi_{\tilde x,\tilde{x}'})^{-1}(y)=
\omega^{N_2,\tilde x} _{\rm Id} \circ(\omega^{N_1,\tilde x} _{\rm Id} )^{-1}(y).
\end{equation}
Since each mapping $\omega^{N_j,\tilde{x}}_ {\rm Id}$, $j=1,2,$ is a homeomorphism from
$V(\tilde{x})$ onto its image, the mapping $\omega^{N_2} _{N_1}$ is a homeomorphism
from $V(x_1)$ onto $V(x_2)$. The same holds for all preimages of
$V(x_2)$ under $\omega^{N_2} _{N_1}$. Hence, $\omega^{N_2} _{N_1}$
%that satisfies the condition for being
is a covering map. The commutativity of the part of the diagram
that involves the mappings $\omega^{N_1} _{\rm Id} $, $\omega^{N_2} _{\rm Id} $, and $\omega_{N_1}^{N_2}$ follows from equation \eqref{eq1a'}.
%is clear by construction.

The existence of $\omega_{N_1}^{\pi_1(X,q_0)}$ and the equality ${\sf P}= \omega_{N_1}^{\pi_1(X,q_0)} \circ \omega^{N_1} _{\rm Id}$ follows by applying the above arguments with $N_2= \pi_1(X,q_0)$. The equality ${\sf P}= \omega_{N_2}^{\pi_1(X,q_0)} \circ \omega^{N_2} _{\rm Id}$ follows in the same way. Since
\begin{align}
{\sf P} =  &\omega_{N_2}^{\pi_1(X,q_0)} \circ \omega_{N_1}^{N_2}\circ  \omega^{N_1} _{\rm Id} \nonumber \\
=  &\omega_{N_1}^{\pi_1(X,q_0)}\quad\quad\;\;\,  \circ \omega^{N_1} _{\rm Id},\nonumber
\end{align}
we have
\begin{equation}\nonumber
\omega_{N_2}^{\pi_1(X,q_0)} \circ \omega_{N_1}^{N_2} = \omega_{N_1}^{\pi_1(X,q_0)}
\end{equation}

We will also use the notation  $\omega^N \stackrel{def}= \omega^N _{\rm Id}$ and $\omega_N \stackrel{def}= \omega_N ^{\pi_1(X,q_0)}$ for a subgroup $N$ of $\pi_1(X,q_0)$.

Let again $N_1 \leq N_2$ be subgroups of $\pi_1(X,q_0)$.
Consider the covering $\omega_{N_1}^{N_2} :\tilde X \diagup ({\rm Is}^{{\tilde q}_0})^{-1}(N_1) \to \tilde X \diagup ({\rm Is}^{{\tilde q}_0})^{-1}(N_2)$. Let $\beta$ be a simple relatively closed curve in  $\tilde X \diagup ({\rm Is}^{\tilde q_0})^{-1}(N_2)$.
Then $(\omega_{N_1}^{N_2})^{-1}(\beta)$ is the union of simple relatively closed curves in $\tilde X \diagup ({\rm Is}^{{\tilde q}_0})^{-1}( N_1)$ and $\omega_{N_1}^{N_2}: (\omega_{N_1}^{N_2})^{-1}(\beta) \to \beta$ is a covering.
Indeed, we cover $\beta$ by small discs $U_k$ in $\tilde X \diagup ({\rm Is}^{{\tilde q}_0})^{-1}( N_1)$ such that for each $k$ the restriction of $\omega_{N_1}^{N_2}$
to each connected component of $(\omega_{N_1}^{N_2})^{-1}(U_k)$ is a homeomorphism onto $U_k$, and $U_k$ intersects $\beta$ %$(\omega^{\pi_1(X,q_0)}_{\mathcal{E}''})^{-1}(L)$
along a connected set. Take any $k$ with $U_k \cap \beta \neq \emptyset$. Consider the preimages $(\omega_{N_1}^{N_2})^{-1}(U_k)$.
Restrict  $(\omega_{N_1}^{N_2})$  to the intersection of each preimage $(\omega_{N_1}^{N_2})^{-1}(U_k)$ with $(\omega_{N_1}^{N_2})^{-1}(\beta)$. We obtain a homeomorphism onto $U_k \cap \beta$.
It follows that the map $(\omega_{N_1}^{N_2})$ is a covering from each connected component of $(\omega_{N_1}^{N_2})^{-1}(\beta)$ onto $\beta$.

\noindent{\bf The extremal length of monodromies.}
Let as before $X$ be a connected finite open Riemann surface with base point $q_0$, and $\tilde{q}_0$ a point in the universal covering $\tilde X$ for which ${\sf P}(\tilde{q}_0)=q_0$ for the covering map ${\sf P}:\tilde{X}\to X$.

Recall that for an arbitrary point $q\in X$ the free homotopy class of an element $e$ of the fundamental group $\pi_1(X,q)$ can be identified with the conjugacy class of elements of $\pi_1(X,q)$ containing $e$ and is denoted by $\widehat{e}$. Notice that for $e_0\in \pi_1(X,q_0)$ and a curve $\alpha$ in $X$ with initial point $q_0$ and terminating point $q$ the free homotopy classes of $e_0$ and of $e={\rm Is}_{\alpha}(e_0)$ coincide, i.e. $\widehat{e}=\widehat{e}_0$.
Consider a simple smooth relatively
closed curve $L$ in $X$.  We will say that a free homotopy class of curves $\widehat {e_0}$ intersects $L$ if each representative of $\widehat{e_0} $ intersects $L$. Choose an orientation of $L$. The intersection number of $\widehat {e_0}$ with the oriented curve $L$ is the intersection number with $L$ of some (and, hence, of each) smooth loop representing $\widehat {e_0}$ that intersects $L$ transversally.
This intersection number is the sum of the intersection numbers over all intersection points. The intersection number at an intersection point equals  $+ 1$ if the orientation determined by the tangent vector to the curve representing $\widehat {e_0}$
followed by the tangent vector to $L$ is the orientation of $X$, and equals $-1$ otherwise.

Let $A$ be an annulus equipped with an orientation (called positive orientation) of simple closed dividing curves in $A$.
% $\gamma$. %that is oriented. %in $A$ fixed.
(A relatively closed curve in a surface $X$ is called dividing, if $X\setminus \gamma$ consists of two connected components.)
A continuous mapping $\omega:A\to X$ is said to represent a conjugacy class $\widehat{ e}$ of elements of the fundamental group $\pi_1(X,q)$ for a point $q\in X$,
if the composition
%for each simple closed curve $\gamma$ in $A$ that is homologous to a boundary %component of $A$
%the restriction of
$\omega\circ\gamma$ represents $\widehat{ e}$ for each positively oriented dividing curve $\gamma$ in $A$ .
% or the inverse $\widehat{ e^{-1}}$.

Let $A$ be an annulus with base point $p$ with a chosen positive orientation of simple closed dividing curves in $A$.
%equipped with a simple closed dividing curve $\gamma$ through $p$.
% and distinguished generator of $\pi_1(A,p)$ chosen,
Let $\omega$ be a continuous mapping from $A$ to a finite Riemann surface ${X}$ with base point $q$ such that $\omega(p)=q$. We write $\omega:(A,p) \to ({X},q)$. The mapping is said to represent the element $e$ of the fundamental group $ \pi_1({X},q)$
if $\omega\circ \gamma$ represents $e$ for some (and hence for each) positively oriented simple closed dividing curve $\gamma$ in $A$ with base point $q$.

We associate to each element $e_0 \in  \pi_1(X,q_0)$ of the free group $\pi_1(X,q_0)$ the annulus  $X(\langle e_0 \rangle)=\tilde{X}\diagup ({\rm Is}^{\tilde{q}_0})^{-1}(\langle e_0 \rangle)$ with base point $q_{\langle e_0\rangle}= \omega_{\rm Id}^{\langle e_0 \rangle}(\tilde{q}_0)$ and the covering map $\omega_{\langle e_0 \rangle} \stackrel{def}=\omega^{\pi_1(X,q_0)}_{\langle e_0 \rangle}:X(\langle e_0 \rangle) \to X$.
By the commutative diagram \ref{fig3} the equality  $\omega_{\langle e_0 \rangle }(q_{\langle e_0 \rangle})= q_0$ holds.
We choose the orientation of simple closed dividing curves in $X(\langle e_0 \rangle)=\tilde{X}\diagup ({\rm Is}^{\tilde{q}_0})^{-1}(\langle e_0 \rangle)$
so that for a curve $\tilde \gamma$ in $\tilde X$ with terminating point $\tilde{q}_0$ and initial point ${\rm Is}^{\tilde{q}_0}(e_0)$ the curve
$\gamma_{\langle e_0 \rangle}\stackrel{def}=\omega_{\langle e_0 \rangle}(\tilde{ \gamma})$ is positively oriented.
The locally conformal mapping
$\omega^{\pi_1(X,q_0)}_{\langle e_0 \rangle}:(X(\langle e_0 \rangle), q_{\langle e_0\rangle})\to (X,q_0)$ represents $e_0$.
This follows from the equality $\omega_{\langle e_0 \rangle}(\gamma_{\langle e_0 \rangle })=\omega_{\langle e_0 \rangle}(\omega^{\langle e_0 \rangle}(\tilde{\gamma}))={\sf P}(\tilde{\gamma})=\gamma$, since ${\sf P}(\tilde{\gamma})$ represents $e_0$.

Take a curve $\alpha$ in $X$ that joins $q_0$ and $q$, and a point $\tilde{q}=\tilde{q}({\alpha})\in \tilde X$ such that $\alpha$ and $\tilde{q}$ are compatible, i.e. ${\rm Is}^{\tilde{q}}={\rm Is}_{\alpha}\circ {\rm Is}^{\tilde{q}_0}$ (see equation \eqref{eq1''}).
Put $e={\rm Is}_{\alpha}(e_0)$.
By equation \eqref{eq1''} $({\rm Is}^{\tilde{q}})^{-1}(e)=({\rm Is}^{\tilde{q}_0})^{-1}(e_0)$, hence,
$\tilde{X}\diagup ({\rm Is}^{\tilde{q}})^{-1}(\langle e\rangle)=\tilde{X}\diagup ({\rm Is}^{\tilde{q}_0})^{-1}(\langle e_0\rangle)=X(\langle e_0\rangle)$.
The locally conformal mapping $\omega_{\langle e_0\rangle}: X(\langle e_0\rangle)\to X$ takes the point $q_{\langle e\rangle}\stackrel{def}=\omega^{\langle e_0 \rangle}(\tilde{q})$ to $q\in X$. Moreover,
$\omega_{\langle e_0\rangle}: (X(\langle e_0\rangle),q_{\langle e\rangle})  \to (X,q)$ represents $e\in\pi_1(X,q)$. This can be seen by repeating the previous arguments.

Let $\alpha$ be an arbitrary curve in $X$ joining $q_0$ with $q$, and $\tilde{q}\in {\sf P}^{-1}(q)$ be arbitrary (i.e. $\alpha$ and $\tilde{q}$ are not required to be compatible).
Let $e\in \pi_1(X,q)$.
Denote the projection $\tilde{X}\to\tilde{X}\diagup ({\rm Is}^{\tilde{q}})^{-1}(\langle e\rangle)$ by $\omega^{\langle e\rangle,\tilde{q}}$, and the projection $\tilde{X}\diagup ({\rm Is}^{\tilde{q}})^{-1}(\langle e\rangle)\to X$ by $\omega_{\langle e\rangle,\tilde{q}}$. Put $q_{\langle e\rangle,\tilde{q}}=\omega^{\langle e\rangle,\tilde{q}}(\tilde{q})$.
For any such choice we choose the orientation of simple closed dividing curves on $\tilde{X}\diagup ({\rm Is}^{\tilde{q}})^{-1}(\langle e\rangle)$ so that $\omega^{\langle e\rangle,\tilde{q}}$ maps
any curve $\tilde{\gamma}$ in $\tilde X$ with initial point $({\rm Is}^{\tilde{q}})^{-1}(\langle e\rangle)(\tilde{q})$ and terminal point $\tilde{q}$ to a positively oriented dividing curve. We will call it the standard orientation of dividing curves in $\tilde{X}\diagup ({\rm Is}^{\tilde{q}})^{-1}(\langle e\rangle)$.
The mapping $\omega_{\langle e\rangle,\tilde{q}}: \big(\tilde{X}\diagup ({\rm Is}^{\tilde{q}})^{-1}(\langle e\rangle),q_{\langle e\rangle,\tilde{q}} \big)\to (X,q)$
represents $e$.

Since the mapping $\;\;\omega_{\langle e_0 \rangle }:(X(\langle e_0 \rangle),(q_0)_{\langle e_0\rangle}) \to (X,q_0)\;\;$ represents $\;e_0\;$,
the mapping $\omega_{\langle e_0 \rangle }:X(\langle e_0 \rangle) \to X$ represents the free homotopy class $\widehat{e_0}$. The following simple lemma will be useful.
\begin{lemm}\label{lem0'}
The annulus $X(\langle e_0 \rangle)$ has
smallest extremal length among annuli which admit a holomorphic mapping to $X$, that represents the conjugacy class $\widehat{e_0}$.
\end{lemm}
\noindent In other words, $X(\langle e_0 \rangle)$ is the "thickest" annulus with the property stated in Lemma \ref{lem0'}.

\medskip

\noindent {\bf Proof}. Take an annulus $A$ with a choice of positive orientation of simple closed dividing curves. Suppose $A \overset{\omega}{-\!\!\!\longrightarrow}X$
is a holomorphic mapping that represents $\widehat{e_0}$. The annulus
$A$ is conformally equivalent to a round annulus in the plane, hence, we may assume that $A$ has the form $A=\{z\in \mathbb{C}: \, r<|z|<R\}$ for  $0\leq r<R\leq \infty$ and the positive orientation of dividing curves is the counterclockwise one.

Take a positively oriented simple closed dividing curve $\gamma^{A}$ in $A$.
% that is homologous to a boundary component of
%$A$. After maybe, reorienting $\gamma^{A}$,
Its image $\omega\circ \gamma^{A}$ under $\omega$ represents the class $\widehat{e_0}$. Choose a point $q^{A}$ in $\gamma^{A}$, and put $q=\omega(q^{A})$.
Then $\gamma^{A}$ represents a generator of $\pi_1(A,q^{A})$ and $\gamma=\omega\circ\gamma^{A} $ represents an element $e$ of $\pi_1(X,q)$ in the conjugacy class $\widehat{e_0}$.
Choose a curve $\alpha$ in $X$ with initial point $q_0$ and terminating point $q$,
and a point $\tilde q$ in $\tilde X$ so that  $\alpha$ and $\tilde q$ are compatible,
and, hence, for $e={\rm Is}_{\alpha}(e_0)$ the equality $({\rm Is}^{\tilde{q}_0})^{-1}(e_0)=({\rm Is}^{\tilde{q}})^{-1}(e)$ holds.
Let $L$ be the relatively closed arc $\{p\cdot r: r\in \mathbb{R}\} \cap A$ in $A$ that contains $p$.
After a homotopy of $\gamma^{A}$ with fixed base point, we may assume that its base point $q^{A}$ is the only point of $\gamma^{A}$ that is contained in $L$.
The restriction $\omega|(A\setminus L)$ lifts to a mapping $\tilde{\omega}:(A\setminus L)\to \tilde X$, that extends continuously to the two strands $L_{\pm}$ of $L$. (Here $L_-$ contains the initial point of $\gamma^{A}$.) Let $p_{\pm}$ be the copies of $p$ on the two strands $L_{\pm}$. We choose the lift $\tilde{\omega}$ so that $\tilde{\omega}(p_+)=\tilde{ q}$.
Since the mapping $(A,q^{A})\to (X,q)$ represents $e$, we obtain $\tilde{\omega}(p_-)=\sigma(\tilde{ q})$ for $\sigma= ({\rm Is}^{\tilde{q}})^{-1}(e)$.
Then for each $z\in L$ the covering transformation $\sigma$ maps the point $\tilde{z}_+\in\tilde{\omega}(L_+)$ for which ${\sf P}(\tilde{z}_+)=z$ to the point $\tilde{z}_-\in\tilde{\omega}(L_-)$ for which ${\sf P}(\tilde{z}_-)=z$. Hence $\omega$ lifts to
a holomorphic mapping $\iota: A\to X(\langle e_0 \rangle)$.
By Lemma 7 of \cite{Jo2} $\lambda(A) \geq \lambda( X(\langle e_0 \rangle))$. \hfill $\Box$

\medskip

For each point $q\in X$ and each element $e\in \pi_1(X,q)$ we denote by $A(\widehat {e})$ the conformal class of the ''thickest'' annulus that admits a holomorphic mapping into $X$ that represents $\widehat{e}$. We saw that $\lambda(A(\widehat{e_0}))=\lambda(\tilde{X}\diagup ({\rm Is}^{\tilde{q}_0})^{-1}(\langle e_0 \rangle))$ for $e_0\in \pi_1(X,q_0)$.
By the same reasoning as before $\lambda(A(\widehat{e}))=\lambda(\tilde{X}\diagup ({\rm Is}^{\tilde{q}'})^{-1}(\langle e \rangle))$ for each
$\tilde{q}' \in \tilde X$ and each element $e \in \pi_1(X,{\sf P}(\tilde{q}'))$.
Hence, if $e_0$ and $e$ are conjugate, then $\lambda(\tilde{X}\diagup ({\rm Is}^{\tilde{q}_0})^{-1}(\langle e_0 \rangle))=\lambda(\tilde{X}\diagup ({\rm Is}^{\tilde{q}})^{-1}(\langle e \rangle))$ for any $\tilde{q}_0\in {\sf P}^{-1}(q_0)\subset \tilde X$ and any $\tilde{q}\in {\sf P}^{-1}(q)$.
Notice that $A(\reallywidehat{e^{-1}})=A(\widehat e)$ for each $e\in \pi_1(X,q),\, q\in X$.

\section{Holomorphic mappings into the twice punctured plane}
\label{sec:2}

The following lemma will be crucial for the estimate of the $\mathcal{L}_-$-invariant of the monodromies of holomorphic mappings from a finite open Riemann surface to $\mathbb{C}\setminus \{-1,1\}$.

\begin{lemm}\label{lem2}
Let $f:X \to \mathbb{C}\setminus \{-1,1\}$ be a non-contractible holomorphic function
on a connected finite open Riemann surface $X$, such that $0$ is a regular value of ${\rm Im}f$. Assume that $L_0$ is a simple relatively closed curve in $X$ such that $f(L_0)\subset (-1,1)$. Let $q \in L_0$ and
%$e\in \pi_1(X,q)$, Put
$q'=f(q)$.

If for an element $e\in \pi_1(X,q)$
the free homotopy class $\widehat e$ intersects $L_0$,
then either the reduced word $f_*(e) \in \pi_1(\mathbb{C}\setminus \{-1,1\},q')$ is a  non-zero power of a standard generator of $\pi_1(\mathbb{C}\setminus \{-1,1\}, q')$ or the inequality
\begin{equation}\label{eq1a}
\mathcal{L}_-(f_*(e)) \leq 2\pi \lambda(A(\widehat{e}))
\end{equation}
holds.
\end{lemm}
Notice that we make a normalization in the statement of the Lemma by  requiring that $f$ maps $L_0$ into the interval $(-1,1)$, not merely into $\mathbb{R}\setminus \{-1,1\}$.

Lemma \ref{lem2}
will be a consequence of the following lemma.

\begin{lemm}\label{lem1}
Let $X$, $f$, $L_0$, $q\in L_0$ be as in Lemma \ref{lem2}, and $e\in \pi_1(X,q)$.
Let $\tilde q$ be an arbitrary point in ${\sf P}^{-1}(q)$.
Consider the annulus $A\stackrel{def}=\tilde{X}\diagup ({\rm Is}^{\tilde{q}})^{-1}(\langle e \rangle)$ and the holomorphic projection $\omega_{A}\stackrel{def}=\omega_{\langle e\rangle,\tilde{q}}$.
Put $q_{A}\stackrel{def}=\omega^{\langle e\rangle,\tilde{q}}(\tilde{q})$ and let $L_{A}\subset A$ be the connected component of $(\omega_{A})^{-1}(L_0)$ that contains $q_{A}$. Then the mapping $\omega_{A}:(A,q_{A})\to (X,q)$ represents $e$.

If $\widehat e$ intersects $L_0$, then $L_{A}$ is a relatively closed curve in $A$ that has limit points on both boundary components of $A$, and the lift $f\circ \omega_{A}$ is a holomorphic function
on $A$ that maps $L_{A}$ into $(-1,1)$.
\end{lemm}

\noindent{\bf Proof of Lemma \ref{lem1}.}
Let $\gamma:[0,1]\to X$ be a curve with base point $q$ in $X$ that represents $e$, and let $\tilde{\gamma}$ be the lift of $\gamma$ to $\tilde X$ with terminating point $\tilde{\gamma}(1)$ equal to $\tilde q$.
Put $\sigma\stackrel{def}= ({\rm Is}^{\tilde{q}})^{-1}(e)$.
Then the initial point
$\tilde{\gamma}(0)$ equals $\sigma(\tilde{q})$.

All connected components of ${\sf P}^{-1}(L_0)$ are relatively closed curves in $\tilde{X}\cong \mathbb{C}_+$ (where $\mathbb{C}_+$ denotes the upper half-plane) with limit points on the boundary of $\tilde X$. Indeed, the lift $f\circ {\sf P}$ of $f$ to $\tilde X$ takes values in $(-1,1)$ on ${\sf P}^{-1}(L_0)$. Hence, $|\exp(\pm \,i\, f\circ {\sf P})|=1$ on ${\sf P}^{-1}(L_0)$. A compact connected component of ${\sf P}^{-1}(L_0)$ would bound a relatively compact topological disc in $\tilde{X}=\mathbb{C}_+$, and by the maximum principle $|\exp(\pm \,i\, f\circ {\sf P})|=1$ on the disc. This would imply that $ f\circ {\sf P}$ is constant on $\tilde X$ in contrary to the assumptions.

%by the maximum principle applied to $\exp(f\circ {\sf P})$ the latter cannot contain %compact connected components since $f$ is not constant.

Let $\tilde{L}_{\tilde{q}}$ be the connected component of ${\sf P}^{-1}(L_0)$ that contains $\tilde{q}$.  The point $\sigma(\tilde{q})$ cannot be contained in $\tilde{L}_{\tilde{q}}$. Indeed, assume the contrary. Then the arc  $\tilde{\gamma '}$ on  $\tilde{L}_{\tilde{q}}$ joining $\sigma(\tilde{q})$ and $\tilde{q}$ is homotopic in $\tilde X$ with fixed endpoints to $\tilde{\gamma}$.
The projection $\gamma '={\sf P}(\tilde{\gamma '})$ is contained in $L_0$ and is homotopic in $X$ with fixed endpoints to $\gamma$. Since $\gamma$ represents $e$ and $e$ is a primitive element of the fundamental group $\pi_1(X,q)$, this is possible only if $L_0$ is compact (and after orienting it) $L_0$ represents $e$.
A small translation of $\gamma'$ to a side of $L_0$ gives a curve in $X$ that does not intersect $L_0$ and represents the free homotopy class $\widehat e$ of $e$.
This contradicts the fact that $\widehat{e}$ intersects $L_0$. Since $\sigma(L_0)$ is also a connected component of ${\sf P}^{-1}(L_0)$, the curves $\tilde{L}_{\tilde{q}}$ and $\sigma(\tilde{L}_{\tilde{q}})$ are disjoint.

Each of the two connected components $\tilde{L}_{\tilde{q}}$ and $\sigma(\tilde{L}_{\tilde{q}})$ divides $\tilde X$. Let $\Omega$ be the domain on $\tilde X$ that is bounded by $\tilde{L}_{\tilde{q}}$ and $\sigma(\tilde{L}_{\tilde{q}})$ and parts of the boundary of $\tilde X$. After a homotopy of $\tilde{\gamma}$  that fixes the endpoints we may assume that $\tilde{\gamma}((0,1))$ is contained in $\Omega$. Indeed, for each connected component of  $\tilde{\gamma}((0,1))\setminus \Omega$ there is a homotopy with fixed endpoints that moves the connected component to an arc on $\tilde{L}_{\tilde{q}}$ or  $\sigma(\tilde{L}_{\tilde{q}})$. A small perturbation yields a curve $\tilde{\gamma}'$ which is homotopic with fixed endpoints to $\tilde{\gamma}$ and has interior contained in $\Omega$. Notice that by the same reasoning as above, $\tilde{\gamma}'((0,1))$ does not meet any $\sigma^k(\tilde{L}_{\tilde{q}})$.

The curve $\omega^{\langle e\rangle,\tilde{q}}(\tilde{\gamma}')$ is a closed curve on $A$ that represents a generator of the fundamental group of $A$ with base point $q_A$. Moreover, $\omega_A\circ \omega^{\langle e\rangle,\tilde{q}}(\tilde{\gamma}')=  \omega_{\langle e\rangle,\tilde{q}}\circ \omega^{\langle e\rangle,\tilde{q}}(\tilde{\gamma}')={\sf P}(\tilde{\gamma}')$ represents $e$. Hence, the mapping $\omega_A:(A,q_A)\to (X,q)$ represents $e$.

The curve $\omega^{\langle e\rangle,\tilde{q}}(\tilde{\gamma}')$ intersects $L_{\langle e\rangle}=\omega^{\langle e\rangle,\tilde{q}}(\tilde{L}_{\tilde{q}})$ exactly once.
Hence, $L_{\langle e\rangle}$ has limit points on both boundary circles of $A$ for otherwise $L_{\langle e\rangle}$ would intersect one of the components of $A\setminus \omega^{\langle e\rangle,\tilde{q}}(\tilde{\gamma}')$ along a set which is relatively compact in $A$, and $\tilde{\gamma}'$ would have intersection number zero with         $L_{\langle e\rangle}$. It is clear that $f\circ \omega_A(L_A)=f(L_0)\subset (-1,1)$.
The lemma is proved. \hfill $\Box$

\medskip

\noindent {\bf Proof of Lemma \ref{lem2}.}
Let $\omega_A:(A,q_A)\to (X,q)$ be the holomorphic mapping from Lemma \ref{lem1} that represents $e$, and let $L_A\ni q_A$ be the relatively closed curve in $A$ with limit set on both boundary components of $A$. Consider a positively oriented dividing curve  ${\gamma}_A:[0,1]\to A$ with base point $\gamma(0)=\gamma(1)=q_A$ such that $\gamma_A((0,1))\subset A\setminus L_A$. The curve $\gamma=\omega_A(\gamma_A)$ represents $e$. The mapping $f\circ \omega_A$ is holomorphic on $A$ and $f\circ \omega_A(\gamma_A)=f(\gamma)$ represents $f_*(e)\in \pi_1(\mathbb{C}\setminus\{-1,1\}, q')$ with $q'=f\circ \omega_A(q_A)=f(q)\in (-1,1)$. Hence,
$f\circ \omega_A(\gamma_A)$ also represents the element $(f_*(e))_{tr}\in \pi_1^{tr}(\mathbb{C}\setminus\{-1,1\})$ in the relative fundamental group  $\pi_1(\mathbb{C}\setminus \{-1,1\},(-1,1))=\pi_1^{tr}(\mathbb{C}\setminus\{-1,1\})$ corresponding to $f_*(e)$.

We prove now that $\Lambda_{tr}(f_*(e))\leq \lambda(A)$.
Let $A_0\Subset A$ be any relatively compact annulus in $A$ with smooth boundary such that $q_A\in A_0$. If $A_0$ is sufficiently large, then the connected component $L_{A_0}$ of $L_A\cap A_0$ that contains
$q_A$ has endpoints on different boundary components of $A_0$. The set $A_0 \setminus L_{A_0}$ is a curvilinear rectangle.
The open horizontal curvilinear sides are the strands of the cut that are reachable from the curvilinear rectangle moving counterclockwise, or clockwise, respectively. The open vertical curvilinear sides are obtained from the boundary circles of $A_0$ by removing an endpoint of the arc $L_{A_0}$. Since $f\circ\omega_A$ maps $L_{A}$ to $(-1,1)$, the restriction of  $f\circ\omega_A$ to $A_0 \setminus L_{A_0}$ represents $(f_*(e))_{tr}$. Hence,
\begin{equation}\label{eq2a}
\Lambda_{tr}(f_*(e))\leq \lambda(A_0 \setminus L_{A_0})\,.
\end{equation}
Moreover,
\begin{equation}\label{eq2b}
\lambda(A_0 \setminus L_{A_0})\leq \lambda(A_0)\,.
\end{equation}
This is a consequence of the following facts. First,  $\lambda(A_0 \setminus L_{A_0})$ is equal to the extremal length $\lambda(\Gamma(A_0 \setminus L_{A_0}))$ in the sense of Ahlfors \cite{A1} of the family $\Gamma(A_0 \setminus L_{A_0})$
of curves in the curvilinear rectangle $A_0 \setminus L_{A_0}$ that join the two horizontal sides of the curvilinear rectangle. Further, $\lambda(A_0)$ is equal to the extremal length $\lambda(\Gamma(A_0)) $  \cite{A1} of the family $\Gamma(A_0)$  of curves in $A_0$ that are free homotopic to simple closed positively oriented dividing curves in $A_0$. Finally,
by \cite{A1}, Ch.1 Theorem 2, the  inequality
\begin{align}\label{eq2}
\lambda(\Gamma(A_0 \setminus L_{A_0})) \leq  \lambda(\Gamma(A_0))\,
\end{align}
holds.
We obtain the inequality
$\Lambda_{tr}(f_*(e))   \leq \lambda(A_0)\,$
for each annulus $A_0\Subset A$,
hence, since $A$ belongs to the class $A(\widehat{e})$ of conformally equivalent annuli,
\begin{align}\label{eq2c}
\Lambda_{tr}(f_*(e))   \leq \lambda(A(\widehat{e}))\,,
\end{align}
and the Lemma follows from Theorem F.
\hfill $\Box$

\medskip

\noindent {\bf The monodromies along two generators.}
In the following Lemma we combine the information on the monodromies along two generators of the fundamental group $\pi_1(X,q)$. We allow the situation when the monodromy along one generator or along both generators of the fundamental group of $X$ is a power of a standard generator of $\pi_1(\mathbb{C}\setminus \{-1,1\},f(q))$.

\begin{lemm}\label{lem3}
Let $f:X \to \mathbb{C} \setminus \{-1,1\}$ be a holomorphic function on a connected open Riemann surface $X$ such that $0$ is a regular value of the imaginary part of $f$.
Suppose $f$ maps a simple relatively closed curve $L_0$ in $X$ to $(-1,1)$, and $q$ is
a point in $L_0$. Let $e^{(1)}$ and $e^{(2)}$ be primitive elements of $\pi_1(X,q)$. Suppose that for each $e= e^{(1)},\;e= e^{(2)}$, and $e=e^{(1)}\,e^{(2)}$,
the free homotopy class $\widehat e$ intersects $L_0$.
Then either $f_*(e^{(j)}),\, j=1,2,\,$ are (trivial or non-trivial) powers of the same standard generator of $\pi_1(\mathbb{C} \setminus \{-1,1\},q')$ with $q'=f(q) \in (-1,1)$, or each of them is the product of at most two elements $w_1$ and $w_2$ of $\pi_1(\mathbb{C} \setminus \{-1,1\},q')$ with
\begin{equation}\nonumber
\mathcal{L}_-(w_j) \leq  2\pi \lambda_{e^{(1)},e^{(2)}},\, j=1,2,
\end{equation}
where
\begin{equation}\nonumber
\lambda_{e^{(1)},e^{(2)}} \stackrel{def}=\max\{\lambda(A(\reallywidehat{e^{(1)}})),\, \lambda(A(\reallywidehat{e^{(2)}})),\, \lambda(A(\reallywidehat{e^{(1)}\,e^{(2)}}))\}.
\end{equation}
Hence,
\begin{equation}\label{eq3e'}
\mathcal{L}_-(f_*(e^{(j)})) \leq 4 \pi \lambda_{e^{(1)},e^{(2)}}, \, j=1,2.
\end{equation}
\end{lemm}

\noindent {\bf Proof.}
If the monodromies $f_*(e^{(1)})$ and $f_*(e^{(2)})$ are not powers of a single standard generator (the identity is considered as zeroth power of a standard generator) we obtain the following.
At most two of the elements, $f_*(e^{(1)})$,  $f_*(e^{(2)})$, and $f_*(e^{(1)}\, e^{(2)})= f_*(e^{(1)})\, f_*(e^{(2)})$, are powers of a standard generator, and if two of them are powers of a standard generator, then they are non-zero powers of different standard generators. If two of them are non-zero powers of standard generators, then the third has the form
$a_{\ell }^{k} a_{\ell '}^{k'}$ with $a_{\ell }$ and $a_{\ell ' }$
being different generators and $k$ and $k'$ being non-zero integers.
By Lemma \ref{lem2} the $\mathcal{L}_-$ of the third element
does not exceed $2\pi \lambda_{e^{(1)},e^{(2)}}$. On the other hand it equals   $\log(3|k'|)+ \log(3|k|)$. Hence, $\mathcal{L}_-(a_{\ell }^{k}) = \log(3|k|)\leq 2\pi \lambda_{e^{(1)},e^{(2)}}$
and  $\mathcal{L}_-(a_{\ell' }^{k'}) = \log(3|k'|)\leq 2\pi \lambda_{e^{(1)},e^{(2)}}$.

If two of the elements $f_*(e^{(1)})$,  $f_*(e^{(2)})$, and $f_*(e^{(1)}\, e^{(2)})= f_*(e^{(1)})\, f_*(e^{(2)})$, are not powers of a standard generator, then the $\mathcal{L}_-$ of each of the two elements does not exceed $2\pi \lambda_{e^{(1)},e^{(2)}}$. Since the $\mathcal{L}_-$ of an element coincides with the $\mathcal{L}_-$ of its inverse, the third element is the product of two elements with $\mathcal{L}_-$ not exceeding $2\pi \lambda_{e^{(1)},e^{(2)}}$.
Since for $x,x'\geq 2$ the inequality $\log(x+x')\leq \log x +\log x'$ holds,
the $\mathcal{L}_-$ of the product does not exceed the sum of the $\mathcal{L}_-$ of the factors. Hence the $\mathcal{L}_-$ of the third element does not exceed $4\pi \lambda_{e^{(1)},e^{(2)}}$. Hence, inequality \eqref{eq3e'} holds.
\hfill $\Box$

\medskip

The following proposition states
the existence of suitable connected components of the zero set of the imaginary part of certain analytic functions on tori with a hole and on planar domains. For any subset $\mathcal{E}'$ of $\pi_1(X;q_0)$ we denote by $(\mathcal{E}')^{-1}$ the set of all elements that are inverse to elements in  $\mathcal{E}'$. Recall that $\mathcal{E}_j$ is the set of primitive elements of $ \pi_1(X,q_0)$ which can be written as product of at most $j$ elements of $\mathcal{E}\cup (\mathcal{E})^{-1} $ for the set $\mathcal{E}$ of generators of $ \pi_1(X,q_0)$ chosen in the introduction.

\begin{prop}\label{prop2a} Let $X$ be a torus with a hole or a planar domain with base point $q_0$ and fundamental group $\pi_1(X,q_0)$, and let $\mathcal{E}$  be a set of generators  of $\pi_1(X,q_0)$ that is associated to a standard bouquet of circles for $X$.
%was chosen in Section 1.
Let $f: X \to \mathbb{C}\setminus \{-1,1\}$ be a non-contractible holomorphic mapping
such that $0$ is a regular value of ${\rm{Im}} f$. Then there exist a simple relatively closed curve $L_0\subset X$ such that
$f(L_0) \subset \mathbb{R}\setminus \{-1,1\},$ and a set $\mathcal{E}_2'  \subset \mathcal{E}_2 \subset\pi_1(X,q_0)$ of primitive elements of $\pi_1(X,q_0)$,
such that the following holds.
Each element $e_{j,0} \in \mathcal{E} \subset \pi_1(X,q_0)$ is the product of at most two elements of $\mathcal{E}_2'\cup (\mathcal{E}_2')^{-1} $.  Moreover, for each $e_0 \in \pi_1(X,q_0)$ which is the product of one or two elements from $\mathcal{E}_2'$ the free homotopy class $\widehat{ e_0}$ has positive intersection number with $L_0$ (after suitable orientation of $L_0$).

If $X$ is a torus with a hole or $X$ equals $\mathbb{P}^1$ with three holes, we may chose $\mathcal{E}_2'$ consisting of two elements, one of them contained in $\mathcal{E}$, the other is either contained in $\mathcal{E}\cup \mathcal{E}^{-1}$ or is a product of two elements of $\mathcal{E}$.
\end{prop}
Notice the following facts. By Theorem E a mapping $f:X\to \mathbb{C}\setminus \{-1,1\}$ is contractible if and only if  for each $e_0\in \pi_1(X,q_0)$ the monodromy $f_*(e_0)$ is equal to the identity. The mapping $f$ is reducible if and only if the mondromy mapping $f_*:\pi_1(X,q_0)\to \pi_1(\mathbb{C}\setminus\{-1,1\},f(q_0))$ is conjugate to a mapping into a subgroup $\Gamma$ of $\pi_1(\mathbb{C}\setminus\{-1,1\},f(q_0))$ that is generated by a single element that is represented by a curve which separates one of the points $1,-1$ or $\infty$ from the other points. In other words, $\Gamma$ is (after identifying fundamental groups with different base point up to conjugacy) generated by a conjugate of one of the elements $a_1$, $a_2$ or $a_1a_2$ of $\pi_1(\mathbb{C}\setminus\{-1,1\},0)$.

If $f$ is irreducible, then it is not contractible, and, hence, the preimage $f^{-1}(\mathbb{R})$ is not empty.

Denote by $M_1$ a M\"obius transformation which permutes the points  $-1,\,1,\, \infty$ and maps the interval $(-\infty,-1)$ onto $(-1,1)$, and let $M_2$ be a M\"obius transformation which permutes the points $-1,\,1,\, \infty$ and maps the interval $(1,\infty)$ onto $(-1,1)$. Let $M_0\stackrel{def}= \mbox {Id}$.

The main step for the proof of Theorem \ref{thm1} is the following Proposition \ref{prop2}.

Recall that $\lambda_j(X)$ was defined in the introduction. Since for $e_0\in \pi_1(X,q_0)$
the equality $\lambda(\tilde{X}\diagup ({\rm Is}^{\tilde{q}_0})^{-1}(\langle e_0\rangle))=\lambda(A(\widehat{e_0}))$ holds, $\lambda_j(X)$ is the maximum of $\lambda(A(\widehat{e_0}))$ over $e_0\in \mathcal{E}_j$.

\begin{prop}\label{prop2} Let $X$ be a connected finite open Riemann surface with base point $q_0$, and let $\mathcal{E}$ be the set of generators of $\pi_1(X,q_0)$ that was chosen in Section 1.
Suppose
$f:X \to \mathbb{C}\setminus \{-1,1\}$ is an irreducible holomorphic mapping, such that
$0$ is a regular value of ${\rm{Im}}f$. Then for one of the functions $M_l \circ f,\, l=0,1,2,\,$
which we denote by $F$, there exists a point $q \in X$ (depending on $f$), such that
the point
$q' \stackrel{def}=F(q)$ is contained in  $(-1,1)$, and a curve $\alpha$ in $X$ joining $q_0$ with $q$, such that the following holds.
For each element $e_j\in\mbox{Is}_{\alpha}(\mathcal{E})$
the monodromy $F_*(e_j)$ is the product of at most four elements of $\pi_1(\mathbb{C}\setminus\{-1,1\},q')$ of $\mathcal{L}_-$ not exceeding $2 \pi \lambda_7(X)$ and, hence,
\begin{equation}\label{eq3b}
\mathcal{L}_-(F_*(e_j)) \leq 8\pi \lambda_7(X) \;\; \mbox{for each}\;\, j.
\end{equation}

If $X$ is a torus with a hole the proposition holds with $\lambda_7(X)$ replaced by $\lambda_3(X)$. If $X$ is a planar domain the proposition holds with $\lambda_4(X)$ instead of $\lambda_7(X)$.
\end{prop}
Notice, that
all monodromies of contractible mappings are equal to the identity, hence the inequality \eqref{eq3b} holds automatically for contractible mappings.

\noindent We postpone the proof of the two propositions and prove first the Theorem \ref{thm1}.

\smallskip

\noindent {\bf Proof of Theorem \ref{thm1}.}
Let $X$ be a connected finite open Riemann surface (possibly of second kind)
with base point $q_0$. Consider an arbitrary open Riemann surface $X^0 \Subset X$ which is relatively compact in $X$ and is a deformation retract of $X$.
%Using Theorem E and Proposition \ref{prop2} we want to bound the
%number of irreducible free homotoy classes of mappings from $X$ to %$\mathbb{C}\setminus \{-1,1\}$ that contain a holomorphic mapping.
Consider a free homotopy class of mappings from $X$ to $\mathbb{C}\setminus \{-1,1\}$, that is represented by an irreducible holomorphic mapping ${\sf f}:X\to \mathbb{C}\setminus \{-1,1\}$. Notice that the restriction ${\sf f}|X^0$ is also irreducible.
%Recall that irreducible mappings are not contractible.
Take a small enough positive number $\varepsilon$,
%be a small enough real number $\varepsilon$,
such that the function $({\sf f} - i \varepsilon)\mid X^0$ takes values in $\mathbb{C}\setminus \{-1,1\}$ and $0$ is a regular value of its imaginary part. Put $f= ({\sf f} - i \varepsilon)\mid X^0$. If $\varepsilon$ is small enough, then the irreducible mapping $f$ on $X^0$ is free homotopic to ${\sf f}\mid X^0$.
We identify the fundamental groups of $X$ and of $X^0$ by the inclusion mapping from $X^0$ to $X$.

Proposition \ref{prop2} applied to the mapping $f:X^0\to \mathbb{C}\setminus\{-1,1\}$
provides a M\"obius transformation $M_l$, that maps one of the components of $\mathbb{R}\setminus \{-1,1\}$ onto $(-1,1)$, and further a point $q\in X^0$ and a curve $\alpha$ in $X^0$ with initial point $q_0$ and terminating point $q$, such for the mapping
%Associate to $f$ the function
$F=M_l \circ f$ the point $q'\stackrel{def}= F(q)$ is contained in $(-1,1)$, and for the generators $e_j\stackrel{def}={\rm Is}_{\alpha}(e_{j,0}),\;  e_{j,0}\in \mathcal{E},$ of $\pi_1(X^0,q)$ the inequalities \eqref{eq3b} hold.
After identifying the fundamental groups $\pi_1(X,q)$ and $\pi_1(X,q_0)$ with different base point by an isomorphism that is defined up to conjugation, Theorem E states that the free homotopy class of $F$ corresponds to a conjugacy class of homomorphisms
$$
\pi_1(X,q_0)\cong \pi_1(X,q)\to \pi_1(\mathbb{C}\setminus \{-1,1\},q')\,,
$$
%$$
%\cong \pi_1(\mathbb{C}\setminus \{-1,1\},q')
%$$
that is represented by a homomorphism $h$ for which $\mathcal{L}_-(h(e_{j,0})) \leq 8\pi \lambda_7(X^0)$ for each $e_{j,0}\in \mathcal{E}$.
More explicitly, there exists
a smooth mapping $\tilde{F}:X^0\to \mathbb{C}\setminus\{-1,1\}$, that is free homotopic to $F$, maps $q_0$ to $q'$, and satisfies the inequality
\begin{equation}\label{eq3b+}
\mathcal{L}_-(\tilde{ F}_*(e_{j,0})) \leq 8\pi \lambda_7(X^0)
\end{equation}
for each $e_{j,0}\in \mathcal{E}$.
The existence of the smooth mapping $\tilde{F}$ can be seen explicitly as follows.
Write $e=\mbox{Is}_{\alpha}(e_{0}) \in \pi_1(X,q)$  for each $e_{0}\in \pi_1(X,q_0)$.
%for a curve $\alpha$ with initial point $q_0$ and
%terminating point $q$.
Parameterise $\alpha$ by the interval $[0,1]$.
The  image of $\alpha$  %in $\mathbb{C} \setminus  \{-1,1\}$
under the mapping $F$ is the curve $\beta= F \circ \alpha$ in $\mathbb{C} \setminus \{-1,1\}$  with initial point $F(q_0)$ and terminating point $F(q)=q'$. Then $F_*(e_{0})=(\mbox{Is}_{\beta})^{-1}(F_*(e))$.
Choose a homotopy $F_t,\, t \in [0,1]$, that joins the mapping $F_0 \stackrel{def}=F$ with a (smooth) mapping $F_1$ denoted by $\tilde F$, so that $F_t(q_0)=\beta(t),\, t\in [0,1]$.
%$\tilde{F} \stackrel{def}=F_1$ and moves
The value $\beta(t)$
%$\beta(t)\stackrel{def}=F_t(q_0)$
moves from the point $\beta(0)=F(q_0)$ to $\beta(1)=q'$ along the curve $\beta$.
Then $\tilde F(q_0)=q'$ and $\tilde F_*(e_{0})=F_*(e)$ for each $e_0 \in \pi_1(X,q_0)$.
Indeed,
denote by $\beta_t$ the curve that runs from $\beta(t)$ to $\beta(1)$ along $\beta$. Then $\beta_0=\beta$ and $\beta_1$ is a constant curve.
Let $\gamma_0$ be a curve that represents $e_0$.
The base point of the curve $F_t(\gamma_0)$ equals $F_t(q_0)=\beta(t)$. Hence, we obtain a continuous family of curves $\beta_t^{-1} F_t(\gamma_0) \beta_t$ with base point $\beta(1)=F(q)$. For $t=1$ the curve is equal to $F_1(\gamma_0)=\tilde{F}(\gamma_0)$, for $t=0$ the curve is equal to $\beta^{-1} F_0(\gamma_0) \beta=F_0(\alpha^{-1} {\rm Is}_{\alpha}(\gamma_0)\alpha)=F_0(\gamma)$. Since the two curves
$F_1(\gamma_0)$ and $F_0(\gamma)$
are homotopic and $F_1=\tilde F$, $F_0=F$, we obtain $\tilde{F}_*(e_0)=F_*(e)$.
The inequalities \eqref{eq3b}
%gives for each $e_{j,0}\in \mathcal{E}$
imply the inequalities
\eqref{eq3b+}.

For each irreducible holomorphic mapping ${\sf f}:X\to \mathbb{C}\setminus \{-1,1\}$ we found a M\"obius transformation $M_l$ and a mapping $\tilde F :X^0\to \mathbb{C}\setminus \{-1,1\}$ that satisfies the condition ${\tilde F} (q_0)=q'\in (-1,1)$ and the inequalities \eqref{eq3b+}, and is free homotopic on $X^0$ to $M_l(({\sf f}-i \varepsilon)|X^0)$ for a small number $\varepsilon$, and, hence, is free homotopic to
$M_l({\sf f}|X^0)$. Using a deformation retraction we obtain a mapping $\tilde{F}^X:X\to \mathbb{C}\setminus \{-1,1\}$ that is free homotopic on $X^0$ to $\tilde{F}$ and, hence, to $M_l({\sf f}|X^0)$. Identifying the fundamental groups $\pi_1(X^0,q_0)$ and $\pi_1(X,q_0)$ by the homomorphism induced by the inclusion and applying Theorem E, we obtain for each irreducible holomorphic mapping ${\sf f}:X\to \mathbb{C}\setminus \{-1,1\}$ a M\"obius transformation $M_l$ and
a smooth mapping $\tilde{F}^X:X\to \mathbb{C}\setminus \{-1,1\}$ that is free homotopic to $M_l(\sf{ f})$ on $X$, and satisfies the condition ${\tilde F}^X (q_0)=q'\in (-1,1)$ and the inequalities \eqref{eq3b+}.

If ${\sf f}:X\to  \mathbb{C}\setminus \{-1,1\}$ is a contractible mapping, it is free homotopic to the function $\tilde{ F}^X\equiv 0$ on $X$, %which is constant on $X$ $X^0$ which is equal to $f\mid X^0$
and the inequalities \eqref{eq3b+} are automatically satisfied for the monodromies of $\tilde{F}^X$.

The number of free homotopy classes of mappings $X\to \mathbb{C}\setminus\{-1,1\}$, that contain a smooth mapping $\tilde{ F}^X$, which satisfies the condition $\tilde{F}^X(q_0)=q'\in (-1,1)$ and
the inequalities \eqref{eq3b+}, are estimated from above as follows.
By Lemma 1 of \cite{Jo4} there are at most
$\frac{1}{2}e^{24\pi \lambda_7(X^0)}+1\leq \frac{3}{2}e^{24\pi \lambda_7(X^0)}$
different reduced words $w \in \pi_1(\mathcal{C}\setminus\{-1,1\}),0)$ (including the identity) with $\mathcal{L}_-(w)\leq 8\pi \lambda_7(X^0)$.
Identify standard generators of $\pi_1(\mathbb{C}\setminus \{-1,1\},q')$ with standard generators of $\pi_1(\mathbb{C}\setminus \{-1,1\},0)$
by the canonical isomorphism.
%We obtain, that $\tilde F_*$ is an element of a set containing
We saw, that there are at most $(\frac{3}{2}e^{24\pi \lambda_7(X^0)})^{2g+m}$ different homomorphisms $h:\pi_1(X^0,q_0)\cong \pi_1(X,q_0) \to \pi_1(\mathbb{C}\setminus \{-1,1\},q') \cong \pi_1(\mathbb{C}\setminus \{-1,1\},0)$ with $\mathcal{L}_-(h(e))\leq 8\pi \lambda_7(X^0)$ for each element $e$ of the set of generators $\mathcal{E}$ of $\pi_1(X^0,q_0)$.
By Theorem E there are at most $(\frac{3}{2}e^{24\pi \lambda_7(X^0)})^{2g+m}$ different
free homotopy classes of mappings $X\to \mathbb{C}\setminus \{-1,1\}$, that contain a smooth mapping $\tilde{ F}^X$ which satisfies the condition $\tilde{F}^X(q_0)=q'\in (-1,1)$ and
the inequalities \eqref{eq3b+}.

For each irreducible or contractible mapping $\sf{f}$ on $X$ one of the mappings $M_l\circ \sf{f}$, $l=0,1,2,$
is free homotopic to a mapping $\tilde{ F}^X$ which satisfies the condition $\tilde{F}^X(q_0)=q'\in (-1,1)$ and
the inequalities \eqref{eq3b+}. The $(M_l)^{-1}\circ\tilde{ F}^X$ represent at most
$3(\frac{3}{2}e^{24 \pi \lambda_7(X^0)})^{2g+m}$ free homotopy classes of irreducible or contractible mappings $X\to \mathbb{C}\setminus \{-1,1\}$.
% containing a holomorphic mapping $\sf f$.
Theorem 1 is proved with the upper bound $3(\frac{3}{2}e^{24 \pi \lambda_7(X^0)})^{2g+m}$ for an arbitrary relatively compact domain $X^0\subset X$ that is a deformation retract of $X$.

It remains to prove that $\lambda_7(X)= \inf\{\lambda_7(X^0): X^0 \Subset X \; \mbox{is\, a\, deformation\, retract\, of\,} X\,\}$.
We have to prove that for each $e_0 \in \pi_1(X,q_0)$ the quantity $\lambda(A(\widehat{e_0}))=\lambda(\widetilde{X} \diagup ({\rm Is}^{\tilde{q}_0})^{-1}(\langle e_0 \rangle))$ is
equal to the infimum of $\lambda( \widetilde{X^0}\diagup ({\rm Is}^{\tilde{q}_0})^{-1}(\langle e_0 \rangle))$ over all $X^0$ being open relatively compact subsets of $X$ which are deformation retracts of $X$. Here $\widetilde{X^0}$ is the universal covering of $X^0$, and the fundamental groups of $X$ and $X^0$ are identified. $\widetilde{X^0}$ ($\widetilde{X}$, respectively) can be defined as set of homotopy classes of arcs in $X^0$ (in $X$, respectively) joining $q_0$ with a point $q\in X^0$ (in $X$ respectively) equipped with the complex structure induced by the projection to the endpoint of the arcs, and the point $\tilde{q}_0$ corresponds to the class of the constant curve.
The isomorphism $({\rm Is}^{\tilde{q}_0})^{-1}$ from $\pi_1(X^0,q_0)$ to the group of covering transformations on $\widetilde{X^0}$ is defined in the same way as it was done for $X$ instead of $X^0$.
These considerations imply that there is a holomorphic mapping from $\widetilde{X}^0 \diagup ({\rm Is}^{\tilde{q}_0})^{-1}(\langle e_0 \rangle)$ into $\widetilde{X} \diagup ({\rm Is}^{\tilde{q}_0})^{-1}(\langle e_0 \rangle)$. Hence, the extremal length of the first set is not smaller than the extremal length of the second set.

Vice versa, take any annulus $A^0$ which is a relatively compact subset of $A(\widehat{e_0})$ and is a deformation retract of $A(\widehat{e_0})$. Its projection to $X$ is relatively compact in $X$, hence, it is contained in a relatively compact deformation retract $X^0$ of $X$. Hence, $A^0$ can be considered as subset of
$\widetilde{X^0} \diagup ({\rm Is}^{\tilde{q}_0})^{-1}(\langle e_0 \rangle)$, and, hence, $\lambda(\widetilde{X^0} \diagup ({\rm Is}^{\tilde{q}_0})^{-1}(\langle e_0 \rangle)) \leq \lambda(A^0)$. Since $\lambda (A(\widehat{e_0})) = \inf \{\lambda (A^0): A^0 \Subset A(\widehat{e_0}) \mbox{\,is\, a\, deformation\, retract\, of\,}A(\widehat{e_0}) \} $ we are done. \hfill $\Box$

\medskip

We proved a slightly stronger statement, namely, the number of homotopy classes of mappings $X\to \mathbb{C}\setminus\{-1,1\}$ that contain a contractible holomorphic mapping or an irreducible holomorphic mapping does not exceed $(\frac{3}{2}e^{24\pi \lambda_7(X)})^{2g+m}$.

\medskip

\noindent {\bf Proof of Proposition \ref{prop2a}.}
Denote the zero set $\{x\in X: \mbox{Im}f(x)=0\}$ by $L$. Since $f$ is not contractible,  $L\neq \emptyset$.

\noindent {\bf 1. A torus with a hole.} Assume first that $X$ is a torus with a hole with  base point $q_0$.
For notational convenience we denote by $e_0$ and ${e'_0}$ the two elements of the set of generators $\mathcal{E}$ of $\pi_1(X,q_0)$ that is associated to a standard bouquet of circles for $X$.
%was chosen in Section 1.
We claim that there is a connected component $L_0$ of $L$ such that (after suitable orientation) the intersection number of the free homotopy class  of one of the elements of $\mathcal{E}$,        %$\pi_1(X,q_0)$,
say of $\widehat{e_0}$ with $L_0$ is positive, and
%$L_0$ has positive
the intersection number with one of the classes $\reallywidehat{e'_0}$, or $\reallywidehat{(e'_0)^{-1}}$, or $\reallywidehat{e_0 \,e'_0}$ with $L_0$ is positive.

The claim is easy to prove in the case when there is a component of $L_0$ which is a simple closed curve that is not contractible and not contractible to the hole of $X$. Indeed, consider the inclusion of $X$ into a closed torus $X^c$ and the homomorphism on fundamental groups $\pi_1(X,q_0)\to \pi_1(X^c,q_0)$ induced by the inclusion. Denote by $e_0^c$ and ${e'_0}^c$ the images of $e_0$ and $e_0'$ under this homomorphism. Notice that $e_0^c$ and ${e'_0}^c$ commute. The (image under the inclusion of the) curve $L_0$ is a simple closed non-contractible curve in $X^c$. It represents the free homotopy class of an element $(e_0^c)^j({e'_0}^c)^k$ for some integers $j$ and $k$ which are not both equal to zero. Hence, $L_0$ is not null-homologous in $X^c$, and by the Poincar\'{e} Duality Theorem
for one of the generators, say for $e_0^c$, the representatives of the free homotopy class $\reallywidehat{e_0^c}$  have non-zero intersection number %with $L_0$, or
with $L_0$.
%(Then $k_0\neq 0$.)
After suitable orientation of $L_0$, we may assume that this intersection number is positive. There is a representative of the class $\reallywidehat{e_0^c}$ which is contained in $X$, hence,
$\widehat{e_0}$ has positive intersection number with $L_0$.

Suppose all compact connected components of $L$ are contractible or
contractible to the hole of $X$.
Consider a relatively compact domain ${X}''\Subset {X}$ in ${X}$ with smooth boundary which is a deformation retract of ${X}$ such that for each connected component of $L$ at most one component of its intersection with ${X}''$ is not contractible to the hole of ${X}''$. (See the paragraph on ''Regular zero sets''.) There is at least one component of $L\cap {X}''$ that is not contractible to the hole of ${X}''$. Indeed, otherwise the free homotopy class of each element of $\mathcal{E}$ could be represented by a loop avoiding $L$, and, hence, the monodromy of $f$ along each
element of $\mathcal{E}$ would be conjugate to the identity, and, hence, equal to the identity, i.e. contrary to the assumption, $f:X\to \mathbb{C}\setminus \{-1,1\} $ would be free homotopic to a constant.

Take a component $L_0''$ of $L\cap {X}''$ that is not contractible to the hole of ${X}''$.
There is an arc of $\partial {X}''$ between the endpoints of $L_0''$ such that
the union $\tilde{L}_0$ of the component $L_0''$ with this arc is a closed curve in ${X}$ that is
not contractible and not contractible to the hole. Hence, for one of the elements of $\mathcal{E}$, say for $e_0$, the intersection number of the free homotopy class $\widehat{e_0}$ with the closed curve $\tilde{L}_0$ is positive after orienting the curve $\tilde{L}_0$ suitably. We may take a representative $\gamma$ of $\widehat{e_0}$ that is
contained in ${X}''$. Then $\gamma$ has positive intersection number with $L_0''$.
Denote the connected component of $L$ that contains $L_0''$ by $L_0$. All components of $L_0\cap {X}''$ different from $L_0''$ are contractible to the hole of ${X}''$. Hence, $\gamma$ has intersection number zero with each of these components.
Hence, $\gamma$ has positive intersection number with $L_0$ since ${\gamma}\subset {X}''$.
We proved that the class  $\widehat{e_0}$
has positive intersection number with $L_0$.

If $\reallywidehat{e'_0}$ also has non-zero intersection number with $L_0$  we define $e''_0=(e'_0)^{\pm 1}$ so that the intersection number of  $\reallywidehat{e''_0}$ with $L_0$ is positive. If $\reallywidehat{e'_0}$  has zero intersection number with $L_0$ we put $e''_0=e_0 \,e'_0$. Then again the intersection number of $\reallywidehat{e''_0}$ with $L_0$ is positive.  Also, the intersection number of  $\reallywidehat{e_0 \, e''_0}$ with $L_0$ is positive. The set $\mathcal{E}_2'\stackrel{def}=\{e_0,e_0''\}$ satisfies the condition required in the proposition.
We obtained Proposition \ref{prop2a} for a torus with a hole.

\noindent {\bf 2. A planar domain.} Let $X$ be a planar domain.
The domain $X$ is conformally equivalent to  a disc with $m$ smoothly bounded holes,
%$\mathcal{C}_j,\, j=1,\ldots m,$
equivalently, to the Riemann sphere with $m+1$ smoothly bounded holes, $\mathbb{P}^1 \setminus \bigcup_{j=1}^{m+1}\mathcal{C}_j$, where $\mathcal{C}_{m+1}$ contains the point $\infty$. As before the base point of $X$ is denoted by $q_0$, and for each $j=1,\ldots,m,$ the generator $e_{j,0}\in \mathcal{E}\subset \pi_1(X,q_0)$
is represented by a curve that surrounds $\mathcal{C}_j$ once counterclockwise.
Since $f$ is not contractible, there must be a connected component of $L$ that has limit points on some $\mathcal{C}_j$ with $j\leq m$. Indeed, otherwise
%there would be representatives
the free homotopy class of each generator could be represented by a curve that avoids $L$.
This would imply that all monodromies are equal to the identity.
We claim that
there exists a component $L_0$ of $L$
with limit points on the boundary of two components $\partial {\mathcal{C}}_{j'}$ and
$\partial {\mathcal{C}}_{j''}$ for some $j', j'' \in \{1,\ldots,m+1\}$ with $j''\neq j'$.

Indeed, assume the contrary.
Then, if a component of $L$ has limit points on a component $\partial {\mathcal{C}}_{j},\, j\leq m,$ then all its limit points are on $\partial {\mathcal{C}}_{j}$.
Take a smoothly bounded simply connected domain ${\mathcal{C}}'_{j}\Subset X\cup {{\mathcal{C}}_{j}}$
that contains the closed set ${{\mathcal{C}}_{j}} $ , so that its boundary $\partial {\mathcal{C}}'_{j}$ represents $\reallywidehat{e_{j,0}}$.
Then all components $L'_k$ of $L\setminus {\mathcal{C}}'_{j} $ with an endpoint on $\partial \mathcal{C}'_{j}$ have another endpoint on this circle. The two endpoints of $L_k'$ on $\partial {\mathcal{C}}'_{j}$ divide $\partial {\mathcal{C}}'_{j}$ into two connected components.
The union of $\overline{L_k'}$ with each of the two components of $\partial {\mathcal{C}}'_{j}\setminus \overline{L_k'}$ is a simple closed curve in $\mathbb{C}$, and, hence, by the Jordan Curve Theorem it bounds a relatively compact topological disc in $\mathbb{C}$. One of these discs contains $ {\mathcal{C}}'_{j}$, the other does not.
%Fix a point $p'_j\in \partial \mathcal{C}'_{j} \setminus L$.
Assign to each component $L_k'$ of $L \setminus {\mathcal{C}}'_{j}   $ with both endpoints on $\partial \mathcal{C}'_{j} $ the closed arc $\alpha_k$  in $\partial \mathcal{C}'_{j}$ with the same endpoints as $L_k'$, whose union with $L_k'$ bounds a relatively compact topological disc in $\mathbb{C}$ that does not contain $\mathcal{C}'_{j}$. These discs are partially ordered by inclusion,
since the $L_k'$ are pairwise disjoint. Hence, the arcs $\alpha_k$ are partially ordered by inclusion. For an arc $\alpha_k$  which contains no other of the arcs (a minimal arc) the curve $f \circ \alpha_k$ except its endpoints is contained in  $\mathbb{C}\setminus \mathbb{R}$. Moreover,
the endpoints of  $f \circ \alpha_k$ lie on ${f(\overline{L_k'})}$, which is contained in one connected components of $\mathbb{R}\setminus \{-1,1\}$, since $\overline{L_k'}$ is connected. Hence, the curve $f \circ \alpha_k$ is homotopic  in $\mathbb{C}\setminus \{-1,1\}$ (with fixed endpoints) to a curve in $\mathbb{R} \setminus \{-1,1\}$. The function $f$ either maps all points on $\partial \mathcal{C}'_{j}\setminus \alpha_k$ that are close to $\alpha_k$ to the open upper half-plane or maps them all to the open lower half-plane. (Recall, that zero is a regular value of $\mbox{Im}f$.) Hence, for an open arc $\alpha'_k \subset \partial \mathcal{C}'_{j} $ that contains $\alpha_k$ the curve $f \circ \alpha'_k $ is homotopic  in  $\mathbb{C}\setminus \{-1,1\}$ (with fixed endpoints) to a curve in $\mathbb{C}\setminus \mathbb{R}$.

Consider the arcs $\alpha_k$ with the following property. For an open arc $\alpha'_k$  in $\partial \mathcal{C}'_{j} $ which contains the closed arc $\alpha_k$ the mapping $f \circ \alpha'_k$ is homotopic in $\mathbb{C}\setminus \{-1,1\}$ (with fixed endpoints) to a curve contained in  $\mathbb{C}\setminus \mathbb{R}$. Induction on the arcs by inclusion shows that this property is satisfied for all maximal arcs among the $\alpha_k$ and, hence, $f \mid \partial \mathcal{C}'_{j} $ is contractible in  $\mathbb{C}\setminus \{-1,1\}$. %We saw that
Hence, if the claim was not true, then for each hole $\mathcal{C}_j,\, j\leq m$, whose boundary contains limit points of a connected component of $L$, the monodromy along the curve $\mathcal{C}'_j$ (with any base point) that represents $\reallywidehat{e_{j,0}}$
would be trivial. Then all monodromies would be trivial, which contradicts the fact that the mapping is not contractible. The contradiction proves the claim.

With $j'$ and $j''$ being the numbers of the claim and $j'\leq m$ we consider the set $\mathcal{E}_2' \subset \mathcal{E}_2$  which consists of the following primitive elements:
$e_{j',0}$, the element $(e_{j'',0})^{-1}$ provided $j'' \neq m+1$, and  $e_{j',0}\, e_{j,0}$ for all $j=1,\ldots,m,\, j\neq j', j \neq j''$. The free homotopy class of each element of $\mathcal{E}_2'$  has intersection number $1$ with $L_0$ after suitable orientation of the curve $L_0$.
Each product of at most two different elements of $\mathcal{E}_2'$ is a primitive element of $\pi_1(X,q)$ and is contained in $\mathcal{E}_4$. Moreover, the intersection number with $L_0$ of the free homotopy class of each product of at most two different elements of $\mathcal{E}_2'$
%any such element
equals $1$ or $2$.
Each element of $\mathcal{E}$ is the product of at most two elements of $\mathcal{E}_2'\cup (\mathcal{E}_2')^{-1}$.

The proposition is proved for the case of planar domains $X$.
% with $\lambda$ replaced by $\lambda_4$.
\hfill $\Box$

\medskip

\noindent {\bf Proof of Proposition \ref{prop2}.}

\noindent {\bf 1.  A torus with a hole.}
Consider the curve $L_0$ and the set $\mathcal{E}_2' \subset \pi_1(X,q_0)$ obtained in Proposition \ref{prop2a}. For one of the functions $M_l \circ f$, denoted by $F$, the image $F (L_0)$ is contained in $(-1,1)$.
Let $e_0$, $e_0'$ be the two elements of $\mathcal{E}$.
Move the base point $q_0$ to a point $q\in L_0$ along a curve $\alpha$ in $X$, and consider the generators $e=\mbox{Is}_{\alpha} (e_0)$ and
$e'=\mbox{Is}_{\alpha} ({e'_0})$ of $ \pi_1(X,q)$, and the set $\mbox{Is}_{\alpha}(\mathcal{E}_2')\subset \pi_1(X,q)$. Then $e$ and $e'$ are products of at most two elements of $\mbox{Is}_{\alpha}(\mathcal{E}_2')$.
Since the free homotopy class of an element of $\pi_1(X,q_0)$ coincides with the free homotopy class of the element of $\pi_1(X,q)$ obtained by applying $\mbox{Is}_{\alpha}$,
the free homotopy class of each product of one or two elements of $\mbox{Is}_{\alpha}(\mathcal{E}_2')$ intersects $L_0$.
We may assume as in the proof of Proposition \ref{prop2a} that  $\mbox{Is}_{\alpha}(\mathcal{E}_2')$ consists of the elements $e$ and $e''$, where $e''$ is either equal to ${e'}^{\pm 1}$, or equals the product of $e$ and $e'$.
%at most two elements among the $e$ and $e'$ and their inverses.
Lemma \ref{lem3} applies to the pair $e$, $e''$,
the function $F$,
and the curve $L_0$.
%By the conditions of Proposition \ref{prop2} and Lemma 4 of \cite{Jo5}
Since $F$ is irreducible, the monodromies of $F$ along
$e$ and $e''$ are not powers of a single standard generator of the fundamental group of $\pi_1(\mathbb{C}\setminus \{-1,1\},q')$.
%Since for each $\tilde{e}_0\in \pi_1(X,q_0)$ the equality
%$A(\widehat{\tilde{e_0}})=A(\mbox{Is}_{\alpha}(\widehat{\tilde{e_0}}))$
%holds,this implies that
Hence, the monodromy along each of the $e$ and $e''$ is the product of at most two elements of $\mathcal{L}_-$ not exceeding $2 \pi\lambda_{e,e''}$. Therefore,
the monodromy of $F$ along each of the $e$ and $e''$ has $\mathcal{L}_-$ not exceeding $4 \pi\lambda_{e,e''}$. Notice that $\lambda_{e,e''}=\lambda_{e_0 e_0''}\leq \lambda_3(X)$, since $e_0''$ is the product of at most two factors, each an element of $\mathcal{E}\cup \mathcal{E}^{-1}$. Since $e'$ is the product of at most two different elements among the $e$ and $e''$ and their inverses, we obtain Proposition \ref{prop2} for $e$ and $e'$, in particular  $\mathcal{L}_-(F_*(e))$ and $\mathcal{L}_-(F_*( e'))$ do not exceed $8\pi \lambda_3(X)$.
Proposition \ref{prop2} is proved for tori with a hole.

\noindent {\bf 2. A planar domain.} Consider the curve $L_0$ and the set $\mathcal{E}_2'$ of Proposition \ref{prop2a}.
Move the base point $q_0$ along an arc $\alpha$ to a point $q \in L_0$. Then $f(q) \in \mathbb{R}\setminus \{-1,1\}$ and for one of the mappings $f$, $M_1 \circ f$, or $M_2 \circ f$, denoted by $F$, the inclusion $F(L_0)\subset (-1,1))$ holds, hence, $q'\stackrel{def}=F(q)$ is contained in $(-1,1)$. Denote  $e_j =  \mbox{Is}_{\alpha}(e_{j,0})$ for each element $e_{j,0}\in \mathcal{E}$. The $e_j$ form the basis $\mbox{Is}_{\alpha}(\mathcal{E})$ of $\pi_1(X,q)$. The set $\mbox{Is}_{\alpha}(\mathcal{E}_2')$ consists of primitive elements of $\pi_1(X,q)$ such that the free homotopy class of each product of one or two elements of $\mbox{Is}_{\alpha}(\mathcal{E}_2')$ intersects $L_0$. Moreover, each element of $\mbox{Is}_{\alpha}(\mathcal{E})$ is the product of one or two elements of  $\mbox{Is}_{\alpha}(\mathcal{E}_2') \cup (\mbox{Is}_{\alpha}(\mathcal{E}_2'))^{-1}$.

By the condition of the proposition not all monodromies $F_*(e),\, e \in \mbox{Is}_{\alpha}(\mathcal{E}_2'),$ are (trivial or non-trivial) powers of the same standard generator of $\pi_1(\mathbb{C}\setminus \{-1,1\},q')$.
Apply Lemma \ref{lem3} to all pairs of elements of $\mbox{Is}_{\alpha}(\mathcal{E}_2')$  whose monodromies are not (trivial or non-trivial) powers of the same standard generator of $\pi_1(\mathbb{C} \setminus \{-1,1\},q')$.
Since the product of at most two different elements of $\mbox{Is}_{\alpha}(\mathcal{E}_2')$ is contained in  $\mbox{Is}_{\alpha}(\mathcal{E}_4)$, Lemma \ref{lem3} shows that the monodromy $F_*(e)$ along each element  $e \in \mbox{Is}_{\alpha}(\mathcal{E}_2')$ is the product of at most two factors, each with $\mathcal{L}_-$ not exceeding $2 \pi \lambda_4(X)$.
Since each element of $\mbox{Is}_{\alpha}(\mathcal{E}) $ is a product of at most two factors in $\mathcal{E}_2' \cup (\mathcal{E}_2')^{-1}$,
%, each factor or its inverse an element of $E$,
the  monodromy $F_*(e_j)$ along each generator $e_j$ of $\pi_1(X,q)$ is the product of at most $4$ factors of $\mathcal{L}_-$ not exceeding  $2 \pi \lambda_4(X)$, and, hence, each monodromy $F_*(e_j)$ has $\mathcal{L}_-$ not exceeding $8 \pi \lambda_4(X)$.
Proposition \ref{prop2} is proved for planar domains.

\noindent {\bf 3.1. The general case. Diagrams of coverings.}
We will use diagrams of coverings to reduce this case to the case of a torus with a hole or to the case of the Riemann sphere with three holes.

Let as before $\tilde{q}_0$ be the point in $\tilde X$ with ${\sf P}(\tilde{q}_0)=q_0$
chosen in Section \ref{sec:1a} . Let $N$ be a subgroup of the fundamental group $\pi_1(X,q_0)$ and let $\omega^N:\tilde{X}\to \tilde{X}\diagup({\rm Is}^{\tilde{q}_0})^{-1}(N)=X(N)$ be the projection defined in Section \ref{sec:1a}. Put $(q_0)_N \stackrel{def}= \omega^N(\tilde{q}_0)$. For an element $e_0\in N\subset \pi_1(X,q_0)$ we denote by $(e_0)_N$ the element of $\pi_1(X(N),(q_0)_N)$ that is obtained as follows. Take a curve $\gamma$ in $X$ with base point $q_0$ that represents $e_0\in N$. Let
$\tilde \gamma$ be its lift to $\tilde X$ with terminating point $\tilde{q}_0$.
Then $\gamma_N\stackrel{def}=\omega^N(\tilde{\gamma})$ is a closed curve in $X(N)=\tilde{X}\diagup ({\rm Is}^{\tilde{q}_0})^{-1}(N) $ with base point $(q_0)_N$. The element of $\pi_1(X(N), (q_0)_N)$ represented by $\gamma_N$ is the required element $(e_0)_N$. All curves $\gamma'_N$  representing $(e_0)_N$ have the form $\omega^N(\tilde{\gamma}')$ for a curve $\tilde{\gamma}'$ in $\tilde X$ with terminating point $\tilde{q}_0$ and initial point $({\rm Is}^{\tilde{q}_0})^{-1}(e_0)(\tilde{q}_0)$.
%${\sf P}((\tilde{\gamma}')$ being a closed curve in $X$ that represents $e_0$.
Since $\omega_N\circ\omega^N={\sf P}$, the curve $\omega_N(\gamma'_N)={\sf P}(\tilde{\gamma}')={\gamma}'$
%for a curve ${\gamma}'$ in $X$ that
represents $e_0$ for each curve $\gamma'_N$ in $X(N)$ that represents $(e_0)_N$. We obtain $(\omega_N)_*((e_0)_N)=e_0$. For two subgroups $N_1\leq N_2$ of $\pi_1(X,q_0)$ we obtain $(\omega^{N_2}_{N_1})_*((e_0)_{N_1})=(e_0)_{N_2}$, $e_0\in N_1$ (see the commutative diagram Figure \ref{fig3}). %In particular,

Let $\tilde q$ be another base point of $\tilde X$ and let $\tilde \alpha$
be a curve in $\tilde X$ with initial point $\tilde{q}_0$ and terminating point $\tilde q$. Let again $N$ be a subgroup of $\pi_1(X,q_0)$.  Put $q_N\stackrel{def}=\omega^N(\tilde{q})$.
The curve $\alpha_N=\omega^N(\tilde{\alpha})$ in $X(N)$, and the base point $\tilde q$ of $\tilde X$ are compatible, hence, $({\rm Is}^{\tilde{q}_0})^{-1}(N)=({\rm Is}^{\tilde{q}})^{-1}({\rm Is}_{\alpha_N}(N))$ and $X({\rm Is}_{\alpha_N}(\sigma)  )=\tilde{X}\diagup ({\rm Is}^{\tilde{q}})^{-1}({\rm Is}_{\alpha_N}(\sigma)),\, \sigma\in N,$ is canonically isomorphic to $X(N)=\tilde{X}\diagup ({\rm Is}^{\tilde{q}_0})^{-1}(N) $.

We will use the previous notation $\omega^{N_2}_{N_1}$ also for the projection
%$X({\rm Is}_{\alpha_{N_j}}(N_j)    )=\tilde{X}\diagup ({\rm Is}^{\tilde{q}})^{-1}({\rm %Is}_{\alpha_{N_j}}(N_j)) $  will
%be denoted as before by $\omega_{N_1}^{N_2}$.
$\tilde{X}\diagup ({\rm Is}^{\tilde{q}})^{-1}({\rm Is}_{\alpha_{N_1}}(N_1))\to \tilde{X}\diagup ({\rm Is}^{\tilde{q}})^{-1}({\rm Is}_{\alpha_{N_2}}(N_2)) $, $N_1\leq N_2$ being subgroups of $\pi_1(X,q_0)$ ($N_1$ may be the identity and $N_2$ may be $\pi_1(X,q_0)$.)

%For a subgroup $N$ of $\pi_1(X,q_0)$ we put
%$q_N\stackrel{def}=\omega^N(\tilde{q})$ and
Put $\alpha\stackrel{def}={\sf P}(\tilde{\alpha})$. For an element $e_0\in \pi_1(X,q_0)$ we put $e\stackrel{def}={\rm Is}_{\alpha}(e_0)\in \pi_1(X,q)$ and denote by $e_N$ the element of $\pi_1(X(N),q_N)$, that is represented by $\omega^N(\tilde{\gamma})$ for a curve $\tilde \gamma$ in $\tilde X$ with terminating point $\tilde q$ and projection ${\sf P}(\tilde{\gamma})$ representing $e$. Again $(\omega^{N_2}_{N_1})_*(e_{N_1})=e_{N_2}$ for subgroups $N_1\leq N_2$ of $\pi_1(X,q)$ and $e\in N_1$, in particular $( \omega_{N})_*(e_{N})=e$ for a subgroup $N$ of $\pi_1(X,q)$ and $e\in N$.

\noindent {\bf 3.2. The estimate for a chosen pair of monodromies.}
%{\bf Reduction to the case of a torus with a hole or a planar domain.}
Since the mapping $f:X\to \mathbb{C}\setminus\{-1,1\}$ is irreducible, there exist two elements $e_0',\, e_0''\in \pi_1(X,q_0)$ such that the monodromies $f_*(e_0')$ and $f_*(e_0'')$ are not powers of a single conjugate of a power of one of the elements $a_1$, $a_2$ or $a_1a_2$. The fundamental group of the Riemann surface $X(\langle e_0',e_0''\rangle)$ is a free group in the two generators $(e'_0)_{\langle e_0',e_0''\rangle}$ and $(e''_0)_{\langle e_0',e_0''\rangle}$, hence, $X(\langle e_0',e_0''\rangle)$ is either a torus with a hole or is equal to $\mathbb{P}^1$ with three holes. Moreover, the system $\mathcal{E}_{\langle e_0',e_0''\rangle}=\{(e'_0)_{\langle e_0',e_0''\rangle}, \,(e''_0)_{\langle e_0',e_0''\rangle}\}$  of generators of the fundamental group $\pi_1(X(\langle e_0',e_0''\rangle),(q_0)_{\langle e_0',e_0''\rangle})$ is associated to a standard bouquet of circles for $X(\langle e_0',e_0''\rangle)$.
This can be seen as follows. The set of generators $\mathcal{E}$ of $\pi_1(X,q_0)$ is associated to a standard bouquet of circles for $X$. For each $e_0\in \mathcal{E}$ we denote the circle of the bouquet that represents $e_0$ by $\gamma_{e_0}$.
For each $e_0\in\mathcal{E}$ we lift the circle $\gamma_{e_0}$ with base point $q_0$ to an arc $\widetilde{\gamma_{e_0}}$ in
$\tilde X$ with terminating point
%of the lift of each circle being
$\tilde{q_0}$. Let $D$ be a small disc in $X$ around $q_0$, and $\widetilde{D_0}$, $\widetilde{D_{e_0}}, e_0\in \mathcal{E},$ be the preimages of $D$ under the projection %$\omega^{\langle e'_0,e''_0\rangle}:\tilde{X}\to X(\langle e'_0,e''_0\rangle)$,
${\sf P}:\tilde{X}\to X$,
that contain $\tilde{q_0}$, or the initial point of $\widetilde{\gamma_{e_0}}$, respectively.
We assume that $D$ is small enough so that the mentioned preimages of $D$ %in $\tilde X$
are pairwise disjoint. Put $D_{\langle e'_0,e''_0\rangle}=\omega^{\langle e'_0,e''_0\rangle}(\widetilde{D_0})$.
For $e_0\neq e'_0,e''_0$
the image
$\omega^{\langle e'_0,e''_0\rangle}(\widetilde{D_0}\cup \widetilde{\gamma_{e_0}}\cup \widetilde{D_{e_0}} )$
is the union of an arc $\omega^{\langle e'_0,e''_0\rangle} (\widetilde{\gamma_{e_0}})$ in $X(\langle e'_0,e''_0\rangle)$ with two disjoint discs, each containing an endpoint of the arc and one of them equal to $D_{\langle e'_0,e''_0\rangle}$.
For $e_0=e'_0,e''_0$ the image $\omega^{\langle e',e''\rangle}(\widetilde{D_0}\cup \widetilde{\gamma_{e_0}}\cup \widetilde{D_{e_0}} )$ is the union of $D_{\langle e'_0,e''_0\rangle}$ with the loop $(\gamma_{e_0})_{\langle e'_0,e''_0\rangle}\stackrel{def}= \omega^{\langle e'_0,e''_0\rangle}( \widetilde{\gamma_{e_0}})$. For $e_0=e'_0,e''_0$ the loop $(\gamma_{e_0})_{\langle e'_0,e''_0\rangle}$ in $X(\langle e'_0,e''_0\rangle)$
has base point $(q_0)_{\langle e'_0,e''_0\rangle}=\omega^{\langle e'_0,e''_0\rangle}(\tilde{q_0})$
%in $X(\langle e',e''\rangle)$  with base point $q_0$
and represents the generator $(e_0)_{\langle e'_0,e''_0\rangle}$ of the fundamental group of $\pi_1(X(\langle e'_0,e''_0\rangle),(q_0)_{\langle e'_0,e''_0\rangle})$.
Since
the bouquet of circles $\cup_{e_0\in \mathcal{E}}\, \gamma_{e_0}$ is a standard bouquet of circles for $X$, the union $(\gamma_{e'_0})_{\langle e'_0,e''_0\rangle} \cup (\gamma_{e''_0})_{\langle e'_0,e''_0\rangle}$ is a standard bouquet of circles in $X(\langle e'_0,e''_0\rangle)$.  This can be seen by looking at the intersections of the loops with a circle that is contained in $D_{\langle e'_0,e''_0\rangle}$ and surrounds $(q_0)_{\langle e'_0,e''_0\rangle}$. By the commutative diagram of coverings the intersection behaviour is the same as for the images of these objects under $\omega_{\langle e'_0,e''_0\rangle}$.
Hence, since $(\gamma_{e'_0})_{\langle e'_0,e''_0\rangle}$ and $(\gamma_{e''_0})_{\langle e'_0,e''_0\rangle}$ represent the generators $({e'_0})_{\langle e'_0,e''_0\rangle}$ and  $({e''_0})_{\langle e'_0,e''_0\rangle}$ of $\mathcal{E}_{\langle e_0',e_0''\rangle}$, the
union $(\gamma_{e'_0})_{\langle e'_0,e''_0\rangle} \cup (\gamma_{e''_0})_{\langle e'_0,e''_0\rangle}$ is a standard bouquet of circles for $X(\langle e'_0,e''_0\rangle)$.

The set $X(\langle e_0',e_0''\rangle)$ is either a torus with a hole or is equal to $\mathbb{P}^1$ with three holes.
Apply Proposition \ref{prop2a} to the Riemann surface $X(\langle e_0',e_0''\rangle)$ with base point $(q_0)_{\langle e'_0,e''_0\rangle}$, the holomorphic mapping $f_{\langle e_0',e_0''\rangle}=f\circ\omega_{\langle e_0',e_0''\rangle}$ into $\mathbb{C}\setminus \{-1,1\}$, and the set of generators $\mathcal{E}_{\langle e_0',e_0''\rangle}$ of the fundamental group $\pi_1(X(\langle e'_0,e''_0\rangle),(q_0)_{\langle e'_0,e''_0\rangle})$. We obtain a relatively closed curve $L_{\langle e_0',e_0''\rangle}$ on which the function $f_{\langle e_0',e_0''\rangle}$ is real,
and a set $(\mathcal{E}_{\langle e_0',e_0''\rangle})_2'=\{({\sf e}_0')_{\langle e_0',e_0''\rangle},({\sf e}_0'')_{\langle e_0',e_0''\rangle}\}$ which contains one of the elements of $\mathcal{E}_{\langle e_0',e_0''\rangle}$. The second element of  $(\mathcal{E}_{\langle e_0',e_0''\rangle})_2'$ is either equal to second element of $\mathcal{E}_{\langle e_0',e_0''\rangle}$ or to its inverse, or to the product of the two elements (in any order) of $\mathcal{E}_{\langle e_0',e_0''\rangle}$.  (We will usually refer to the product $ ({ e}_0')_{\langle e_0',e_0''\rangle}\,({e}_0'') _{\langle e_0',e_0''\rangle}=({ e}_0'\,{ e}_0'')_{\langle e_0',e_0''\rangle} $,
but we may change the product $({ e}_0'\,{ e}_0'')_{\langle e_0',e_0''\rangle} $ to the product $({ e}_0''\,{ e}_0')_{\langle e_0',e_0''\rangle} $, without changing the arguments and the estimate of the $\mathcal{L}_-$ of the monodromies of the elements of $\mathcal{E}_2'$.)
%will not depend on the order of the factors.)
The free homotopy classes $\reallywidehat{ ({\sf e}_0')_{\langle e_0',e_0''\rangle}}$, $\reallywidehat{ ({\sf e}_0'')_{\langle e_0',e_0''\rangle}}$, and $\reallywidehat{ ({\sf e}_0')_{\langle e_0',e_0''\rangle}\,({\sf e}_0'') _{\langle e_0',e_0''\rangle}}=\reallywidehat{({\sf e}_0'\,{\sf e}_0'')_{\langle e_0',e_0''\rangle}} $ intersect $L_{\langle e_0',e_0''\rangle}$.

Choose a point $q_{\langle e_0',e_0''\rangle}\in L_{\langle e_0',e_0''\rangle}$ and a point $\tilde{q}\in\tilde{X}$ with $\omega^{\langle e_0',e_0''\rangle}(\tilde{q})= q_{\langle e_0',e_0''\rangle}$. Let $\tilde{\alpha}$ be a curve in $\tilde X$ with initial point $\tilde{q}_0$ and terminating point $\tilde q$, and $\alpha_{\langle e_0',e_0''\rangle}= \omega^{\langle e_0',e_0''\rangle}(\tilde{\alpha})$. Put ${\sf e}'_{\langle e_0',e_0''\rangle}={\rm Is}_{\alpha_{\langle e_0',e_0''\rangle}}(({\sf e}_0')_{\langle e_0',e_0''\rangle})$ and ${\sf e}''_{\langle e_0',e_0''\rangle}={\rm Is}_{\alpha_{\langle e_0',e_0''\rangle}}(({\sf e}_0'')_{\langle e_0',e_0''\rangle})$. %By Proposition \ref{prop2}
For one out of three M\"obius transformations $M_l$ the mapping $F_{\langle e_0',e_0''\rangle}=M_l\circ f_{\langle e_0',e_0''\rangle}=M_l\circ f\circ\omega_{\langle e_0',e_0''\rangle}$ takes $L_{\langle e_0',e_0''\rangle}$ to $(-1,1)$, and hence $F_{\langle e_0',e_0''\rangle}$
takes a value $q'=F_{\langle e_0',e_0''\rangle}(q_{\langle e_0',e_0''\rangle})\in (-1,1)$ at $q_{\langle e_0',e_0''\rangle}$. By Lemma \ref{lem3} each of the $(F_{\langle e_0',e_0''\rangle})_*({\sf e}'_{\langle e_0',e_0''\rangle})$ and $(F_{\langle e_0',e_0''\rangle})_*({\sf e}''_{\langle e_0',e_0''\rangle})$ is the product of at most two elements of  $\pi_1(\mathbb{C}\setminus\{-1,1\},q')$     of $\mathcal{L}_-$ not exceeding $2\pi \lambda_3(X(\langle e_0',e_0''\rangle))$, hence,
\begin{align*}
\mathcal{L}_-&((F_{\langle e_0',e_0''\rangle})_*({\sf e}'_{\langle e_0',e_0''\rangle}))\leq 4\pi\lambda_3(X(\langle e_0',e_0''\rangle)),\\
\mathcal{L}_-&((F_{\langle e_0',e_0''\rangle})_*({\sf e}''_{\langle e_0',e_0''\rangle}))\leq 4\pi\lambda_3(X(\langle e_0',e_0''\rangle))\,.
\end{align*}
It follows that each of the $(F_{\langle e_0',e_0''\rangle})_*({ e}'_{\langle e_0',e_0''\rangle})$ and $(F_{\langle e_0',e_0''\rangle})_*({ e}''_{\langle e_0',e_0''\rangle})$ is the product of at most four elements of  $\pi_1(\mathbb{C}\setminus\{-1,1\},q')$     of $\mathcal{L}_-$ not exceeding $2\pi \lambda_3(X(\langle e_0',e_0''\rangle) )$, hence,
\begin{align*}
\mathcal{L}_-&((F_{\langle e_0',e_0''\rangle})_*({ e}'_{\langle e_0',e_0''\rangle}))\leq 8\pi\lambda_3(X(\langle e_0',e_0''\rangle)),\\
\mathcal{L}_-&((F_{\langle e_0',e_0''\rangle})_*({ e}''_{\langle e_0',e_0''\rangle}))\leq 8\pi \lambda_3(X(\langle e_0',e_0''\rangle))\,.
\end{align*}
It remains to take into account that for a subgroup $N$ of $\pi_1(X,q_0)$
the equation $(F_N)_*({{ e}}_{\langle e_0',e_0''\rangle})=F_*({e})$ holds for each ${e}\in {\rm Is}_{\alpha}(N)$, and
$\lambda_j(X(N))\leq \lambda_j(X)$ for each natural number $j$.

\noindent {\bf 3.3. Other generators. Intersection of free homotopy classes with a component of the zero set.}
%{Reduction to the case of a torus with a hole or a planar domain.}
Take any element $e\in {\rm Is}_{\alpha}(\mathcal{E})$  that is not in $\langle e',e''\rangle$.
%different from $e'$ and $e''$.
Then the monodromy $F_*(e)$ is either equal to the identity, or
one of the pairs  $(F_*(e),F_*({\sf e}'))  $ or  $(F_*(e),F_*({\sf e}''))  $ consists of two elements of $\pi_1(\mathbb{C}\setminus\{-1,1\},q')$ that are not powers of the same standard generator $a_j$, $j=1$ or $2$. Interchanging if necessary ${\sf e}'$ and ${\sf e}''$, we may suppose this option holds for the pair
$(F_*(e),F_*({\sf e}'))  $. Moreover, changing if necessary ${\sf e}'$ to its inverse $({\sf e}')^{-1}$, we may assume that ${\sf e}'$ is either an element of ${\rm Is}_{\alpha}(\mathcal{E})$ or it is a product of two elements of ${\rm Is}_{\alpha}(\mathcal{E})$. The quotient $X(\langle e,{\sf e}'\rangle)=\tilde{X}\diagup ({\rm Is}^{\tilde q})^{-1}(\langle e, {\sf e}'\rangle)$ is a Riemann surface whose fundamental group is a free group in two generators. Hence  $X(\langle e,{\sf  e}'\rangle)$ is either a torus with a hole or is equal to $\mathbb{P}^1$ with three holes.

We consider a diagram of coverings as follows. Let first $X(\langle {\sf  e}'\rangle)= \tilde{X}\diagup ({\rm Is}^{\tilde q})^{-1}(\langle {\sf e}'\rangle)$ be the annulus
with base point $q_{\langle {\sf e}'\rangle}=\omega^{\langle {\sf e}'\rangle}(\tilde{q})$,    that admits a mapping $\omega_{\langle {\sf e}'\rangle}: X(\langle {\sf  e}'\rangle)\to X$ that represents ${\sf  e}'$. By Lemma \ref{lem1} the connected component $L_ {\langle {\sf e}'\rangle}$ of $(\omega_{\langle {\sf e}'\rangle}^{\langle {\sf e}', {\sf e}''\rangle}    )^{-1}(L_{\langle {\sf e}', {\sf e}''\rangle})$ that contains $q_{\langle {\sf e}'\rangle}=\omega^{\langle {\sf e}'\rangle}(\tilde{q})$ is a relatively closed curve in $X(\langle {\sf  e}'\rangle)$ with limit points on both boundary components. The free homotopy class of the generator ${\sf  e}'_{\langle {\sf e}'\rangle}$ of $\pi_1(X(\langle {\sf e}'\rangle),q_{\langle {\sf e}'\rangle})$   intersects $L_ {\langle {\sf e}'\rangle}$. The mapping $F_{\langle {\sf e}'\rangle}=M_l \circ f\circ \omega_{\langle {\sf e}'\rangle} $ maps  $L_ {\langle {\sf e}'\rangle}$ into $(-1,1)$, and $F_{\langle {\sf e}'\rangle}(q_{\langle {\sf e}'\rangle})= F\circ{\sf P}(\tilde{q})=q'$.

Next we consider the quotient $X(\langle {\sf  e}', e \rangle)= \tilde{X}\diagup ({\rm Is}^{\tilde q})^{-1}(\langle {\sf e}',e\rangle)$ whose fundamental group is again a free group in two generators. The image $L_{\langle {\sf e}', e\rangle}\stackrel{def}= \omega_{\langle {\sf e}'\rangle}^{\langle {\sf e}', e\rangle}(L_{\langle {\sf e}'\rangle})$ is a connected component of the preimage of $(-1,1)$ under $ F_{\langle {\sf e}', e\rangle}$. Indeed, $L_{\langle {\sf e}', e\rangle}$ is connected as image of a connected set under a continuous mapping, and $F_{\langle {\sf e}', e\rangle}(\omega_{\langle {\sf e}'\rangle}^{\langle {\sf e}', e\rangle}(L_{\langle {\sf e}'\rangle}))= F\circ \omega_{\langle {\sf e}',e\rangle}\circ
\omega_{\langle {\sf e}'\rangle}^{\langle {\sf e}', e\rangle}(L_{\langle {\sf e}'\rangle})= F_{\langle {\sf e}'\rangle}(L_{\langle {\sf e}'\rangle})\subset (-1,1)$. Moreover, since the mapping $\omega_{\langle {\sf e}'\rangle}^{\langle {\sf e}', e\rangle}:X({\langle {\sf e}'\rangle})\to X({\langle {\sf e}',e\rangle})$ is a covering, its restriction $({\rm Im}F\circ \omega_{\langle {\sf e}'\rangle})^{-1}(0)\to ({\rm Im}F\circ \omega_{\langle {\sf e}', e\rangle})^{-1}(0)$  is a covering. Hence the image under $\omega_{\langle {\sf e}'\rangle}^{\langle {\sf e}', e\rangle}$ of a connected component of $({\rm Im}F\circ \omega_{\langle {\sf e}'\rangle})^{-1}(0)$
is open and closed in $({\rm Im}F\circ \omega_{\langle {\sf e}', e\rangle})^{-1}(0)$. Hence, $L_{\langle {\sf e}', e\rangle}\stackrel{def}= \omega_{\langle {\sf e}'\rangle}^{\langle {\sf e}', e\rangle}(L_{\langle {\sf e}'\rangle})$  is a connected component of the preimage of $(-1,1)$ under $ F_{\langle {\sf e}', e\rangle}$.
%In the same way
Put $q_{\langle {\sf e}', e\rangle}=\omega_{\langle {\sf e}'\rangle}^{\langle {\sf e}', e\rangle}(q_{\langle {\sf e}'\rangle})=\omega_{\langle {\sf e}'\rangle}^{\langle {\sf e}', e\rangle}\circ \omega^{\langle {\sf e}'\rangle}(\tilde{q})=\omega^{\langle {\sf e}', e\rangle}(\tilde{q})$. Note that
%$\omega_{\langle {\sf e}', e\rangle}(q_{\langle {\sf e}', e\rangle})=q$, hence
$F_{\langle {\sf e}', e\rangle}(q_{\langle {\sf e}', e\rangle})=F\circ\omega_{\langle {\sf e}', e\rangle}(q_{\langle {\sf e}', e\rangle})=
F(q)=q'$.

The free homotopy class $\reallywidehat{{\sf e}'_{\langle {\sf e}',e \rangle}}$
%$= \reallywidehat{{\sf e}'_{\langle {\sf e}',e \rangle}}$
in $X(\langle {\sf e}',e \rangle)$
that is related to ${\sf e}'$  intersects $L_{\langle {\sf e}',e \rangle} $. Indeed,
consider any loop $\gamma'_{\langle{\sf e}',e  \rangle}$ in $X(\langle {\sf e}',e \rangle)$ with some base point $q'_{\langle {\sf e}',e \rangle}$, that represents $\reallywidehat{{\sf e}'_{\langle {\sf e}',e \rangle  }}$.
There exists a loop ${\gamma}'_{\langle {\sf e}' \rangle}$ in $ X(\langle {\sf e}' \rangle)$ which represents $\reallywidehat{{\sf e}'_{\langle{\sf e}' \rangle}}$ such that $\omega_{\langle {\sf e}'\rangle}^{\langle {\sf e}',e\rangle  }({\gamma}'_{\langle {\sf e}'\rangle})=\gamma'_{\langle {\sf e}',e\rangle}$.
Such a curve ${\gamma}'_{\langle {\sf e}' \rangle}$ can be obtained as follows. There is a loop $\gamma''_{\langle {\sf e}',e\rangle}$ in $X(\langle {\sf e}',e\rangle  )$    with base point ${q}_{\langle {\sf e}',e\rangle}$ that represents $({\sf e}'  )_{\langle {\sf e}',e\rangle}$,  and a curve $\alpha'_{\langle {\sf e}',e\rangle}$ in $X(\langle{\sf e}',e\rangle)$, such that $\gamma'_{\langle {\sf e}',e\rangle}$ is homotopic with fixed endpoint to $(\alpha'_{\langle {\sf e}',e\rangle})^{-1}\, \gamma''_{\langle{\sf e}',e\rangle}\,    \alpha'_{\langle {\sf e}',e\rangle}$. Consider the lift $\tilde{\gamma}''$ of $\gamma''_{\langle {\sf e}',e\rangle}$ to $\tilde X$ with terminating point $\tilde{q}$, and the lift $\tilde{\alpha}'$ of $\alpha'_{\langle {\sf e}',e\rangle}$ with initial point $\tilde{q}$. The initial point of $\gamma''_{\langle {\sf e}',e\rangle}$ equals $\sigma(\tilde{q})$ for the covering transformation $\sigma= ({\rm Is}^{\tilde{q}})^{-1}({\sf e}' )=
({\rm Is}^{\tilde{q}_0})^{-1}({\sf e}' _0)$. (See equation \eqref{eq1''}.)
The initial point of the curve $\sigma((\tilde{\alpha}')^{-1}) \tilde{\gamma}''\tilde{\alpha}'$ is obtained from its terminating point by applying the covering transformation $\sigma$. Hence,
$\omega^{\langle {\sf e}'\rangle}(\sigma((\tilde{\alpha}')^{-1}) \tilde{\gamma}''\tilde{\alpha}')$ is a closed curve in $X(\langle{\sf e}'  \rangle )$ that represents $\reallywidehat{{\sf e}' _{\langle {\sf e}'\rangle}}$ and projects to $(\alpha'_{\langle {\sf e}',e\rangle})^{-1}\, \gamma''_{\langle {\sf e}',e\rangle}\,    \alpha'_{\langle {\sf e}',e\rangle}$ under $\omega_{\langle {\sf e}' \rangle}^{\langle {\sf e}',e\rangle  }$.
Since $\gamma'_{\langle {\sf e}',e  \rangle}$ is homotopic to $(\alpha'_{\langle {\sf e}',e \rangle})^{-1}\, \gamma''_{\langle {\sf e}',e\rangle}\,    \alpha'_{\langle {\sf e}',e\rangle}$ with fixed base point, it also has a lift to $X(\langle {\sf e}'\rangle)$ which represents $\reallywidehat{{\sf e}' _{\langle {\sf e}'\rangle}}$.

Since $\reallywidehat{{\sf e}'_{\langle  {\sf e}'  \rangle}}$ intersects $L_{\langle {\sf e}' \rangle}$,
the loop
${\gamma}'_{\langle {\sf e}' \rangle}$ has an intersection point $p'_{\langle {\sf e}' \rangle}$ with $L_{\langle {\sf e}'\rangle}$. The point $p'_{\langle {\sf e}',e \rangle  }=  \omega_{\langle {\sf e}'\rangle}^{\langle {\sf e}',e\rangle  }(p'_{\langle {\sf e}'\rangle}    )$ is contained in $\gamma'_{\langle {\sf e}',e \rangle  }$ and in $L_{\langle {\sf e}',e\rangle}$.
We proved that the free homotopy class $\reallywidehat{{\sf e}'_{\langle {\sf e}',e\rangle } }$     in $X(\langle {\sf e}',e \rangle  )$ intersects $L_{\langle {\sf e}',e\rangle  }$.

\noindent {\bf 3.4. A system of generators associated to a standard bouquet of circles.}
We claim that
%, maybe, after changing $\,{\sf e}'\,$ to $\,({\sf e}')^{-1}\;$,
the system of generators  $\;{\sf e}'_{\langle {\sf e}',e\rangle }, \;\; e_{\langle {\sf e}',e\rangle }\;$ of $\;\pi_1(X(\langle {\sf e}',e\rangle),q_{\langle {\sf e}',e\rangle })$ is associated to a standard bouquet of circles for $X(\langle {\sf e}',e\rangle)$. If ${\sf e}'\in \mathcal{E}$ the claim can be obtained as in paragraph 3.2. Suppose ${\sf e}'=e'e''$ for $e', e''\in \mathcal{E}$. Consider the system $\mathcal{E}'$ of generators of $\pi_1(X,q)$ that is obtained from $\mathcal{E}$ by replacing $e'$ by $e'e''$.
If $e'$ and $e''$ correspond to a handle of $X$, then $\mathcal{E}'$
is also associated to a standard bouquet of circles for $X$, see
Figure 4a for the case when $e'$ is represented by an $\alpha$-curve and $e''$ is represented by a $\beta$-curve. The situation when $e'$ is represented by a $\beta$-curve and $e''$ is represented by an $\alpha$-curve is similar. The claim is obtained as in paragraph 3.2.

\begin{figure}[h]
\begin{center}
\includegraphics[width=11cm]{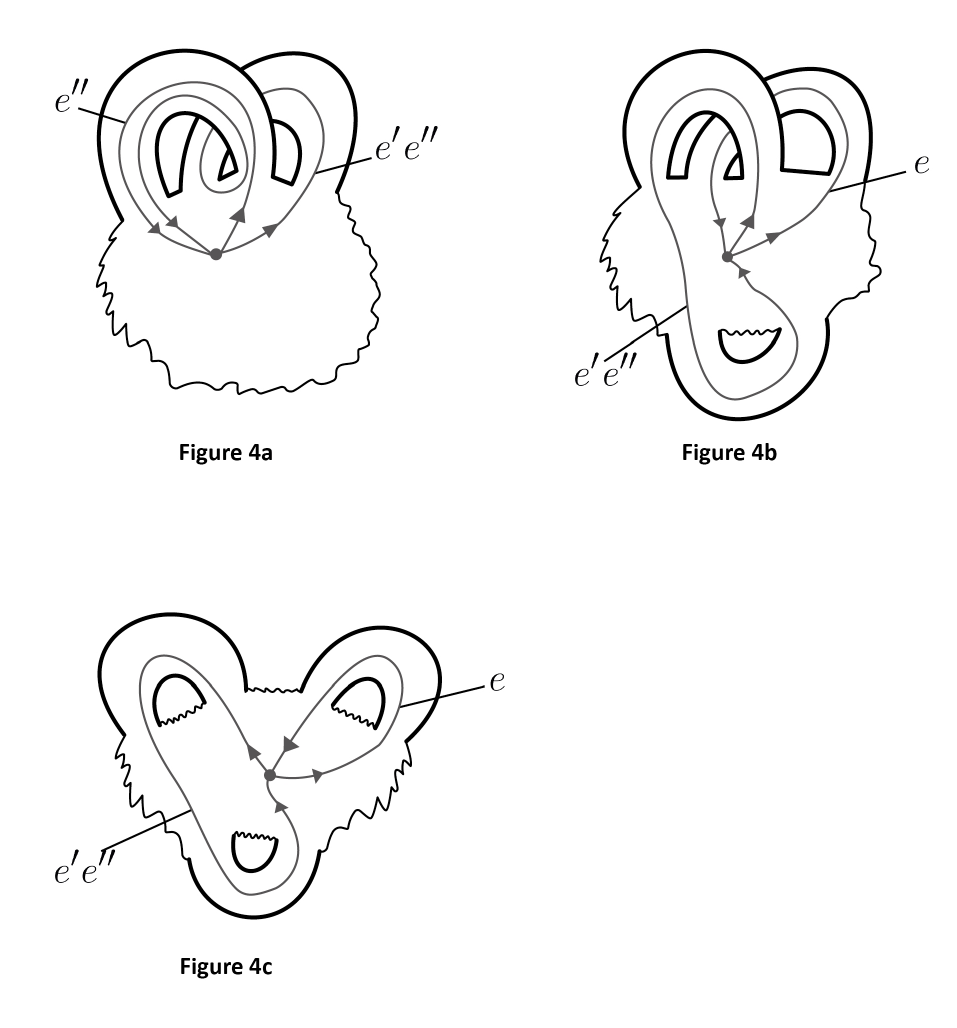}
\end{center}
\caption{Standard bouquets of circles}\label{fig4}
\end{figure}

Suppose one of the pairs $(e,e')$ or $(e,e'')$ corresponds to a handle of $X$.
We assume that $e$ corresponds to an $\alpha$-curve and $e'$ corresponds to a $\beta$-curve
%the pair $(e,e')$ corresponds to
of a handle of $X$
(see Figure 4b). The remaining cases are treated similarly, maybe, after replacing $e'e''$ by $e''e'$ (see paragraph 3.2). With our assumption $\mathcal{E}'$ is associated to a bouquet of circles that is a deformation retract for $X$, but it is not a standard bouquet of circles. Nevertheless, the pair $(e_{\langle {\sf e}',e\rangle }, {\sf e}'_{\langle {\sf e}',e\rangle })$ with ${\sf e}'=e'e''$ is associated to a standard bouquet of circles for $X(\langle {\sf e}',e\rangle)$. This can be seen as before. Consider the bouquet of circles corresponding to $\mathcal{E}'$ and take its union with a disc $D$ around $q$. Lift this set to $\tilde X$. We obtain
the union of a collection of arcs in $\tilde X$ with terminating point $\tilde q$, with a collection of discs in $\tilde X$ around $\tilde q$ and around the initial points of the arcs. Take the union of the arcs and the discs.  The image in $X(\langle {\sf e}',e\rangle)$ of this union under the projection $\omega^{\langle {\sf e}',e\rangle }$ is the union of the two loops $(\gamma_{e})_{\langle {\sf e}',e\rangle}\cup (\gamma_{{\sf e}'})_{\langle {\sf e}',e\rangle}$, the disc $D_{\langle {\sf e}',e\rangle }$ %$\stackrel{def}= \omega^{\langle {\sf e}',e\rangle }(\widetilde{D_0})
%=\omega^{\langle {\sf e}',e\rangle }(\widetilde{D_{e}}), e=e',e''$,
and a set, that is contractible to $D_{\langle {\sf e}',e\rangle }$. Looking at the intersection of the two loops with a small circle contained in  $D_{\langle {\sf e}',e\rangle }$ and surrounding $q_{\langle {\sf e}',e\rangle }$, we see as before that
%,  and contracting the arcs contained in the projection to
%the point $q_{\langle {\sf e}',e\rangle }$.
$(\gamma_{e})_{\langle {\sf e}',e\rangle}\cup (\gamma_{{\sf e'}})_{\langle {\sf e}',e\rangle}$ is a standard bouquet of circles for $X(\langle {\sf e}',e\rangle)$. In this case $X(\langle {\sf e}',e\rangle)$ is a torus with a hole.

In the remaining case no pair of generators among $e$, $e'$, and $e''$ corresponds to a handle.
In this case again $\mathcal{E}'$ does not correspond to a standard bouquet of circles for $X$. But $\{e_{\langle {\sf e}',e\rangle } , (e'e'')_{\langle {\sf e}',e\rangle }\}$ (maybe, after changing $e'e''$ to $e''e'$)
corresponds to a standard bouquet of circles for $X(\langle {\sf e}',e\rangle )$. (See Figure 4c for the case when walking along a small circle around $q$ counterclockwise, we meet the incoming and outgoing rays of representatives of the three elements of $\mathcal{E}$ in the order $e,e',e''$. If the order is different the situation is similar, maybe, after replacing $e'e''$ by $e''e'$.) In this case $X(\langle {\sf e}',e\rangle )$ is a planar domain.

\noindent{\bf 3.5. End of the proof.}
Consider first the case when $X(\langle {\sf e}',e \rangle  )$ is a torus with a hole.
Since $\reallywidehat{{\sf e}'_{\langle {\sf e}',e \rangle }}$ intersects $L_{\langle {\sf e}',e \rangle }$, we see as in the proof when $X$ itself is a torus with a hole, that the curve $L_{\langle {\sf e}',e \rangle }$ cannot be contractible or contractible to the hole, and
the intersection number must be different from zero. Then the intersection number with $L_{\langle {\sf e}',e \rangle }$ of the free homotopy class of one the choices $e_{\langle {\sf e}',e \rangle }^{\pm 1}$ or $({\sf e}'e)_{\langle {\sf e}',e \rangle }$, denoted by  ${\sf e}'''_{\langle {\sf e}',e \rangle }$, is not zero and has the same sign. By lemma \ref{lem3} each of the
$(F_{\langle {\sf e}',e \rangle })_*({\sf e}'_{\langle {\sf e}',e \rangle })$ and $(F_{\langle {\sf e}',e \rangle })_*({\sf e}'''_{\langle {\sf e}',e \rangle })$ is the product of at most two elements of $\pi_1(\mathbb{C}\setminus\{-1,1\},q')$
%$\pi_1(X_{\langle {\sf e}',e \rangle },q_{\langle {\sf e}',e \rangle })$
with $\mathcal{L}_-$ not exceeding
\begin{equation}\label{eq++++}
2\pi\lambda_{{\sf e}'_{\langle {\sf e}',e \rangle },{\sf e}'''_{\langle {\sf e}',e \rangle }}\leq 2\pi \lambda_5(X),
\end{equation}
since ${\sf e}'$ is the product of at most two elements of $\mathcal{E}\cup \mathcal{E}^{-1}$ and ${\sf e}'''$ is the product of at most three elements of $\mathcal{E}\cup \mathcal{E}^{-1}$.
The element $e$ is the product of at most two different elements among the ${\sf e}'$ and ${\sf e}'''$ or their inverses. Hence, the monodromy $F_*(e)=(F _{\langle {\sf e}',e \rangle })_*(e_{\langle {\sf e}',e \rangle })$ is the product of at most four elements with $\mathcal{L}_-$ not exceeding \eqref{eq++++}.
Hence,
\begin{equation}\label{eq+++1}
F_*(e)\leq 8\pi \lambda_5(X)\,.
\end{equation}

Consider now the case when $X(\langle {\sf e}',e \rangle)$ equals $\mathbb{P}^1$ with three holes. Since ${\sf e}'_{\langle {\sf e}',e \rangle}$ and $e_{\langle {\sf e}',e \rangle}$ correspond to a standard bouquet of circles for $X(\langle {\sf e}',e \rangle)$, the curves representing ${\sf e}'_{\langle {\sf e}',e \rangle}$ surround counterclockwise one of the holes, denoted by  $\mathcal{C}'$ , and the curves representing $e_{\langle {\sf e}',e \rangle}$ surround counterclockwise another hole, denoted by  $\mathcal{C}''$.
After applying a M\"obius transformation we may assume that the remaining hole, denoted by $\mathcal{C}_{\infty}$, contains the point $\infty$.
%Denote the third hole by $\mathcal{C}''$.
There are several possibilities for the behaviour of the curve $L_{\langle {\sf e}',e \rangle}$.
%Denote by $\mathcal{C}'$ $\mathcal{C}''$, $\mathcal{C}_{\infty}$....
Since $\reallywidehat{{\sf e}'_{\langle {\sf e}',e \rangle}}$ intersects $L_{\langle {\sf e}',e \rangle}$, the curve $L_{\langle {\sf e}',e \rangle}$ must have limit points on $\mathcal{C}'$. The first possibility is that $L_{\langle {\sf e}',e \rangle}$ has limit points on $\partial\mathcal{C}'$ and $\partial\mathcal{C}''$, the second possibility is, $L_{\langle {\sf e}',e \rangle}$
has limit points on $\mathcal{C}'$ and $\mathcal{C}_{\infty}$, the third possibility is, $L_{\langle {\sf e}',e \rangle}$ has all limit points on $\mathcal{C}'$, and $\mathcal{C}''$ is contained in the bounded connected component of $\mathbb{C}\setminus ( L_{\langle {\sf e}',e \rangle}\cup \mathcal{C}')$.

In the first case the free homotopy classes $\reallywidehat{{\sf e}'_{\langle {\sf e}',e \rangle}}$ and $\reallywidehat{e_{\langle {\sf e}',e \rangle}^{-1}}$ have positive intersection number with the suitably oriented curve $L_{\langle {\sf e}',e \rangle}$. In the second case the free homotopy classes $\reallywidehat{{\sf e}'_{\langle {\sf e}',e \rangle}}$ and $\reallywidehat{ ({{\sf e}'e})_{\langle {\sf e}',e \rangle} }$ have positive intersection number with the suitably oriented curve $L_{\langle {\sf e}',e \rangle}$. In the third case the free homotopy classes of ${\sf e}'_{\langle {\sf e}',e \rangle}$,  $ ({{\sf e}'}^2 e)_{\langle {\sf e}',e \rangle} $ and of their product intersect $L_{\langle {\sf e}',e \rangle}$. The first two cases were treated in paragraph 2 of this section.
The statement concerning the third case is proved as follows.

Any curve that is contained in the complement of $\mathcal{C}'\cup L_{\langle {\sf e}',e \rangle}    $ has either winding number zero around $\mathcal{C}'$ (as a curve in the complex plane $\mathbb{C}$), or its winding number around $\mathcal{C}'$ coincides with the winding number around $\mathcal{C}''$.
On the other hand the representatives of the free homotopy class of ${\sf e}'_{\langle {\sf e}',e \rangle} $
have winding number $1$ around $\mathcal{C}'$ and winding number $0$ around $\mathcal{C}''$.
The representatives of the free homotopy class of $({{\sf e}'}^2 e)_{\langle {\sf e}',e \rangle}$ have winding number $2$ around $\mathcal{C}'$, and  winding number $1$ around $\mathcal{C}''$.
% and winding number zero around each other hole $\mathcal{C}_l,\, l\leq m$.
By the same argument the free homotopy class of the product of ${\sf e}'_{\langle {\sf e}',e \rangle}$ and $({{\sf e}'}^2 e)_{\langle {\sf e}',e \rangle}$ intersects $L_{\langle {\sf e}',e \rangle}$.

We let ${\sf e}'''_{\langle {\sf e}',e \rangle}$ be equal to $e_{\langle {\sf e}',e \rangle}^{-1}$ in the first case, equal to $ ({{\sf e}'e})_{\langle {\sf e}',e \rangle} $ in the second case, and equal to $({{\sf e}'}^2 e)_{\langle {\sf e}',e \rangle}$ in the third case.

By lemma \ref{lem3} each of the
$(F_{\langle {\sf e}',e \rangle })_*({\sf e}'_{\langle {\sf e}',e \rangle })$ and $(F_{\langle {\sf e}',e \rangle })_*({\sf e}'''_{\langle {\sf e}',e \rangle })$ is the product of at most two elements of $\pi_1(\mathbb{C}\setminus\{-1,1\},q')$
%$\pi_1(X_{\langle {\sf e}',e \rangle },q_{\langle {\sf e}',e \rangle })$
with $\mathcal{L}_-$ not exceeding
\begin{equation}\label{eq++++'}
2\pi\lambda_{{\sf e}'_{\langle {\sf e}',e \rangle },{\sf e}'''_{\langle {\sf e}',e \rangle }}\leq 2\pi \lambda_7(X),
\end{equation}
We used that ${\sf e}'$ is the product of at most two elements of $\mathcal{E}\cup \mathcal{E}^{-1}$, $e\in \mathcal{E}\cup \mathcal{E}^{-1}$ and ${\sf e}'''$ is the product of at most five elements of $\mathcal{E}\cup \mathcal{E}^{-1}$.
Since $e$ is the product of at most two different elements among the $({\sf e}')^{\pm 1}$ and $({\sf e}''')^{\pm 1}$, the monodromy $F_*(e)=(F _{\langle {\sf e}',e \rangle })_*(e_{\langle {\sf e}',e \rangle })$ is the product of at most four elements with $\mathcal{L}_-$ not exceeding \eqref{eq++++'}.
Hence,
\begin{equation}\label{eq+++1}
F_*(e)\leq 8\pi \lambda_7(X)\,.
\end{equation}
The proposition is proved.
\hfill $\Box$

\section{($\sf{g},\sf{m}$)-bundles over Riemann surfaces}\label{sec:3}
We will consider bundles whose fibers are smooth surfaces or Riemann surfaces of type
$(\sf{g},\sf{m})$.
\begin{defn}\label{def9.1}{\rm (Smooth oriented $(\textsf{g},\textsf{m})$ fiber bundles.)}
Let $X$ be a smooth oriented manifold of dimension $k$, let
${\mathcal X}$ be a smooth (oriented)
manifold of dimension $k+2$ and ${\mathcal P} : {\mathcal X} \to X$
an orientation preserving smooth proper submersion such that for each point $x \in X$ the
fiber ${\mathcal P}^{-1} (x)$ is a closed oriented surface of genus $\sf{g}$.
Let $\mathbold{E}$ be a smooth submanifold of $\mathcal{X}$ that intersects each fiber $\mathcal{P}^{-1}(x)$ along a set $E_x$ of $\sf{m}$ distinguished points.
Then the tuple ${\mathfrak F}_{\textsf{g},\textsf{m}} = ({\mathcal
X} , {\mathcal P} , \mathbold{E}, X)$ is called
a smooth (oriented) fiber bundle over $X$ with fiber a smooth closed oriented surface of
genus $\textsf{g}$ with $\textsf{m}$ distinguished points
(for short, a smooth oriented $(\textsf{g},\textsf{m})$-bundle).
\end{defn}
If ${\sf m}=0$ the set $\mathbold E$ is the empty set and we will often denote the bundle
by $({\mathcal X} , {\mathcal P} , X)$.
If ${\sf{m}}>0$ the mapping $x\to E_x$ locally defines $\sf m$ smooth sections.
$(\textsf{g},0)$-bundles will also be called genus $\textsf{g}$ fiber bundles.
For $\textsf{g}=1$ and $\textsf{m}=0$ the bundle is also called an elliptic fiber
bundle.  \index{fiber bundle !  elliptic}

\smallskip

In the case when the base manifold is a Riemann surface,  a
holomorphic $(\textsf{g}$,$\textsf{m})$ fiber bundle over $X$
is defined as follows.
\begin{defn}\label{def2}
Let $X$ be a Riemann surface, let ${\mathcal X}$ be a complex surface, and ${\mathcal P}$
a holomorphic proper submersion from ${\mathcal X}$ onto $X$, such that
each fiber $\mathcal{P}^{-1}(x)$ is a closed Riemann surface of genus $\sf g$.
Suppose $\mathbold{E}$ is a complex one-dimensional submanifold of $\mathcal{X}$ that intersects each fiber $\mathcal{P}^{-1}(x)$ along a
set $E_x$ of ${\sf m}$ distinguished points. Then the
tuple ${\mathfrak F}_{\textsf{g},\textsf{m}} = ({\mathcal X} , {\mathcal P}, \mathbold{E}, X)$ is called a
holomorphic $(\textsf{g}$,$\textsf{m})$ fiber bundle over $X$.
\end{defn}
Notice that for ${\sf{m}}>0$ the mapping $x\to E_x$
locally defines $\sf{m}$ holomorphic sections.

Two smooth oriented (holomorphic, respectively) $(\textsf{g},\,\textsf{m})$ fiber
bundles, ${\mathfrak F}^0 = ({\mathcal X}^0 ,
{\mathcal P}^0 , \mathbold{E}^0, X)$ and ${\mathfrak F}^1 = ({\mathcal X}^1 , {\mathcal
P}^1 ,\mathbold{E}^1, X)$, are called smoothly isomorphic (holomorphically isomorphic, respectively,) if there are smooth (holomorphic,
respectively) homeomorphisms $\Phi: \mathcal{X}^0\to \mathcal{X}^1$ and $\phi:X^0\to X^1$
such that for each $x \in X^0$ the mapping $\Phi$ maps the fiber $(\mathcal{P}^0 )^{-1}(x)$
onto the fiber $(\mathcal{P}^1 )^{-1}(\phi(x))$ and the set of distinguished points in
$(\mathcal{P}^0 )^{-1}(x)$ to the set of distinguished points in $(\mathcal{P}^1 )^{-1}(\phi(x))$.
Holomorphically isomorphic bundles will be considered the same holomorphic bundles.

Two smooth (oriented) $(\textsf{g},\textsf{m})$ fiber bundles
over the same oriented
smooth base manifold $X$,  ${\mathfrak F}^0 = ({\mathcal X}^0 ,
{\mathcal P}^0 ,\mathbold{E}^0, X)$, and ${\mathfrak F}^1 = ({\mathcal X}^1 , {\mathcal
P}^1 ,\mathbold{E}^1, X)$, are called (free) isotopic if for an open interval $I$ containing $[0,1]$ there
is a smooth $(\textsf{g},\textsf{m})$ fiber bundle $({\mathcal Y} , {\mathcal P} , \mathbold{E},  X \times
I)$ over the base $X \times I$ (called an isotopy) with the following property.  For each $t \in [0,1]$ we put
${\mathcal Y}^t = {\mathcal P}^{-1} (X \times \{t\})$ and $\mathbold{E}^t= \mathbold{E}\cap {\mathcal P}^{-1} (X \times \{t\})$. The bundle ${\mathfrak F}^0\,$ is
equal to
$\,\left( {\mathcal Y}^0 , {\mathcal P} \mid {\mathcal Y}^0 , \mathbold{E}^0,   X
\times \{0\} \right)\,$, and the bundle $\,{\mathfrak F}^1\,$ is
equal to $\,\bigl( {\mathcal Y}^1 , {\mathcal P} \mid {\mathcal
Y}^1 , \mathbold{E}^1, X \times \{1\} \bigl)\,$.

Two smooth $({\sf g},{\sf m})$-bundles are smoothly isomorphic if and only if they are isotopic (see \cite{Jo5}).

\index{$\mathfrak{F}_{{\sf g},{\sf m}}$} \index{fiber-bundle ! ${\sf g},{\sf m}$-fiber bundle}
Denote
by $S$ a reference surface of genus $\sf g$ with a set $E\subset S$ of $\sf{m}$ distinguished points.
By Ehresmann's Fibration Theorem each smooth $(\sf{g},\sf{m})$-bundle ${\mathfrak
F}_{\textsf{g},\textsf{m}} = ({\mathcal
X} , {\mathcal P} , \mathbold{E}, X)$ with set of distinguished points $E_x \stackrel{def}=\mathbold{E}\cap{\mathcal P}^{-1}(x) $ in the fiber over $x$ is
locally smoothly trivial, i.e. each point in $X$ has a neighbourhood $U\subset X$
such that the restriction of the bundle to $U$ is isomorphic to the trivial bundle $\big(U\times S, {\rm {pr}}_1, U\times E,U\big)$ with set $\{x\}\times E$ of distinguished points in the fiber $\{x\}\times S$ over $x$. Here ${\rm pr}_1: U\times S\to U $ is the projection onto the first factor.

The idea of the proof of Ehresmann's Theorem is the following. Choose smooth coordinates on $U$
by a mapping from a rectangular box to $U$.
Consider
smooth vector fields $v_j$ on $U$, which form a basis of the tangent
space at each point of $U$. Take smooth vector fields $V_j$ on ${\mathcal
P}^{-1} (U)$ that are tangent to $\mathbold{E}$ at points of this set and are mapped to
$v_j$ by the differential of
${\mathcal P}$. Such vector fields can easily be obtained locally.
To obtain the globally defined vector fields $V_j$ on ${\mathcal P}^{-1} (U)$
one uses partitions of unity. The required diffeomorphism $\varphi_U$ is
obtained by composing the flows of these vector fields (in any fixed
order).

In this way a trivialization of the bundle can be obtained over any simply connected domain.

Let $q_0$ be a base point in $X$ and $\gamma(t),\, t \in [0,1],$ a smooth curve in $X$ with
base point $q_0$ that represents an element $e$ of the fundamental group
$\pi_1(X,q_0)$. Let $\varphi^t:\mathcal{P}^{-1}(q_0) \to
\mathcal{P}^{-1}(\gamma(t)),\, t \in [0,1]$, $\varphi^0=\mbox{Id}$, be a smooth family
of diffeomorphisms that map the set of distinguished points  in $\mathcal{P}^{-1}(q_0)$ to
the set of distinguishes points in $\mathcal{P}^{-1}(\gamma(t))$. To obtain such a family we may restrict the bundle to the closed curve given by $\gamma$ and lift the restriction to a bundle over the real axis $\mathbb{R}$. The family of diffeomorphisms may be obtained by considering Ehrenpreis' vector field for the lifted bundle and take the flow of this vector field.
The mapping
$\varphi^1$ obtained for $t=1$ is an orientation preserving self-homeomorphism of the fiber over $q_0$ that
preserves the set of distinguished points. Its isotopy class depends only on the homotopy
class of the curve and the isotopy class of the bundle. The isotopy class of its inverse
$(\varphi^1)^{-1}$ is called the monodromy of the bundle along $e$. Assign to each element $e\in \pi_1(X,q_0)$ the monodromy of the bundle along $e$. We obtain a homomorphism from
$\pi_1(X,q_0)$ to the group of isotopy classes of
self-homeomorphisms of the fiber over $q_0$ that preserve the set of distinguished points. The modular group ${\rm{Mod}}(\textsf{g},\textsf{m})$ is the group of isotopy classes of
self-homeomorphisms of a reference Riemann surface of genus $\textsf{g}$ that map a reference set of $\textsf{m}$
distinguished points to itself.

The following theorem holds (see e.g. \cite{FaMa} and \cite{Jo5} ).

\medskip

\noindent {\bf Theorem G.} {\it Let $X$ be a connected finite open smooth oriented surface.
The set of isotopy classes
of smooth oriented $(\sf{g},\sf{m})$ fiber bundles over $X$
is in one-to-one correspondence to the set of conjugacy classes of
homomorphisms from the fundamental group $\pi_1 (X,q_0)$ into the
modular group ${\rm{Mod}}(\textsf{g},\textsf{m})$.}

\medskip

Let $2{\sf g}-2 + {\sf m}>0$.
A holomorphic $({\sf g},{\sf m})$-bundle is called locally holomorphically trivial if it is locally holomorphically isomorphic to the trivial  $({\sf g},{\sf m})$-bundle. All fibers of a locally holomorphically trivial $({\sf g},{\sf m})$-bundle ${\mathfrak{F}}=({\mathcal{X}}, {\mathcal{P}}, {\mathbold{E}}, {X})$
are conformally equivalent to each other. For a locally trivial holomorphic $({\sf g},{\sf m})$-bundle there exists a finite unramified covering $\hat{\sf P} :\hat{X}\to X$ and a lift
$\hat{\mathfrak{F}}=(\hat{\mathcal{X}}, \hat{\mathcal{P}}, \hat{\mathbold{E}}, \hat{X})$
of $\mathfrak{F}$ to $\hat X$ such that $\hat{\mathfrak{F}}$ is holomorphically isomorphic to the trivial bundle. This can be seen as follows. Consider the lift $\tilde{\mathfrak{F}}$ of the bundle $\mathfrak{F}$ to the universal covering ${\sf P}:\tilde X\to X$ of $X$, i.e.
$\tilde{\mathfrak{F}}=(\tilde{\mathcal{X}}, \tilde{\mathcal{P}}, \tilde{\mathbold{E}}, \tilde{X})$, where the fiber $\tilde{\mathcal{P}}^{-1}(\tilde{x})$ with distinguished points $\tilde{\mathbold{E}}\cap \tilde{\mathcal{P}}^{-1}(\tilde{x})$ is conformally equivalent to the fiber ${\mathcal{P}}^{-1}({x})$ with distinguished points ${\mathbold{E}}\cap {\mathcal{P}}^{-1}({x})$ with $x={\sf P} (\tilde{x})$. The bundle $\tilde{\mathfrak{F}}$ is locally holomorphically trivial. Since $\tilde X$ is simply connected, $\tilde{\mathfrak{F}}$ is holomorphically trivial on $\tilde X$, hence,  there is a biholomorphic mapping $\Phi: \tilde{\mathcal{X}}\to \tilde{X}\times S$ that maps $\tilde{\mathcal{P}}^{-1}(\tilde{x})$ to $\{\tilde{x}\}\times S$ for each $\tilde{x}\in\tilde{X}$. Here $S$ is the fiber $\tilde{\mathcal{P}}^{-1}(\tilde{q}_0)$ over a chosen point $\tilde{q_0}$ over the base point $q_0\in X$. The mapping $\Phi^{-1}$ provides a uniquely determined holomorphic family of conformal mappings $ \varphi_{\tilde{x}}: S=\tilde{\mathcal{P}}^{-1}(\tilde{q}_0) \to \tilde{\mathcal{P}}^{-1}(\tilde{x})$, $\tilde{x}\in \tilde{X}$, that map the set of distinguished points in one fiber to the set of distinguished points in the other fiber, such that the total space ${\mathcal{X}}$ of the bundle $\mathfrak{F}$ is holomorphically equivalent to the quotient of $\{\tilde{x}\}\times S$
by the following equivalence relation $\sim$. A point $\zeta_1$ in the fiber $\tilde{\mathcal{P}}^{-1}(\tilde{x}_1)$ is equivalent by the relation $\sim$ to a point $\zeta_2$ in the fiber $\tilde{\mathcal{P}}^{-1}(\tilde{x}_2)$ if ${\sf P}(\tilde{x}_1)={\sf P}(\tilde{x}_2)$ and $\zeta_2=\varphi_{\tilde{x}_2}\varphi_{\tilde{x}_1}^{-1}(\zeta_1)$. The mapping $\varphi_{\tilde{x}_2}\varphi_{\tilde{x}_1}^{-1}(\zeta_1)$ is a holomorphic self-homeomorphism of the fiber $\tilde{\mathcal{P}}^{-1}(\tilde{x}_1)$. The set of such self-homeomorphims is finite.

Consider the set $N$ of elements $e\in \pi_1(X,q_0)$ for which $\varphi_{({\rm Is}^{\tilde{q}_0})^{-1}(e)(\tilde{q}_0)}$
%$\varphi_{\tilde{x}_0}$
is the identity. As before  $({\rm Is}^{\tilde{q}_0})^{-1}$ is the isomorphism from the fundamental group to the group of covering transformations. The set $N$ is a normal subgroup of the fundamental group. It is of finite index, since two cosets $e_1 \,N$ and $e_2\, N$ are equal if $\varphi_{({\rm Is}^{\tilde{q}_0})^{-1}(e_2 e_1^{-1})(\tilde{q}_0)}={\rm Id}$, and there are only finitely many distinct holomorphic self-homeomorphisms of $\tilde{\mathcal{P}}^{-1}(\tilde{q}_0)$. Hence, $\tilde{X}\diagup ({\rm Is}^{\tilde{q}_0})^{-1}(N)$ is a finite unramified covering of $X$ and the lift of the bundle $\mathfrak{F}$ to $\hat{X}$ has the required property.

Vice versa, if for a holomorphic $({\sf g},{\sf m})$-bundle $\mathfrak{F}$ there exists a finite unramified covering $\hat{\sf P}:\hat{X}\to X$, such that the lift
$\hat{\mathfrak{F}}=(\hat{\mathcal{X}}, \hat{\mathcal{P}}, \hat{\mathbold{E}}, \hat{X})$
of $\mathfrak{F}$ to $\hat X$
%, such that $\hat{\mathfrak{F}}$
is holomorphically isomorphic to the trivial bundle, then $\mathfrak{F}$ is locally holomorphically trivial.

A smooth (holomorphic, respectively) bundle is called isotrivial, if it has a finite covering by the trivial bundle. If all monodromy mapping classes of a smooth bundle are periodic, then the bundle is isotopic (equivalently, smoothly isomorphic) to an
isotrivial bundle. This can be seen by the same arguments as above.

We explain now the notion of irreducible smooth $(\textsf{g},\textsf{m})$-bundles. It is based on Thurston's notion of irreducible surface homeomorphisms.
Let $S$ be a connected finite smooth oriented surface. It is either closed or homeomorphic to a surface with a finite number of punctures. We will assume from the beginning that $S$ is either closed or punctured.

A finite non-empty set of mutually disjoint Jordan curves $\{ C_1 ,
\ldots , C_{\alpha}\}$ on a connected closed or punctured oriented surface $S$ is called
admissible if no $C_i$ is homotopic
to a point in $X$, or to a puncture,
or
to a $C_j$ with $i \ne j$. Thurston calls an isotopy class $\mathfrak{m}$ of self-homeomorphisms
of $S$ (in other words, a mapping class on $S$) reducible if there is an admissible system
of curves $\{ C_1 ,\ldots , C_{\alpha}\}$ on $S$ such that
some (and, hence, each) element in $\mathfrak{m}$ maps the system to an isotopic system. In
this case we say that the system  $\{ C_1 ,
\ldots , C_{\alpha}\}$ reduces $\mathfrak{m}$. A mapping class which is not reducible is
called irreducible.

Let $S$ be a closed or punctured surface with set $E$ of distinguished points. We say that $\varphi$ is a self-homeomorphism of $S$ with distinguished points $E$, if $\varphi$ is a self-homeomorphism of $S$ that maps the set of distinguished points $E$ to itself.
Notice that each self-homeomorphism of the punctured surface $S\setminus E$ extends to a self-homeomorphism of the surface $S$ with set of distinguished points $E$.
We will sometimes identify
self-homeomorphisms of $S\setminus E$ and self-homeomorphism of $S$ with set $E$ of distinguished points.

For a (connected oriented closed or punctured) surface $S$ and a finite subset $E$ of $S$ a finite non-empty set of mutually disjoint Jordan curves $\{ C_1 ,
\ldots , C_{\alpha}\}$ in $S\setminus E$ is called
admissible for $S$ with set of distinguished points $E$ if
it is admissible for $S \setminus E$. An admissible system of
curves for $S$ with set of distinguished points $E$ is said to reduce a mapping class $\mathfrak{m}$ on $S$ with set of distinguished
points $E$, if the induced mapping class on $S\setminus E$ is reduced by this system
of curves.

Conjugacy classes of reducible mapping classes can be decomposed in some sense into
irreducible components, and conjugacy classes of reducible mapping classes can be recovered
from the irreducible components up to products of commuting Dehn twists. Conjugacy classes
of irreducible mapping classes are classified and studied.

A Dehn twist \index{Dehn twist} about a simple closed curve $\gamma$
in an oriented surface $S$ is a mapping that is isotopic to the following one. Take a tubular neighbourhood of
$\gamma$ and parameterize it as a round annulus  $A=\{ e^{-\varepsilon} < \vert z \vert < 1\}$  so that $\gamma$ corresponds to $|z|=e^{-\frac{\varepsilon}{2}}$. The mapping is
an orientation preserving self-homeomorphism of $S$
which is the identity outside $A$
and is equal to the mapping $e^{-\varepsilon s +2\pi i t}\to  e^{-\varepsilon s +2\pi i (t+s) }$ for  $e^{-\varepsilon s +2\pi i t}\in A$, i.e. $s\in (0,1)$. Here $\varepsilon$ is a small positive number.

Thurston's notion of reducible mapping classes takes over to families of mapping classes on a surface of type $({\sf{g}},{\sf{m}})$, and therefore to $({\sf{g}},{\sf{m}})$-bundles. Namely,
an admissible system of
curves on a (connected oriented closed or punctured) surface $S$
with set of $\sf m$ distinguished points $E$ is said to reduce a family of mapping classes $\mathfrak{m}_j \in \mathfrak{M}(S; \emptyset, E) $ if it reduces each $\mathfrak{m}_j$.
Similarly, a $({\sf{g}},{\sf{m}})$-bundle with fiber $S$ over the base point $q_0$ and set of distinguished points $E\subset S$ is called reducible if there is an admissible system of curves in the fiber over the base point
%it is isotopic to a surface bundle with an admissible system of
%curves in the fiber over the base point
that reduces all monodromy mapping classes simultaneously. Otherwise the bundle is called
irreducible.

Reducible bundles can be decomposed into irreducible bundle components and the reducible bundle can be recovered from the irreducible bundle components up to commuting Dehn twists
in the fiber over the base point.

Let $X$ be a finite open connected Riemann surface. By a holomorphic (smooth, respectively) $(0,n)$-bundle with a section over $X$ we mean a holomorphic (smooth, respectively) $(0,n+1)$-bundle
$(\mathcal{X},\mathcal{P}, \mathbold{E}, X)$,
such that the complex manifold (smooth manifold, respectively) $\mathbold{E}\subset \mathcal{X}$ is the disjoint union of two complex manifolds (smooth manifolds, respectively) $\mathring{\mathbold{E}}$ and $\mathbold{s}$, where $\mathring{\mathbold{E}}\subset \mathcal{X}$
intersects each fiber $\mathcal{P}^{-1}(x)$ along a set $\mathring{E}_x$ of $n$ points, and $\mathbold{s}\subset \mathcal{X}$ intersects each fiber $\mathcal{P}^{-1}(x)$ along a
single point $s_x$. We will also say, that the mapping $x\to s_x,\, x\in X$, is a holomorphic (smooth, respectively) section of the $(0,n)$-bundle with set of distinguished points $\mathring{E}_x$ in the fiber over $x$.

A special $(0,n+1)$-bundle is a bundle over $X$ of the form $(X\times \mathbb{P}^1, {\rm pr}_1, \mathbold{E}, X)$,
where ${\rm pr}_1 : X\times \mathbb{P}^1 \to X$ is the projection onto the first factor, and the smooth submanifold $\mathbold{E}$ of $X\times \mathbb{P}^1$  is equal to the disjoint union $\mathring{\mathbold{E}}\cup \mathbold{s}_{\infty}$ where $\mathbold{s}_{\infty}$ intersects each fiber $\{x\}\times\mathbb{P}^1$ along the point $\{x\}\times\{\infty\}$,
and the set $\mathring{\mathbold{E}}$ intersects each fiber along $n$ points.
A special $(0,n+1)$-bundle is, in particular, a $(0,n)$-bundle with a section.

Two smooth $(0,n)$-bundles with a section (in particular, two special $(0,n+1)$-bundles) are called isotopic if they are isotopic as $(0,n+1)$-bundles with an isotopy that joins the sections of the bundles. A holomorphic
(smooth, respectively) $(0,n)$-bundle with a section is isotopic to a
holomorphic (smooth, respectively) special  $(0,n+1)$-bundle over $X$ (see \cite{Jo5}).

Theorem \ref{thm2} is a consequence of the following theorem on $(0,3)$-bundles with a section.

\begin{thm}\label{thm3}
Over a connected Riemann surface of genus $g$ with $m+1$ holes there are up to isotopy no more than $(15 \exp( 6 \pi \lambda_{10}(X)))^{6(2g+m)}$
%$\big(3^6\cdot 5^6\cdot\exp(36 \pi \lambda_8(X))\big)^{2g+m}$
irreducible holomorphic $(0,3)$-bundles
with a holomorphic section.
\end{thm}

\smallskip

For a reducible $(0,4)$-bundle the fiber of each irreducible bundle component is a thrice-punctured Riemann sphere. Hence each irreducible bundle component of a reducible $(0,4)$-bundle is isotopic to an isotrivial bundle. For more details see \cite{Jo5}.

\smallskip

Theorem \ref{thm1} (with a weaker estimate) is a consequence of Theorem \ref{thm3}. Indeed, consider holomorphic (smooth, respectively) bundles whose fiber over each point $x\in X$ equals $\mathbb{P}^1$ with set of distinguished points $\{-1,1,f(x),\infty\}$ for a function $f$ which depends holomorphically (smoothly, respectively) on the points $x\in X$ and does not take the values $-1$ and $1$. Then we are in the situation of Theorem \ref{thm1}.
It is not hard to see that the mapping $f$ is reducible, iff this $(0,4)$-bundle
%whose fibers are $\mathbb{C}$ with set of distinguished points $F(z)$
is reducible (see also Lemma 7 of \cite{Jo5}).

The relation between Theorems \ref{thm2} and  \ref{thm3} is given in Proposition \ref{prop1} below.

A holomorphic $(1,1)$-bundle $\mathfrak{F}=\big(\mathcal{X},\mathcal{P}, \mathbold{s},X\big)$ is called a double branched covering of the special
holomorphic %(smooth, respectively)
$(0,4)$-bundle $\;\big(X\times \mathbb{P}^1,{\rm pr}_1,\mathbold{E} ,\,X\big)\;$
if there exists a holomorphic %(smooth, respectively)
mapping $\;P:\mathcal{X}\to X\times\mathbb{P}^1$
that maps each fiber $\mathcal{P}^{-1}(x)$ of the $(1,1)$-bundle
onto the fiber
$\{x\}\times \mathbb{P}^1$
of the $(0,4)$-bundle over the same point $x$, such that the
restriction
${{P}}: \mathcal{P}^{-1}(x) \to \{x\}\times \mathbb{P}^1$ is a holomorphic double
branched covering with branch locus being the set $\{x\}\times (\mathring{E}_x\cup \{\infty\})=\mathbold{E}\cap (\{x\}\times \mathbb{P}^1) $ of distinguished
points in the fiber $\{x\}\times \mathbb{P}^1$, and ${P}$
maps the distinguished point $s_x$  in the fiber $\mathcal{P}^{-1}(x)$ over $x$ to
the point $\{x\}\times \{\infty\}$ in $\{x\}\times \mathbb{P}^1$. We will also denote
$(X\times \mathbb{P}^1  , {\rm pr}_1  ,\mathbold{E}, X)$ by ${{P}}((\mathcal{X},\mathcal{P},\mathbold{s},X))$, and call
the bundle $(\mathcal{X},\mathcal{P},\mathbold{s}, X)$ a lift of $(X\times\mathbb{P}^1, {\rm pr}_1 ,\mathbold{E} ,X)$. Let the fiber of the $(1,1)$-bundle over the base point $q_0\in X$ be
$Y$ with distinguished point $s$, and let  the fiber of the $(0,4)$-bundle over $q_0$ be
$\mathbb{P}^1$ with distinguished points $\mathring{E}\cup\{\infty\}$ for a set $\mathring{E}\subset C_3(\mathbb{C})\diagup \mathcal{S}_3$. Then the monodromy mapping class $\mathfrak{m}_1\in \mathfrak{M}(\mathbb{P}^1;\infty,\mathring{E})$ of the $(0,4)$-bundle along any generator
of the fundamental group of $X$ is the projection of the monodromy mapping class $\mathfrak{m}\in \mathfrak{M}(Y;s,\emptyset)$ of the $(1,1)$-bundle along the same generator. This means that there are representing homeomorphisms $\varphi\in \mathfrak{m}$ and $\varphi_1\in \mathfrak{m}_1$ such that $\varphi_1({ P}(\zeta))= {P}(\varphi(\zeta)),\, \zeta \in Y$. We will also say that $\mathfrak{m}$ is a lift of $\mathfrak{m}_1$. The lifts of a mapping class $\mathfrak{m}_1\in \mathfrak{M}(\mathbb{P}^1;\infty,\mathring{E})$ differ by the involution of $Y$, that interchanges the sheets of the double branched covering.
Hence, each class $\mathfrak{m}_1\in \mathfrak{M}(\mathbb{P}^1;\infty,\mathring{E})$ has exactly two lifts.

\begin{prop}\label{prop1} Let $X$ be a Riemann surface
of genus $ g$ with ${ m}+1\geq 1$ holes with base point $q_0$ and curves $\gamma_j$ representing a set of generators $e_j\in \pi_1(X,q_0)$.\\
\noindent (1) Each holomorphic $(1,1)$-bundle
over $X$ is holomorphically isomorphic to the double branched covering of a special holomorphic
$(0,4)$-bundle over $X$.\\
\noindent (2) Vice versa, for each special holomorphic
$(0,4)$-bundle over $X$ and each collection $\mathfrak{m}^j$ of lifts of the $2{ g} +{ m} $ monodromy mapping classes $\mathfrak{m}_1^j$ of the bundle along the $\gamma_j$ there exists a double branched covering by
a holomorphic $(1,1)$-bundle
with collection of monodromy mapping classes equal to the $\mathfrak{m}^j$.
Each special holomorphic $(0,4)$-bundle has exactly $2^{2g+m}$ non-isotopic holomorphic lifts.\\
\noindent (3) A lift of a special $(0,4)$-bundle is reducible if and only if the special
 $(0,4)$-bundle is reducible.
\end{prop}

The proof of the proposition uses the fact that a holomorphic $(1,1)$-bundle over $X$ is holomorphically isomorphic to a holomorphic bundle whose fiber over each point $x$ is a
quotient $\mathbb{C}\diagup \Lambda_x$ of the complex plane by a lattice $\Lambda_x$ with distinguished point $0\diagup \Lambda_x$. The lattices depend holomorphically on the point $x$.
To represent the fibers as branched coverings depending holomorphically on the points in $X$ we use embeddings of punctured tori into $\mathbb{C}^2$ by suitable versions of the Weierstraß $\wp$-function.
For a detailed proof of Proposition \ref{prop1} see \cite{Jo5}.

\smallskip

\noindent {\bf Preparation of the proof of Theorem \ref{thm3}.}
The proof will be given in terms of braids. Let $C_n(\mathbb{C})=\{(z_1,\ldots,z_n)\in \mathbb{C}^n: z_j\neq z_k \;\mbox{for}\; j\neq k\}$ be the $n$-dimensional configuration space.
The symmetrized configuration space is its quotient $C_n(\mathbb{C})\diagup \mathcal{S}_n$ by the diagonal action of the symmetric group $\mathcal{S}_n$.
We write points of $C_n(\mathbb{C})$ as ordered $n$-tuples $(z_1,\ldots,z_n)$ of points in $\mathbb{C}$, and points of $C_n(\mathbb{C})\diagup \mathcal{S}_n$ as unordered tuples $\{z_1,\ldots,z_n\}$ of points in $\mathbb{C}$.
We regard geometric braids on $n$ strands with
base point $E_n$ as loops in the symmetrized configuration space $C_n(\mathbb{C})\diagup \mathcal{S}_n$ with base point $E_n$, and
braids on $n$ strands ($n$-braids for short)
with base point $E_n\in C_n(\mathbb{C})\diagup \mathcal{S}_n$ as homotopy classes of loops with
base point $E_n$ in $C_n(\mathbb{C})\diagup \mathcal{S}_n$,
equivalently, as
elements of the fundamental group $\pi_1(C_n (\mathbb {C})
\diagup
{\mathcal S}_n, E_n)$ of the symmetrized configuration space
with
base point $E_n$.

Each smooth mapping $\;\;F:X \to C_n ({\mathbb C}) \diagup
{\mathcal S}_n\;\;$ defines a smooth special $\;(0,n+1)$-bundle $\;\;\;(X\times \mathbb{P}^1, {\rm pr}_1, \mathbold{E},X)\,,$ where
$\mathbold{E}\cap (\{x\}\times \mathbb{P}^1)=\{x\}\times  ( F(z) \cup \{\infty\})$.
Vice versa, for each smooth special $(0,n+1)$-bundle $(X\times \mathbb{P}^1, {\rm pr}_1, \mathbold{E},X)$ the mapping that assigns to each point $x \in X$ the set of finite distinguished points in the fiber over $x$ defines a smooth  mapping $F:X \to C_n ({\mathbb C}) \diagup{\mathcal S}_n$. The mapping $F$ is holomorphic iff the bundle is holomorphic. It is called irreducible iff the bundle is irreducible.
Choose a base point $q_0\in X$.
The restriction of the mapping $F$ to each loop with base point $q_0$ defines a geometric braid with base point $F(q_0)$. The braid represented by it is called the monodromy of the mapping $F$
along the element of the fundamental group represented by the loop.

The monodromy mapping classes of a special $(0,n+1)$-bundle are isotopy classes of self-homeomorphisms of
the fiber $\mathbb{P}^1$ over the base point $q_0$  which map the set of finite distinguished points $E_n=F(q_0)$ in this fiber onto itself, and fix $\infty$.
Two smooth mappings $F_1$ and $F_2$ from $X$ to $C_n ({\mathbb C}) \diagup {\mathcal S}_n$, that have equal value $E_n\in C_n(\mathbb{C})\diagup\mathcal{S}_n$ at the base point $q_0$, define special $(0,n+1)$-bundles , that are isotopic with an isotopy that fixes the fiber over $q_0$ and the set of distinguished points in this fiber,
iff their restrictions to each curve in $X$ with base point $q_0$ define braids that differ by an element of the center $\mathcal{Z}_n$ of the braid group $\mathcal{B}_n$ on $n$ strands (in other words, by a power of a full twist).
Indeed, the braid group on $n$ strands modulo its center $\mathcal{B}_n\diagup \mathcal{Z}_n$ is isomorphic to the group of mapping classes
of $\mathbb{P}^1$ that fix $\infty$ and map $E_n$ to itself.

Note that for the group $\mathcal{PB}_3$ of pure braids on three strands the quotient $\mathcal{PB}_3\diagup \mathcal{Z}_3$ is isomorphic to the fundamental group of $\mathbb{C}\setminus \{-1,1\}$. The isomorphism maps the generators $\sigma_j^2\diagup \langle \Delta_3^2 \rangle, \, j=1,2$, of $\mathcal{PB}_3\diagup \mathcal{Z}_3$ to the standard generators $a_j,\, j=1,2$, of the fundamental group $\pi_1(\mathbb{C}\setminus \{-1,1\},0)$.
Here $\langle \Delta_3^2 \rangle$ denotes the group generated by $\Delta_3^2$ which is equal to the center $\mathcal{Z}_3$.

The proof of Theorem \ref{thm3} will go now along the same lines as the proof of Theorem \ref{thm1} with some modifications.
Lemma H, Lemmas \ref{lem4} and \ref{lem3a}, and Theorem I below are given in terms of braids rather than in terms of elements of $\mathcal{B}_3\diagup \mathcal{Z}_3$.

The following lemma and the following theorem were proved in \cite{Jo2}.

\medskip

\noindent {\bf Lemma H.} {\it Any braid $b\in \mathcal{B}_3$ which is not a power of $\Delta_3$ can be written in a unique way in the form
\begin{equation}\label{eq2'}
\sigma_j^k \, b_1 \, \Delta_3^{\ell}\,
\end{equation}
where $j=1$ or $j=2$, $k\neq 0$ is an integer, $\ell$ is a (not necessarily even) integer, and $b_1$ is a word in $\sigma_1^2$ and $\sigma_2^2$ in reduced form. If
$b_1$ is not the identity, then the first term of $b_1$ is a non-zero even
power of $\sigma_2$ if $j=1$, and  $b_1$ is a non-zero even  power of
$\sigma_1$ if $j=2$.}

\medskip

For an integer $j\neq 0$ we put $q(j)=j$ if $j$ is even, and for odd $j$ we
denote by $q(j)$ the even integer neighbour of $j$ that is closest to zero. In other words, $q(j)=j$ if $j\neq 0$ is even, and for each odd integer $j\,,$  $q(j)= j
-{\mbox{sgn}}(j)$, where  ${\mbox{sgn}}(j)$
for a non-zero integral number $j$ equals $1$ if $j$ is positive,
and $-1$ if $j$ is negative. For a braid in form \eqref{eq2'} we put
$\vartheta(b) \stackrel{def}{=}\sigma_j^{q(k)} \, b_1$. If $b$ is a power of $\Delta_3$ we put $\vartheta(b) \stackrel{def}{=}{\rm Id}$.

Let $C_n(\mathbb{R})\diagup \mathcal{S}_n$ be the totally real subspace of
$C_n(\mathbb{C})\diagup \mathcal{S}_n$. It is defined in the same way as
$C_n(\mathbb{C})\diagup \mathcal{S}_n$  by replacing $\mathbb{C}$ by $\mathbb{R}$. Take a base point $E_n \in C_n(\mathbb{R})\diagup \mathcal{S}_n$.
The fundamental group $\;\pi_1(\,C_n (\mathbb {C}) \diagup
{\mathcal
S}_n\,,\; E_n\,)\;$ with base point is isomorphic to the relative fundamental
group
$\;
\pi_1(\,C_n (\mathbb {C}) \diagup {\mathcal S}_n\,,\; C_n
(\mathbb
{R}) \diagup {\mathcal S}_n \,)\,.
$
The elements of the latter group
are homotopy classes of arcs in $\,C_n (\mathbb {C}) \diagup
{\mathcal S}_n\,$ with endpoints in the totally real subspace
$\,C_n (\mathbb {R})\diagup{\mathcal S}_n\,$ of the symmetrized configuration space.

Let $b \in \mathcal{B}_n$ be a braid.
Denote by $b_{tr}$ the element
of the relative fundamental group
$\pi_1(\,C_n
(\mathbb {C}) \diagup {\mathcal S}_n\,,\; C_n (\mathbb {R})
\diagup
{\mathcal S}_n\, )$ that corresponds to $b$ under the mentioned group isomorphism.
For a rectangle $R$ in the plane with sides parallel to the axes we let $f:R
\to \,C_n
(\mathbb {C}) \diagup {\mathcal S}_n\,$ be a mapping which
admits a continuous extension to the closure $\bar R$ (denoted again by $f$) which
maps the (open) horizontal sides into $\,C_n (\mathbb {R})
\diagup
{\mathcal S}_n\,$. We say that the mapping represents $b_{tr}$
if for each maximal vertical line segment contained in $R$
(i.e. $R$
intersected with a vertical line in $\mathbb{C}$) the
restriction of
$f$ to the closure of the line segment represents $b_{tr}$.

The  extremal length of a $3$-braid with totally real horizontal boundary values is defined as
\begin{align}
\Lambda_{tr}(b)=& \inf \{\lambda(R): R\, \mbox{ a rectangle
which
admits a holomorphic map to} \nonumber \\
&C_n (\mathbb {C}) \diagup {\mathcal S}_n \,\mbox{ that
represents}\; b_{tr}\}\,.\nonumber
\end{align}
(see \cite{Jo2}.)
The following theorem holds (see \cite{Jo2}).

\medskip

\noindent {\bf Theorem I.} {\it Let $b \in \mathcal{B}_3$ be a (not necessarily pure) braid which is not a power of $\Delta_3$, and let $w$ be the reduced word representing the image of $\vartheta(b)$ in $\mathcal{PB}_3 \diagup \langle \Delta_3^2\rangle$.
Then
$$
\Lambda_{tr}(b) \geq  \frac{1}{2\pi}\cdot \mathcal{L}_-(w)  \,,
$$
except in the case when $b=\sigma_j^{k}\,\Delta_3^{\ell}$, where $j=1$ or  $j=2$, $k\neq 0$ is an integer number, and $\ell$ is an arbitrary integer. In this case $\Lambda_{tr}(b)=0$.}

\medskip

The set
\begin{align}
\mathcal{H} \stackrel{def}= & \{ \{z_1,z_2,z_3\} \in C_3(\mathbb{C})\diagup \mathcal{S}_3: \mbox{the three points}\; z_1,z_2, z_3  \nonumber\\
& \mbox{are contained in a real line in the complex plane}\}
\end{align}
is a smooth real hypersurface of $C_3(\mathbb{C})\diagup \mathcal{S}_3$. Indeed, let $\{z_1^0,z_2^0,z_3^0\}$ be a point of the symmetrized configuration space. Introduce coordinates near this point by lifting a neighbourhood of the point to $C_3(\mathbb{C})$ with coordinates $(z_1,z_2,z_3)$. Since the linear map  $M(z)\stackrel{def}= \frac{z-z_1}{z_3-z_1},\, z \in \mathbb{C},$ maps the points $z_1$ and $z_3$ to the real axis, the three points $z_1, \, z_2,$ and $z_3$  lie on a real line in the complex plane iff the imaginary part of $z'_2\stackrel{def}=M(z_2)= \frac{z_2-z_1}{z_3-z_1}$ vanishes. The equation $\mbox{Im}\frac{z_2-z_1}{z_3-z_1}=0$ in local coordinates $(z_1,z_2,z_3)$ defines a local piece of a smooth real hypersurface.

For each complex affine self-mapping $M$ of the complex plane we consider the diagonal action $M\big((z_1,z_2,z_3)\big)=\big(M(z_1),M(z_2),M(z_3)\big)$ on points $(z_1,z_2,z_3)\in C_3(\mathbb{C})$, and the diagonal action $M\big(\{z_1,z_2,z_3\}\big)=\{M(z_1),M(z_2),M(z_3)\}$ on points $\{z_1,z_2,z_3\}\in C_3(\mathbb{C})\diagup \mathcal{S}_3$.

The following two lemmas replace Lemma \ref{lem2} in the case of
$(0,3)$-bundles with a section.

\begin{lemm}\label{lem4}
Let $A$ be an annulus with an orientation of simple closed dividing curves. Suppose $F:A \to C_3(\mathbb{C})\diagup \mathcal{S}_3$ is a holomorphic mapping whose image is not contained in $\mathcal{H}$.
Suppose $L_A$ is a simple relatively closed curve in $A$ with limit points on both boundary circles of $A$, and $F(L_A) \subset \mathcal{H}$. Moreover, for a point $q_A\in L_A$ the value $F(q_A)$ is in the totally real subspace $ C_3(\mathbb{R})\diagup \mathcal{S}_3$.
Let $e_A\in \pi_1(A,q_A)$ be the positively oriented generator of the fundamental group of $A$ with base point $q_A$.
If the braid $b\stackrel{def}=F_*(e_A) \in \mathcal{B}_3 $ is different from $\sigma_j^k \, \Delta_3^{2 \ell'} $ with
$j$ equal to $1$ or $2$, and $k\neq 0$ and $\ell'$ being integers, then
\begin{equation}\label{eq+}
\mathcal{L}_-(\vartheta(b)) \leq 2\pi \lambda(A).
\end{equation}
\end{lemm}
Notice that the braids  $\sigma_j^k \, \Delta_3^{\ell} $ for odd $\ell$ are exceptional for Theorem I, but not exceptional for Lemma \ref{lem4}. The reason is that the braid in Lemma \ref{lem4} is related to a mapping of an annulus, not merely to a mapping of a rectangle.
For $t\in [0,\infty)$ we put
$\log_+ t\stackrel{def}=
\begin{cases} \log t& \; t \in [1,\infty) \\
0 & \; t \in [0,1)\;\;\;\; .\\
\end{cases}$

\begin{lemm}\label{lem3a}
If the braid in Lemma {\rm \ref{lem4}} equals $b = \sigma_j^k \, \sigma_{j'}^{k'} \,\Delta_3^{\ell}$ with $j$ and $j'$ equal to $1$ or to $2$,  $j'\neq j$, and $k$ and $k'$ being non-zero integers,
and $\ell$ an even integer, then
\begin{equation}\label{eq++}
\log_+(3[\frac{|k|}{2}]) + \log_+(3[\frac{|k'|}{2}])\leq {\pi} \lambda(A).
\end{equation}
\end{lemm}
\noindent Here for a non-negative number $x$ we denote by $[x]$ the smallest integer not exceeding $x$.

\medskip

\noindent {\bf Proof of Lemma \ref{lem4}.}
By the same argument as in the proof of Lemma  \ref{lem2} we may assume that
the annulus $A$ has smooth boundary, the mapping
$F$ extends continuously to the closure $\overline{A}$, and the curve $L_A$ is a smooth (connected) curve in $\overline{A}$ whose endpoints are on different boundary components of $A$. Perhaps after a diagonal action of a fixed M\"obius transformation on each point of $C_3(\mathbb{C}\diagup \mathcal{S}_3)$,
we may also assume that the value of $F$ at the point $q_A \in L_A$ is equal to the unordered triple
$\{-1,q',1\}\in C_3(\mathbb{R})\diagup \mathcal{S}_3$ for a number $q' \in \mathbb{R}\setminus\{-1,1\}$.
We restrict the mapping $F$ to
$A\setminus L_0$. Let $\tilde R$ be a lift of $A\setminus L_A$ to the infinite strip $\widetilde{\overline{A}}$ that covers $\overline{A}$. We consider $\tilde R$ as curvilinear rectangle with horizontal sides being the two different lifts of $L_A$
and vertical sides being the lifts of the two boundary circles cut at the endpoints of $L_A$.

Take a closed curve $\gamma_A:[0,1]\to A$ in $A$ with base point $q_A\in L_A$, that intersects $L_A$  only at the base point and represents
the element $e_A\in\pi_1(A,q_A)$.
Let $\tilde{\gamma}_A$ be the lift of $\gamma_A$ to $\widetilde{\overline{A}}$ for which  $\tilde{\gamma}_A((0,1))$ is contained in $\tilde R$,
and let $\tilde F=(\tilde {F}_1,\tilde {F}_2,\tilde {F}_3):{\tilde R} \to \mathcal{C}_3(\mathbb{C})$ be a lift of $F$
to a mapping  from ${\tilde R}$ to the configuration space $\mathcal{C}_3(\mathbb{C})$.
The continuous extension of $\tilde F$ to $\overline{\tilde R}$ is also denoted by $\tilde F$.
We may choose the lift so that the value
of $\tilde F$
at the copy of $q_A$ on the lower horizontal side of $\tilde R$
equals $(-1,q',1)$. For each $z\in \overline{\tilde R}$
we consider the complex affine mapping  $\mathfrak{A}_z(\zeta)=a(z)\zeta+b(z) \stackrel{def} = -1 + 2\frac{\zeta-\tilde{F}_1(z)}{\tilde{F}_3(z)-\tilde{F}_1(z)},\, \zeta \in \mathbb{C}$.  Denote by
$\hat{F}(z)\stackrel{def}=\mathfrak{A}_{z}(\tilde{F}(z))=\mathfrak{A}_{z}
\big((\tilde{F}_1(z),\tilde{F}_2(z),\tilde{F}_3(z)\big),\, z \in \overline{\tilde{R}}, $
the result of applying  $\mathfrak{A}_z$ to each of the three points of $\tilde{F}(z)$.
The mapping $\hat{F}(z)=(\hat{F}_1(z),\hat{F}_2(z),\hat{F}_3(z))=(-1,\hat{F}_2(z),1)$
is holomorphic on $\tilde R$.
Let $\tilde{F}_{\rm sym}\stackrel{def}=\{\tilde{F}_1,\tilde{F}_2,\tilde{F}_3\}$ ( $\hat{F}_{\rm sym}\stackrel{def}=\{\hat{F}_1,\hat{F}_2,\hat{F}_3\}$, respectively) be the projection of $\tilde{F}$ ($\hat{F}$, respectively) to a mapping from $\overline{\tilde{R}}$ to the symmetrized configuration space $C_3(\mathbb{C})\diagup \mathcal{C}_3$.
Since $F(L_A)\subset \mathcal{H}$ the mapping $\hat{F}_{\rm sym}$
takes the horizontal sides of $\tilde R$
to the totally real subspace
$C_3(\mathbb{R})\diagup \mathcal{S}_3$ of the symmetrized configuration space. Moreover, $\hat{F}_{\rm sym}$ maps the copy of $q_A$ on the
lower side of $\tilde R$ to $\{-1,q',1\}$.
%Let $\tilde{\hat{F}}$ be the lift of $\hat F$ to the configuration space which
Recall that also $\tilde{F}_{\rm sym}$
takes the value $\{-1,q',1\}$ at the copy of $q_A$ on the
lower side of $\tilde R$.

The restrictions of $\tilde{F}_{\rm sym}$ and of
${\hat{F}}_{\rm sym}$ to the curve $\tilde{\gamma}_A$
represent elements of the relative fundamental group
$\pi_1(C_3(\mathbb{C})\diagup \mathcal{S}_3,C_3(\mathbb{R})\diagup \mathcal{S}_3)$.
The represented elements of the relative fundamental group
differ by a finite number of half-twists.
Indeed, for each $z$, the lifts to $C_3(\mathbb{C})$, $\tilde F(z)$ and ${\hat{F}}(z)$, differ by a complex affine mapping.
Hence, ${\hat{F}}(\tilde {\gamma}_A(t))=b(t)+a(t) \tilde F(\tilde{\gamma}_A(t))$ for continuous functions $a$ and $b$ on $[0,1]$ with $a$ nowhere vanishing, $b(0)=0$, $a(0)=1$, and $b(1)$ and $a(1)$ real valued.
Then the function $b:[0,1]\to \mathbb{C}$ is homotopic with endpoints in $\mathbb{R}$ to the function that is identically equal to zero. The mapping $a:[0,1]\to \mathbb{C}\setminus \{0\}$ is homotopic with endpoints in $\mathbb{R}$ to $\frac{a}{|a|}$. Hence, the mappings $\hat{F}(\tilde{\gamma}_A(t))$ and $\frac{a(t)}{|a|(t)}\tilde{F}(\tilde{\gamma}_A(t))$ from $[0,1]$ to $C_3(\mathbb{C})\diagup \mathcal{S}_3$ are homotopic with endpoints in $C_3(\mathbb{R})\diagup \mathcal{S}_3$. The statements follows.

Let $\omega(z):A\setminus L_A\to R$ be the conformal mapping of the curvilinear rectangle onto the rectangle of the form $R=\{z=x+i y:x \in (0,1),\, y \in (0,{\sf{a}})\}$, that maps the lower curvilinear side of $A\setminus L_0$ to the lower side of $R$. (Note that the number $\sf{a}$ is uniquely defined by $\tilde R$.) For $i'\in \mathbb{Z}$ we put $ \mathring{F}_{i'}(z) \stackrel{def}= e^{-i' \frac{\pi}{\sf a}\omega(z)} \hat{F}_{\rm sym}(z)$.
Then, for some choice of $i'$ the restrictions $\tilde{F}_{\rm sym}\mid {\gamma}_A$
and  $\mathring{F}_{i'}\mid {\gamma}_A$
represent the same element of $\pi_1(C_3(\mathbb{C})\diagup \mathcal{S}_3,C_3(\mathbb{R})\diagup \mathcal{S}_3))$, namely $b_{tr}$. We represented $b_{tr}$
by the holomorphic map $\mathring{F}_{i'}$ from the rectangle $\tilde R$ into $C_3(\mathbb{C})\diagup \mathcal{S}_3$ that maps horizontal sides into $C_3(\mathbb{R})\diagup \mathcal{S}_3$. Hence,
\begin{equation}\label{eq100}
\Lambda(b_{tr})\leq \lambda(\tilde R) = \lambda(A\setminus L_A)\leq \lambda(A)\,.
\end{equation}

For $b\neq \sigma_j^k \, \Delta_3^{\ell}$ with $j$ equal to $1$ or $2$, and $k\neq 0$ and $\ell$ being integers, the statement of Lemma \ref{lem4} follows from Theorem I in the same way as Lemma \ref{lem2} follows from Theorem F. For $b=\sigma_j^k \, \Delta_3^{\ell}$ with $k=0$ the statement is trivial since then $ \vartheta({\rm Id})={\rm Id}$ and $\mathcal{L}_-({\rm Id})=0$.

To obtain the statement in the remaining case $b = \sigma_j^k \, \Delta_3^{2 \ell' +1} $ with
$j$ equal to $1$ or $2$, and $k$ and $\ell'$ being integers, we use Lemma \ref{lem3a}.
Notice that $\sigma_1\, \Delta_3 = \Delta_3 \, \sigma_2$ and $\sigma_2\, \Delta_3 = \Delta_3 \, \sigma_1$. Hence, $b^2 = \sigma_j^k\, \sigma_{j'}^k \,\Delta_3^{4 \ell'+2}$ with $\sigma_j\,\neq  \sigma_{j'}$. Let $\omega_2:A^2\to A$ be the two-fold unbranched covering of $A$ by an annulus $A^2$. The equality $\lambda(A^2)=2 \lambda(A)$ holds. Let ${q}_{A^2}$ be a point in $\omega_2^{-1}(q_A)$, and let ${L}_{q_{A^2}}$ be the lift of $L_A$ to $A^2$ that contains ${q}_{A^2}$.
Denote by ${\gamma}_{A^2}$ the loop $\omega_2^{-1}(\gamma_A)$ with base point $q_{A^2}$.
Then $ {F}\circ \omega_2 \mid {\gamma}_{A^2}$ represents $b^2$ and $(b^2)_{tr}$.
Lemma \ref{lem3a} applied to $ \sigma_j^k\, \sigma_{j'}^k \,\Delta_3^{4 \ell'+2}$ gives the estimate $2\log_+(3[\frac{|k|}{2}]) \leq \pi \lambda(A^2)=2\pi\lambda(A)$.
Since  $\vartheta(b)=\sigma_j^{2[\frac{|k|}{2}]\mbox{sgn}(k)}$,
the inequality \eqref{eq+} follows.
The lemma is proved. \hfill $\Box$

\medskip

\noindent {\bf Proof of Lemma \ref{lem3a}.}
By \cite{Jo2}, Lemma 1 and Proposition 6, statement 2,
\begin{align}\label{eq17}
\Lambda_{tr}( \sigma_j^k \, \sigma_{j'}^{k'} \,\Delta_3^{\ell})
\geq \frac{1}{\pi} (\log_+(3[\frac{|k|}{2}])+ \log_+(3[\frac{|k'|}{2}])).
\end{align}
Since by \eqref{eq100} the inequality $\Lambda_{tr}(\sigma_j^k \, \sigma_{j'}^{k'} \,\Delta_3^{\ell}) \leq \lambda(A)$ holds,
the lemma is proved. \hfill $\Box$
\medskip

We want to emphasize that periodic braids are not non-zero powers of a $\sigma_j$, so the lemma is true also for periodic braids. For each periodic braid $b$ of the form $\sigma_1 \sigma_2= \sigma_1 ^{-1} \, \Delta_3$, $(\sigma_1 \sigma_2)^2=\sigma_1 \, \Delta_3$, $\sigma_2 \sigma_1= \sigma_2 ^{-1} \, \Delta_3$, $(\sigma_2 \sigma_1)^2=\sigma_2 \, \Delta_3$,
and $\Delta_3$ the $\mathcal{L}_-(\vartheta(b))$ vanishes. However, for instance for the conjugate $\sigma_1^{-2k} \Delta_3 \sigma_1^{2k}= \sigma_1^{-2k} \sigma_2^{2k} \Delta_3$ of $\Delta_3$ we have $\mathcal{L}_-(\vartheta(\sigma_1^{-2k} \Delta_3 \sigma_1^{2k}))= 2 \log(3|k|)$. Another example, for the conjugate
$\sigma_2^{-2k}\, \sigma_1 \sigma_2\,  \sigma_2^{2k}$ of $\sigma_1 \sigma_2$ we have
\begin{align}\nonumber
\sigma_2^{-2k}\, \sigma_1 \sigma_2\,  \sigma_2^{2k}
%=  \sigma_2^{-2k-1}\, \sigma_2 \sigma_1 \sigma_2\,  \sigma_2^{2k} =
=  \sigma_2^{-2k-1}\, \Delta_3\,  \sigma_2^{2k} = \sigma_2^{-2k-1}\, \sigma_1^{2k} \,\Delta_3\,.
\end{align}
and  $\mathcal{L}_-(\vartheta(\sigma_2^{-2k}\, \sigma_1 \sigma_2\,  \sigma_2^{2k})) $ equals  $2 \log(3|k|)$.

Notice that the lemmas and Theorem I descend to statements on elements of $\mathcal{B}_3 \diagup \mathcal{Z}_3$ rather than on braids.
For an element $\textsf{b}$ of the quotient $\mathcal{B}_3 \diagup \mathcal{Z}_3$ we put $\vartheta(\textsf{b})= \vartheta(b)$ for any representative $b \in \mathcal{B}_3$ of $\textsf{b}$.

\smallskip

Lemma \ref{lem8a} below is an analog of Lemma \ref{lem3}. It follows from Lemma \ref{lem4} in the same way as Lemma \ref{lem3} follows from Lemma \ref{lem2}.

\begin{lemm}\label{lem8a}
Let $X$ be a connected finite open Riemann surface, and $F:X \to C_3(\mathbb{C}) \diagup \mathcal{S}_3$ be a holomorphic map that is transverse to the hypersurface $\mathcal{H}$ in $ C_3(\mathbb{C}) \diagup \mathcal{S}_3$.
Suppose $L_0$ is a simple relatively closed curve in $X$
such that $F(L_0)$ is contained in $\mathcal{H}$, and for a point $q \in L_0$ the point $F(q)$ is contained in the totally real space $C_3(\mathbb{R}) \diagup \mathcal{S}_3$.
Let $e^{(1)}$ and $e^{(2)}$ be primitive elements of $\pi_1(X,q)$. Suppose that for $e=e^{(1)}$, $e=e^{(2)}$, and $e=e^{(1)}e^{(2)}$ the free homotopy class $\widehat e$ intersects $L_0$.
Then either the two monodromies of $F$ modulo the center  $F_*(e^{(j)})\diagup \mathcal{Z}_3,\, j=1,2,\,$ are powers
of the same element $\sigma_j\diagup \mathcal{Z}_3$
of $\mathcal{B}_3\diagup \mathcal{Z}_3$, or each of them is the product of at most two elements $\textsf{b}_1$ and $\textsf{b}_2$ of $\mathcal{B}_3\diagup \mathcal{Z}_3$ with
\begin{equation}\label{eq101}
\mathcal{L}_-(\vartheta(\textsf{b}_j)) \leq  2\pi \lambda_{e^{(1)},e^{(2)}},\, j=1,2,
\end{equation}
where
\begin{equation}\nonumber
\lambda_{e^{(1)},e^{(2)}} \stackrel{def}=\max\{\lambda(A(\reallywidehat{e^{(1)}})),\, \lambda(A(\reallywidehat{e^{(2)}})),\, \lambda(A(\reallywidehat{e^{(1)}\,e^{(2)}}))\}.
\end{equation}
\end{lemm}

\noindent{\bf Proof.}
Suppose for an element $e\in \pi_1(X,q)$ the free homotopy class $\widehat e$ intersects $L_0$. By Lemma \ref{lem1} there exists an annulus $A$, a point $q_A\in A$, and a holomorphic map $\omega_A:(A,q_A)\to (X,q)$ that represents $e$. Moreover, the connected component of $(\omega_A)^{-1}(L_0)$ that contains $q_A$ has limit points on both boundary components of $A$. Put $F_A=F\circ\omega_A$. By the conditions of Lemma \ref{lem8a} $F_A(L_A)=F(L_0)\subset \mathcal{H}$ and $F_A(q_A)\in C_3(\mathbb{R})\diagup\mathcal{S}_3$.
Let $e_A$ be the generator of $\pi_1(A,q_A)$ for which $\omega_A(e_A)=e$.
The mapping $F_A:A\to C_3(\mathbb{C})\diagup \mathcal{S}_3$, the point $q_A$ and the curve $L_A$  satisfy the conditions of Lemma \ref{lem4}.
Notice that the equality $(F_A)_*(e_A)= F_*(e)$ holds.
Hence, if $F_*(e)$ is not a power of a $\sigma_j$ then inequality \eqref{eq+} holds for $F_*(e)$.

Suppose the two monodromies modulo center $F_*(e^{(j)})\diagup \mathcal{Z}_3,\, j=1,2,\,$
are not (trivial or non-trivial)
powers of the same element $\sigma_j \diagup \mathcal{Z}_3$ of $\mathcal{B}_3\diagup \mathcal{Z}_3$.
Then at most two of the elements, $F_*(e^{(1)})\diagup \mathcal{Z}_3$,  $F_*(e^{(2)})\diagup \mathcal{Z}_3$, and $F_*(e^{(1)}e^{(2)})\diagup \mathcal{Z}_3=F_*(e^{(1)})\diagup \mathcal{Z}_3\cdot F_*(e^{(2)})\diagup \mathcal{Z}_3 $, are powers of an element of the form $\sigma_j\diagup \mathcal{Z}_3$.

If the monodromies modulo center along two elements among $e^{(1)}$,  $e^{(2)}$, and $e^{(1)}e^{(2)}$ are
not (zero or non-zero) powers of a $\sigma_j\diagup \mathcal{Z}_3$ then by Lemma \ref{lem4} for each of these two monodromies modulo center inequality \eqref{eq101} holds, and the third monodromy modulo center is the product of two elements of $\mathcal{B}_3\diagup \mathcal{Z}_3$ for which  inequality \eqref{eq101} holds.
If the monodromies modulo center along two elements among $e^{(1)}$, $e^{(2)}$, and $e^{(1)}e^{(2)}$ have the form $\sigma_{j}^{k} \diagup \mathcal{Z}_3$ and  $\sigma_{j'}^{k'}\diagup\mathcal{Z}_3$,
then the $\sigma_j$ and the $\sigma_{j'}$ are different and $k$ and $k'$ are non-zero.
The third monodromy modulo center  has the form $\sigma_{j }^{\pm k} \sigma_{j'}^{\pm k'} \diagup \mathcal{Z}_3$ (or the order of the two factors interchanged).
Lemma \ref{lem3a} gives the inequality $\log_+(3\max([\frac{|k|}{2}]) + \log_+(3\max([\frac{|k'|}{2}])\leq \pi \lambda_{e^{(1)},e^{(2)}}$.
Since $\mathcal{L}_-(\vartheta(\sigma_{j}^{\pm k}))= \log(3[\frac{k}{2}])$ and  $\mathcal{L}_-(\vartheta(\sigma_{j'}^{\pm k'})) =\log(3[\frac{k'}{2}])$,
inequality \eqref{eq101} follows for the other two monodromies.
The lemma is proved.
\hfill $\Box$

\smallskip
The following lemma holds.
\begin{lemm}\label{lem7}
Let $X$ be a connected finite open Riemann surface, and $F:X\to C_3(\mathbb{C})\diagup \mathcal{S}_3$ a smooth mapping. Suppose for a base point $q_1$ of $X$ each element of $\pi_1(X,q_1)$ can be represented by a curve with base point $q_1$ whose image under $F$ avoids $\mathcal{H}$. Then all monodromies of $F$ are powers of the same periodic braid of period $3$.
\end{lemm}

\noindent {\bf Proof.}
Take the monodromy of $F$ along any curve with base point $q_1$. It has a power that is a pure $3$-braid $b$, and a representative of $b$ avoids $\mathcal{H}$.
Then for some integer $l$ the first and the last strand of $b \,\Delta_3^{2l}$ are fixed, and
a representative of $b \,\Delta_3^{2l}$ avoids $\mathcal{H}$. Hence,  $b\,\Delta_3^{2l}={\rm Id}$ and  $b=\Delta_3^{-2l}$. We saw that the monodromy of $F$ along each element $e\in \pi_1(X,q_1)$ is a periodic braid.

If a representative $f:[0,1]\to C_3(\mathbb{C})\diagup \mathcal{S}_3,\, f(0)=f(1),$ of a $3$-braid $b$ avoids $\mathcal{H}$, then the associated permutation $\tau_3(b)$ cannot be a transposition. Indeed, assume the contrary. Then there is a lift $\tilde f$ of $f$ to $C_3(\mathbb{C})$, for which $(\tilde{f}_1(1),\tilde{f}_2(1),\tilde{f}_3(1))=(\tilde{f}_3(0),\tilde{f}_2(0),\tilde{f}_1(0))$.
Let $L_t$ be the line in $\mathbb{C}$ that contains  $\tilde{f}_1(t)$ and $\tilde{f}_3(t)$, and is oriented so that running along $L_t$ in positive direction we meet first $\tilde{f}_1(t)$ and then $\tilde{f}_3(t)$.
The point $f_2(0)$ is not on $L_0$. Assume without loss of generality, that it is on the left of $L_0$ with the chosen orientation of $L_0$. Since for each $t\in [0,1]$ the three points $\tilde{f}_1(t),\,\tilde{f}_2(t)$ and $\tilde{f}_3(t)$ in $\mathbb{C}$ are not on a real line, the point $\tilde{f}_2(t)$ is on the left of $L_t$ with the chosen orientation. But the unoriented lines $L_0$ and $L_1$ coincide, and their orientation is opposite. This implies $\tilde{f}_2(1)\neq\tilde{f}_2(0)$, which is a contradiction.
We proved that all monodromies are periodic with period $3$.

There is a smooth homotopy $F_s,\, s\in[0,1],$ of $F$, such that $F_0=F$, each $F_s$ is different from $F$ only on a small neighbourhood of $q_1$, each $F_t$ avoids $\mathcal{H}$ on this neighbourhood of $q_1$,
and
$F_1(q_1)$ is the set of vertices of an equilateral triangle with barycenter $0$.
Since $F$ and $F_1$ are free homotopic, their monodromy homomorphisms are conjugate, and it is enough to prove the statement of the lemma for $F_1$.

For notational convenience we will keep the notation $F$ for the new mapping and assume that $F(q_1)$ is the set of vertices of an equilateral triangle with barycenter $0$.
The monodromy $F_*(e)$ along each element $e\in \pi_1(X,q_1)$ is a periodic braid of period $3$. Hence, $\tau_3(F_*(e))$ is a cyclic permutation. Consider the braid $g$ with base point $F(0)$ that corresponds to rotation by the angle $\frac{2\pi}{3}$, i.e. it is represented by the geometric braid $t\to e^{\frac{i2\pi t}{3}}F(0),\, t\in [0,1],$ that avoids $\mathcal{H}$.
There exists an integer $k$ such that $F_*(e)\, g^k$ is a pure braid that is represented
by a mapping that avoids $\mathcal{H}$.
Hence, $F_*(e)\, g^k$ represents $\Delta_3^{2l}$ for some integer $l$. We proved that for each $e\in \pi_1(X,q_1)$ the monodromy  $F_*(e)$ is represented by rotation of $F(0)$ around the origin by the angle $\frac{2\pi j}{3}$ for some integer $j$.
The Lemma is proved. \hfill $\Box$

\medskip

Let as before $X$ be a finite open connected Riemann surface. The following proposition is the main ingredient of the proof of Theorem \ref{thm3}. Let as before $\mathcal{E}\subset \pi_1(X,q_0)$ be the system of generators of the fundamental group with base point $q_0\in X$ that was chosen in Section 1.

\begin{prop}\label{prop4}
Let $(X\times \mathbb{P}^1,{\rm pr}_1,\mathbold{E},X)$
%$(\mathcal{X},\mathcal{P},X)$
be an
irreducible holomorphic special $(0,4)$-bundle
over a finite open Riemann surface $X$, that is not isotopic to a locally holomorphically trivial bundle.  Let $F(x),\, x \in X,$ be the set of finite distinguished points in the fiber over $x$. Assume that $F$ is transverse to $\mathcal{H}$.
Then there exists a complex affine mapping $M$ and a point $q\in X$ such that $M\circ
F(q)$ is contained in $C_3(\mathbb{R})\diagup \mathcal{S}_3$, and for an arc $\alpha$ in $X$ with initial point $q_0$ and terminating point $q$ and each element $e_j\in Is_{\alpha}(\mathcal{E})$ the monodromy modulo center $(M\circ F)_*(e_j)\diagup \mathcal{Z}_3$ can be written as product of at most $6$ elements  $\textsf{b}_{j,k},\, k=1,2,3,4,5,6,$ of $\mathcal{B}_3 \diagup \mathcal{Z}_3$ with
\begin{equation}\label{eq7}
\mathcal{L}_-(\vartheta(\textsf{b}_{j,k})) \leq 2\pi \lambda_{10}(X).
\end{equation}

If $X$ is a torus with a hole the monodromy along each $e_j$ is the product of at most $4$ elements with
$\mathcal{L}_-(\vartheta(\textsf{b}_{j,k})) \leq 2\pi \lambda_3(X)$, and in case of a planar domain the monodromy along each $e_j$ is the product of at most $6$ elements with
$\mathcal{L}_-(\vartheta(\textsf{b}_{j,k})) \leq 2\pi \lambda_8(X)$.

\end{prop}

\medskip

\noindent {\bf Proof of Proposition \ref{prop4}.} Since the bundle is not isotopic to a locally holomorphically trivial bundle, it is not possible that
all monodromies are powers of the same periodic braid, and by Lemma \ref{lem7} the set
\begin{align}\label{eq7a}
L\stackrel{def}= & \{ z \in X: F(z) \in\mathcal{H }\}
\end{align}
is not empty.

\noindent {\bf 1. A torus with a hole.} Let $X$ be a torus with a hole
and let  $\mathcal{E}=\{e'_0 ,\,e''_0\}$ be a set of  generators of $\pi_1(X,q_0)$ that is associated to a standard bouquet of circles for $X$.
There exists a connected component $L_0$ of $L$ which is not contractible and not contractible to the hole. Indeed, otherwise there would be a base point $q_1$ and a curve $\alpha_{q_1}$ that joins $q_0$ with $q_1$, such that for both elements of
$\mbox{Is}_{\alpha_{q_1}}(\mathcal{E})$
there would be representing loops with base point $q_1$ which do not meet $L$, and hence, by Lemma \ref{lem7} the monodromies along both elements would be powers of a single periodic braid of period $3$.

Hence, as in the proof of Proposition \ref{prop2} there exists a component $L_0$ of $L$, which is a simple smooth relatively closed curve in $X$, such that
%perhaps after switching $e'_0$ and $e''_0$ and orienting $L_0$ suitably,
the free homotopy class of one of the elements of $\mathcal{E}$, say of $e'_0$,
%$\widehat e'_0$
has positive intersection number with $L_0$. Put ${\sf e}'_0 =e'_0$.
%There is another element of the fundamental group
Moreover, the intersection number with $L_0$ is positive for the free homotopy class of one of the elements
%for one of the elements
${e_0''}^{\pm 1}$ or $e'_0 e''_0$.  %denoted by $e''_0$
Denote this element by ${\sf e}_0''$. (Since $\reallywidehat{e'_0 e''_0}=\reallywidehat{e''_0 e'_0}$ we may also put ${\sf e}_0''=e''_0 e'_0$ if the free homotopy class of $e'_0 e''_0$ intersects $L_0$.)
Put $\mathcal{E}_2'= \{{\sf e}_0', {\sf e}_0''\}$. The free homotopy class of each element of  $\mathcal{E}_2'$ and of the product of its two elements
intersects $L_0$.
One of the ${\sf e}_0'$ and ${\sf e}_0''$ is an element of $\mathcal{E}$, the other is in $\mathcal{E}\cup \mathcal{E}^{-1}$ or is the product of two
elements of $\mathcal{E}$. Each element of  $\mathcal{E}$ is the product of at most two elements of $\mathcal{E}'_2\cup {\mathcal{E}'_2}^{-1}$.

Move the base point $q_0$ to a point $q \in L_0$ along a curve $\alpha$, and consider the respective generators ${\sf e}'=\mbox{Is}_{\alpha}({\sf e}'_0)$ and ${\sf e}''=\mbox{Is}_{\alpha}({\sf e}''_0)$ of the fundamental group $\pi_1(X,q)$ with base point $q$. Since $F(L_0)\subset \mathcal{H}$ there is a complex affine mapping $M$ such that $M\circ F(q)\in C_3(\mathbb{R})\diagup \mathcal{S}_3$.
Since $F$ is irreducible, the monodromy maps modulo center $(M\circ F)_*({\sf e}')\diagup \mathcal{Z}_3$ and $(M\circ F)_*({\sf e}'')\diagup \mathcal{Z}_3$ are not powers of a single standard generator $\sigma_j\diagup \mathcal{Z}_3$  of $\mathcal{B}_3 \diagup \mathcal{Z}_3$ (see Lemma 7 of \cite{Jo5}).
Hence, the second option of Lemma \ref{lem8a} occurs. We obtain that each of the $(M\circ F)_*({\sf e}')\diagup \mathcal{Z}_3$ and $(M\circ F)_*({\sf e}'')\diagup \mathcal{Z}_3$
is a product of at most two elements ${\sf{b}}_j$ of $\mathcal{B}_3\diagup \mathcal{S}_3$ with $\mathcal{L}_-(\vartheta({\sf{b}}_j))\leq 2\pi\lambda_3(X)$. Hence, $(M\circ F)_*(e')\diagup \mathcal{Z}_3$ and $(M\circ F)_*(e'')\diagup \mathcal{Z}_3$ are products of at most $4$
elements of $\mathcal{B}_3\diagup \mathcal{Z}_3$ with this property.
The proposition is proved for tori with a hole.

\noindent {\bf  2. A planar domain.} Let $X$ be a planar domain. Maybe, after applying a M\"obius transformation, we represent $X$ as the Riemann sphere with holes $\mathcal{C}_j$, $j=1,\ldots,m+1,$ such that $\mathcal{C}_{m+1}$ contains $\infty$. Recall, that the set $\mathcal{E}$ of generators $e_{j,0}, \, j=1,\ldots,m,$ of the fundamental group $\pi_1(X,q_0)$
with base point $q_0$ is chosen so that $e_{j,0}$ is represented by a loop with base point $q_0$ that surrounds $\mathcal{C}_j$ counterclockwise and does not surround any other hole.
There is a connected component $L_0$ of $L$
of one of the following kinds.
Either $L_0$ has limit points on the boundary of two different holes (one of them may contain $\infty$) (first kind),
or a component $L_0$ has limit points on a single hole $\mathcal{C}_j,\, j\leq m+1,$
and $\mathcal{C}_j\cup L_0$ divides the plane $\mathbb{C}$ into two connected components each of which contains a hole (maybe, only the hole containing $\infty$) (second kind),
or there is a compact component $L_0$ that divides $\mathbb{C}$ into two connected components each of which contains at least two holes (one of them may contain $\infty$).
Indeed, suppose each non-compact component of $L$ has boundary points on the boundary of a single hole and the union of the component with the hole does not separate the remaining holes of $X$, and for each compact component of $L$ one of the connected components of its complement in $X$ contains at most one hole.
Then there exists a base point $q_1$, a curve $\alpha_{q_1}$ in $X$ with initial point $q_0$ and terminating point $q_1$, and a representative of each element of  $\mbox{Is}_{\alpha_{q_1}}(\mathcal{E})\subset\pi_1(X,q_1)$ that avoids $L$.
Lemma \ref{lem7} implies that all monodromies modulo center are powers of a single periodic element of $\mathcal{B}_3\diagup \mathcal{Z}_3$ which is a contradiction.

If there is a component $L_0$ of the first kind we may choose the same set of primitive elements $\mathcal{E}_2'\subset \mathcal{E}_2 \subset \pi_1(X,q_0)$ as in the proof of
Proposition \ref{prop2a}
in the planar case.
The free homotopy class of each element of $\mathcal{E}_2'$ and of the product of two such elements intersects $L_0$.
Moreover, each element of $\mathcal{E}$ is the product of at most two elements of $\mathcal{E}_2'$.
Let $\alpha_q$ be a curve in $X$ with initial point $q_0$ and terminating point $q$,
and $M$ a complex affine mapping, such that $(M\circ F)(q)\in C_3(\mathbb{R})\diagup \mathcal{S}_3$. Since $M\circ F$ is irreducible, the monodromies modulo center of $M\circ F$ along the elements of $\mbox{Is}_{\alpha}(\mathcal{E}_2')$ are not (trivial or non-trivial) powers of a single element $\sigma_j\diagup \mathcal{Z}_3$. Hence,
for each element of $\mbox{Is}_{\alpha}(\mathcal{E}_2')$ there
exists another element of $\mbox{Is}_{\alpha}(\mathcal{E}_2')$ so that
the second option of
Lemma \ref{lem8a} holds for this pair of elements of $\mbox{Is}_{\alpha}(\mathcal{E}_2')$. Therefore, the monodromy modulo center of $M\circ F$ along each element of $\mbox{Is}_{\alpha}(\mathcal{E}_2')$ is the product of at most two elements ${\sf b}_j\in\mathcal{B}_3\diagup \mathcal{Z}_3$ of $\mathcal{L}_-$ not exceeding $2\pi \lambda_4(X)$, and the monodromy modulo center of $M\circ F$ along each element
$\mbox{Is}_{\alpha}(\mathcal{E})$ is the product of at most $4$ elements of $\mathcal{B}_3\diagup \mathcal{Z}_3$ with $\mathcal{L}_-(\vartheta({\sf b}_j))$, each not exceeding $2\pi \lambda_4(X)$.

Suppose there is no component of the first kind but a component $L_0$ of the second kind. Assume first that all limit points of $L_0$ are on the boundary of a hole $\mathcal{C}_j$ that does not contain $\infty$. Put $\mathcal{E}_3'=\{e_{j,0}\} \cup_{1\leq k\leq m,\, k\neq j}\{ e_{j,0}^2 e_{k,0}\}$. Each element of $\mathcal{E}_3'$ is a primitive element and is the product of at most three generators contained in the set $\mathcal{E}$.
Further, each element of $\mathcal{E}$ is the product of at most three elements of $\mathcal{E}_3'\cup {\mathcal{E}_3'}^{-1}$.

The free homotopy class of each element of $\mathcal{E}_3'$ and of each product of two different elements of $\mathcal{E}_3'$
%$e_j$ and $e_j^2 e_k,\, 1\leq k\leq m,\, k\neq j$
intersects $L_0$. Indeed,
any curve that is contained in the complement of $\mathcal{C}_j\cup L_0$ has either winding number zero around $\mathcal{C}_j$ (as a curve in the complex plane $\mathbb{C}$), or its winding number around $\mathcal{C}_j$ coincides with the winding number around each of the holes in the bounded connected component of $\mathcal{C}_j$.
On the other hand the representatives of the free homotopy class of $e_{j,0}$
have winding number $1$ around $\mathcal{C}_j$ and winding number $0$ around each other hole that does not contain $\infty$. The representatives of the free homotopy class of $e_{j,0}^2 e_{k,0}$, $k\leq m,\, k\neq j$,
have winding number $2$ around $\mathcal{C}_j$,  winding number $1$ around $\mathcal{C}_k$, and winding number zero around each other hole $\mathcal{C}_l,\, l\leq m$.
The argument for products of two elements of $\mathcal{E}_3'$ is the same.

Choose a point $q\in L_0$, a curve $\alpha$ in $X$ with initial point $q_0$ and terminating point $q$, and a complex affine mapping $M$ such that $M\circ F(q)\in C_3(\mathbb{R})\diagup \mathcal{S}_3$.
Lemma \ref{lem8a} finishes the proof for this case in the same way as in the case when there is a component of first kind. In the present case each
$(M\circ F)_*(\tilde{e})\diagup \mathcal{Z}_3$, $\tilde{e}\in {\rm Is}_{\alpha}(\mathcal{E}_3')$, can be written as a product of at most $2$ factors ${\sf{b}}\in \mathcal{B}_3\diagup \mathcal{Z}_3$ with $\mathcal{L}_-(\vartheta({\sf{b}}))\leq 2\pi \lambda_6(X)$. Hence, each
$(M\circ F)_*(e_j)\diagup \mathcal{Z}_3$, $e_j=\mbox{Is}_{\alpha}(e_{j,0})$, can be written as a product of at most  $6$ factors ${\sf{b}}\in \mathcal{B}_3\diagup \mathcal{Z}_3$ with $\mathcal{L}_-(\vartheta({\sf{b}}))\leq 2\pi \lambda_6(X)$.

Assume that the limit points of $L_0$ are on the boundary of the hole $\mathcal{C}_{\infty}$ that contains $\infty$.
Let $\mathcal{C}_{j_0}$ and $\mathcal{C}_{k_0}$ be holes that are contained in different components of $X \setminus (L_0 \cup \mathcal{C}_{\infty})$, and let $e _{j_0,0}$ and $e _{k_0,0}$ be the elements of $\mathcal{E}$ whose representatives surround $\mathcal{C}_{j_0}$, and $\mathcal{C}_{k_0}$ respectively.
Denote by  $\mathcal{E}'_3$ the set that consists of the elements  $e _{j_0,0}e _{k_0,0}\,$,  $\;e _{j_0,0}^2e _{k_0,0}\,$, and all elements $ e _{j_0,0}e _{k_0,0} \tilde {e}_0$ with $\tilde{e}_0$ running over $ \mathcal{E}\setminus\{ e _{j_0,0},e _{k_0,0}\}$.
Each element of $\mathcal{E}_3'$ is the product of at most $3$ elements of $\mathcal{E}$, and each element of $\mathcal{E}$ is the product of at most $3$ elements of $\mathcal{E}'_3\cup (\mathcal{E}'_3)^{-1}$.

Each element of $\mathcal{E}'_3$ and each product of at most two  different elements of $\mathcal{E}'_3$ intersects $L_0$. Indeed, if a closed curve is contained in one of the components of $X\setminus (L_0\cup \mathcal{C}_{\infty})$ then its winding number around each hole contained in the other component is zero. But for all mentioned elements there is a hole in each component of  $X\setminus (L_0 \cup \mathcal{C}_{\infty})$ such that the winding number of the free homotopy class  of the element around the hole does not vanish.
Lemma \ref{lem8a} applies with the same meaning of $q$, $\alpha$, and $M$ as before.
Again, each
$(M\circ F)_*(e_j)\diagup \mathcal{Z}_3$, $e_j=\mbox{Is}_{\alpha}(e_{j,0})$, can be written as a product of at most  $6$ factors ${\sf{b}}\in \mathcal{B}_3\diagup \mathcal{Z}_3$ with $\mathcal{L}_-(\vartheta({\sf{b}}))\leq 2\pi \lambda_6(X)$.

Notice that in case of $m+1=3$ holes only these two possibilities for the curve $L_0$ may occur. Hence in this case we find a set $\mathcal{E}_3'=\{{\sf e}'_0,{\sf e}''_0\} \subset \pi_1(X,q_0)$, such that one of the elements of $\mathcal{E}_3'$ is the product of at most two elements of $\mathcal{E}\cup \mathcal{E}^{-1}$, and each of the monodromies $F_*({\sf e}'_0)$ and $F_*({\sf e}''_0)$
%is in $\mathcal{E}\cup \mathcal{E}^{-1}$ or
is the product of at most two elements ${\sf{b}}\in \mathcal{B}_3\diagup \mathcal{Z}_3$ with $\mathcal{L}_-(\vartheta({\sf{b}}))\leq 2\pi \lambda_5(X)$.
Moreover, $e$ and $e'$ are products of at most three factors, each an element of  $\mathcal{E}_3'\cup \mathcal{E}_3'^{-1}$.
%and not all factors being equal.

Suppose there are no components of $L$ of the first or the second kind, but there is a connected component $L_0$ of $L$ of the third kind.
Let $\mathcal{C}_{j_0}$ be a hole contained in the bounded component of the complement of $L_0$, and let $\mathcal{C}_{k_0},\, k_0\leq m,$ be a hole that is contained in the unbounded component of $X \setminus L_0$. Let $e _{j_0,0}$ and $e _{k_0,0}$ be the elements of $\mathcal{E}$ whose representatives surround $\mathcal{C}_{j_0}$, and $\mathcal{C}_{k_0}$ respectively.
Consider the set $\mathcal{E}'_4$ consisting of the following elements: $e_{j_0,0} e_{k_0,0}$,
$e_{j_0,0}^2 e_{k_0,0}$, and $e_{j_0,0}^2 e_{k_0,0} \tilde{ e}_0$ for each $\tilde{e}_0\in \mathcal{E}$ different from $e_{j_0,0}$ and $e_{k_0,0}$. Each element of $\mathcal{E}'_4$ is the product of at most $4$ elements of $\mathcal{E}$ and each element of $\mathcal{E}$ is the product of at most $3$ elements of $\mathcal{E}_4'\cup (\mathcal{E}_4')^{-1}$. The product of two different elements of $\mathcal{E}'_4$ is contained in $\mathcal{E}'_8$.

The free homotopy classes of each element of $\mathcal{E}'_4$ and of each product of two different elements of $\mathcal{E}'_4$ intersects $L_0$.
Indeed, if a loop is contained in the bounded connected component of $X\setminus L_0$, its winding number around the holes $\mathcal{C}_j\,, j\leq m,$ contained in the unbounded component is zero. If a loop is contained in the unbounded connected
component of $X\setminus L_0$, its winding numbers around all holes contained in the bounded connected component are equal. But the winding number of $e_{j_0,0} e_{k_0,0}$ and $e_{j_0,0}^2 e_{k_0,0}$ around the hole $\mathcal{C}_{j_0}$ is positive and the winding number around the other holes that are contained in the bounded connected component of $X\setminus L_0$ vanishes, hence the representatives of these two elements cannot be contained in the unbounded  component of $X\setminus L_0$. Since the winding number of representatives of these elements around $\mathcal{C}_ {k_0}$ is positive, the representatives cannot be contained in the bounded component of $X\setminus L_0$.
For representatives of  the elements $e_{j_0,0}^2 e_{k_0,0} \tilde{e}_0$
the winding number around $\mathcal{C}_{j_0}$ equals $2$, the winding number around any other hole in the bounded component of $X\setminus L_0$ is at most $1$, and the winding number around
$\mathcal{C}_{k_0}$ equals $1$. Hence, the free homotopy classes of the mentioned elements must intersect both components of $X \setminus L_0$, hence they intersect $L_0$.

Representatives of any product of two elements of  $\mathcal{E}_4'$  have winding number around $\mathcal{C}_{j_0}$ at least $3$, the winding number around any other hole in the bounded component of $X\setminus L_0$ is at most $1$, and the winding number around
$\mathcal{C}_{k_0}$ equals $2$.  Hence, the free homotopy classes of these elements intersect $L_0$.

For a point $q\in L_0$, a curve $\alpha$ in $X$ joining $q_0$ and $q$,
and a complex affine mapping $M$ for which $M\circ F(q)\in C_3(  \mathbb{R})\diagup \mathcal{S}_3$, an application of Lemma \ref{lem8a} proves that in this case each
$(M\circ F)_*(e_j)\diagup \mathcal{Z}_3$, $e_j=\mbox{Is}_{\alpha}(e_{j,0})$, can be written as a product of at most  $6$ factors ${\sf{b}}\in \mathcal{B}_3\diagup \mathcal{Z}_3$ with $\mathcal{L}_-(\vartheta({\sf{b}}))\leq 2\pi \lambda_8(X)$.
Proposition \ref{prop4} is proved in the planar case with a slightly better constant.

\noindent {\bf 3. The general case.}
Since not all monodromies are powers of a single element of $\mathcal{B}_3\diagup \mathcal{Z}_3$ that is either periodic or reducible, there exists a pair of generators $e_0'$, $e_0''$ in $\mathcal{E}$, such that the monodromies along them are not powers of a single periodic or reducible element. Consider the projection $\omega^{\langle  e_0', e_0''\rangle}: \tilde{X}\to X(\langle  e_0', e_0''\rangle)$.
By the proof for tori with a hole or for $\mathbb{P}^1$ with three holes there exist a relatively closed curve $L_{\langle  e_0', e_0''\rangle}$ in $X(\langle  e_0', e_0''\rangle)$ and a M\"obius transformation $M$,
%(chosen out of three M\"obius transformation $M_l,\, l=0,1,2,$
such that for $F=M\circ f$ the mapping $F_{\langle  e_0', e_0''\rangle}  =F \circ \omega_{\langle  e_0', e_0''\rangle}$ takes $L_{\langle  e_0', e_0''\rangle}$ into $\mathcal{H}$, and takes a chosen point $q_{\langle  e_0', e_0''\rangle}\in L_{\langle  e_0', e_0''\rangle}$ to a point in $C_3(\mathbb{R})\diagup \mathcal{S}_3$.

Choose a point $\tilde{q}\in \tilde X$, for which $\omega^{\langle e'_0, e''_0\rangle}(\tilde{q})=q_{\langle e_0, e'_0\rangle}$. Let $\tilde{\alpha}$ be a curve in $\tilde X$ with initial point $\tilde{q}_0$ and terminating point $\tilde{q}$.
Then $\alpha_{\langle e_0', e''_0\rangle}\stackrel{def}=\omega^{\langle e'_0, e''_0\rangle}(\tilde{\alpha})$ is a curve in $X(\langle e'_0, e''_0\rangle)$ with initial point ${q_0}_{\langle e'_0, e''_0\rangle}$ and terminating point $ q_{\langle e'_0, e''_0\rangle}$, and the curve $\alpha_{\langle e'_0, e''_0\rangle}$ in $X(\langle e'_0, e''_0\rangle)$ and the point $\tilde{q}_0$ in the universal covering $\tilde X$ of $X(\langle e'_0, e''_0\rangle)$ are compatible. Put $\alpha=\omega_{\langle e'_0, e''_0\rangle}(\alpha_{\langle e'_0, e''_0\rangle})$, and for each $e_0\in \pi_1(X,q_0)$ we denote as before the element ${\rm Is}_{\alpha}({ e}_0)$ by $e$.

As in the case of a torus with a hole or $\mathbb{P}^1$ with three holes there are elements ${\sf e}_0'$ and ${\sf e}_0''$, one of them contained in $\mathcal{E}$
%$\cup \mathcal{E}^{-1}$
or equal to the product of at most two factors among the $e'_0$ and $e''_0$, the second either contained in $\mathcal{E}\cup \mathcal{E}^{-1}$, or equal to the product of at most three factors among the $e'_0$ and $e''_0$,
%the other equal to the product of at most $3$ factors among the $e'$ and $e''$,
such that the free homotopy classes of  $({\sf e}_0')_{\langle e'_0, e''_0\rangle}$, of $({\sf e}''_0)_{\langle e'_0, e''_0\rangle}$, and of their product intersect  $ L_{\langle e'_0, e''_0\rangle}$. (For the definitions of  $({\sf e}_0')_{\langle e'_0, e''_0\rangle}$ and of $({\sf e}''_0)_{\langle e'_0, e''_0\rangle}$ see paragraph 3.1.) Moreover, $e'_0$ and $e''_0$ are products of at most three factors, each being either $({\sf e}'_0)^{\pm 1}$ or $({\sf e}''_0)^{\pm 1}$.
Put ${\sf e}'_{\langle e'_0, e''_0\rangle}={\rm Is}_{\alpha_{\langle e'_0, e''_0\rangle} }   (({\sf e}_0')_{\langle e'_0, e''_0\rangle})$, ${\sf e}''_{\langle e'_0, e''_0\rangle}={\rm Is}_{\alpha_{\langle e'_0, e''_0\rangle} }  ( ({\sf e}_0'')_{\langle e'_0, e''_0\rangle})$.
By Lemma \ref{lem8a}  each monodromy $(F_{\langle e'_0, e''_0\rangle})_*({\sf e}'_{\langle e'_0, e''_0\rangle})= F_*({\sf e}')$ and $(F_{\langle e'_0, e''_0\rangle})_*({\sf e}''_{\langle e'_0, e''_0\rangle})= F_*({\sf e}'')$ is the product of at most two elements ${\sf b}_j\in\mathcal{B}_3\diagup \mathcal{Z}_3$
with $\mathcal{L}_-(\vartheta({\sf b}_j))\leq 2\pi \lambda_5(X)$. Since $e'$ and $e''$ are products of at most three elements among $({\sf e}')^{\pm 1}$ and  $({\sf e}'')^{\pm 1}$,
%of $\mathcal{E}\cup \mathcal{E}^{-1}$,
each of the monodromies $F_*(e')$ and $F_*(e'')$ is the product of at most $6$  elements ${\sf b}_j\in\mathcal{B}_3\diagup \mathcal{Z}_3$
with $\mathcal{L}_-(\vartheta({\sf b}_j))\leq 2\pi \lambda_5(X)$.

Take any element $e_0\in \mathcal{E}\setminus\{e'_0,e''_0\}$.  Either the pair of monodromies ($F_*({\sf e}')$, $F_*(e)$) or  the pair of monodromies ($F_*({\sf e}'')$, $F_*(e)$) does not consist of two powers of the same element of $\mathcal{B}_3\diagup \mathcal{S}_3$ that is either periodic or reducible. Suppose this is so for  the pair ($F_*({\sf e}')$, $F_*(e)$).
%${\sf e}''$ the monodromies of $F$ along

Let $L_{\langle {\sf e}'\rangle}$ be the connected component of $(\omega^{\langle {\sf e}', {\sf e}''\rangle}_{\langle {\sf e}'\rangle})^{-1}(L_{\langle {\sf e}',{\sf e}''\rangle})$ that contains $\omega^{\langle {\sf e}'\rangle}(\tilde{q})$. By Lemma \ref{lem1}, applied to the holomorphic projection $\tilde{X}\diagup ({\rm Is}^{\tilde{q}})^{-1}(\langle {\sf e}'\rangle) \to X(\langle {\sf e}', {\sf e}''    \rangle)$, the free homotopy class $\reallywidehat{{\sf e}'_{\langle {\sf e}'\rangle}}$ intersects $L_{\langle {\sf e}'\rangle}$. (For the definition of ${\sf e}'_{\langle {\sf e}'\rangle}$ see paragraph 3.1.)
As in the proof of Proposition \ref{prop2} we consider the Riemann surface $X(\langle e, {\sf e}'\rangle)$ and the curve $L_{\langle e, {\sf e}'\rangle}= \omega^{\langle e, {\sf e}'\rangle}_ {\langle  {\sf e}'\rangle}(L_{\langle {\sf e}'\rangle}) $ (see paragraph 3.3. of the proof of proposition \ref{prop2}). As there we see that
the free homotopy class
$\reallywidehat{ {\sf e}'_{\langle e, {\sf e}'\rangle}}$ intersects $L_{\langle e, {\sf e}'\rangle}$.
The system $e_{\langle e, {\sf e}'\rangle}, {\sf e}'_{\langle e, {\sf e}'\rangle}$ is associated to a standard bouquet of circles for $X(\langle e, {\sf e}'\rangle)$ (though the representing curves of ${\sf  e}'$ in $X$ may not be simple closed curves or may intersect representing curves of $e$). This can be seen in the same way
as in the proof of Proposition \ref{prop2}.
Apply the arguments, used for $X(\langle e',e''\rangle )$ and the generators $e'_{\langle e', e''\rangle},e''_{\langle e', e''\rangle}$ of the fundamental group $\pi_1(X(\langle e', e''\rangle), q_{\langle e', e''\rangle})$,  to
$X(\langle e, {\sf e}'\rangle)$ and the generators  $e_{\langle e, {\sf e}'\rangle}, {\sf e}'_{\langle e, {\sf e}'\rangle}$ of the fundamental group $\pi_1(X(\langle e, {\sf e}'\rangle),   q_{\langle e, {\sf e}'\rangle})$.

In the case when $X(\langle e, {\sf e}'\rangle)$ is a torus with a hole, the intersection number of $\reallywidehat{{\sf e}'_{\langle e, {\sf e}'\rangle} }$ with $L_{\langle e, {\sf e}'\rangle}$ is non-zero. Put $\mathfrak{e}'={\sf e}'$. For one of the choices $e^{\pm 1}$, or ${\sf e}'\, e$, denoted by $\mathfrak{e}''$, the free homotopy classes of $\mathfrak{e}'_{\langle e, {\sf e}'\rangle}$, $\mathfrak{e}''_{\langle e, {\sf e}'\rangle}$, and of their product intersect $L_{\langle e, {\sf e}'\rangle}$. Moreover, $e$ is the product of at most two factors, each being $(\mathfrak{e}')^{\pm 1}$, or $(\mathfrak{e}'')^{\pm 1}$.

In case $X(\langle e, {\sf e}'\rangle)$ is planar, the curve $L_{\langle e, {\sf e}'\rangle}$ must have limit points on the hole that corresponds to the generator ${\sf e}'_{\langle e, {\sf e}'\rangle}$ of the fundamental group $\pi_1(X(\langle e, {\sf e}'\rangle), q_{\langle e, {\sf e}'\rangle})$. We find elements $\mathfrak{e}'$ and $\mathfrak{e}''$ such that $\mathfrak{e}'= {\sf e}' $ and $\mathfrak{e}''$ is either equal to $e^{-1}$,
%${\sf e}'$
or to the product of at most three factors, one being equal to $e$ and the others equal to ${\sf e}' $, and the free homotopy classes of $\mathfrak{e}'_{\langle e, {\sf e}'\rangle}$, $\mathfrak{e}''_{\langle e, {\sf e}'\rangle}$, and their product intersect  $L_{\langle e, {\sf e}'\rangle}$. Moreover, %$\mathfrak{e}=e$
$e$ is the product of at most $3$ factors, each being equal to $(\mathfrak{e}'')^{\pm 1}$ or
$(\mathfrak{e}')^{\pm 1}$.
%, and ${\sf e}' $ is the product of at most three elements
%of $\mathcal{E}\cup \mathcal{E}^{-1}$.

In both cases for $X(\langle e, {\sf e}'\rangle)$ the element $\mathfrak{e}' \mathfrak{e}''$ is the product of at most $10$ elements of $\mathcal{E} \cup \mathcal{E}^{-1}$.
Lemma \ref{lem8a} implies, that $F_*(\mathfrak{e})$ and $F_*(\mathfrak{e}')$ are products of at most two factors $\sf b$ with $\mathcal{L}_-(\vartheta({\sf b}))$ not exceeding $2\pi \lambda_{10}(X)$. Hence, $F_*(e)$ is the product of at most $6$ factors $\sf b$ with $\mathcal{L}_-(\vartheta({\sf b}))$ not exceeding $2\pi \lambda_{10}(X)$.
We obtain the statement of Proposition \ref{prop4} in the general case.
Proposition \ref{prop4} is proved. \hfill $\Box$

\medskip

\noindent {\bf Proof of Theorem \ref{thm3}.}
Let $X$ be a connected Riemann surface of genus $\sf g$ with ${\sf m}+1\geq 1$ holes. Since each holomorphic $(0,3)$-bundle with a holomorphic section on $X$ is isotopic to a holomorphic special $(0,4)$-bundle, we need to estimate the number of isotopy classes of irreducible smooth special $(0,4)$-bundles on $X$, that contain a holomorphic bundle.
By Lemma 4 of \cite{Jo5}
the monodromies of such a bundle are not powers of a single element of $\mathcal{B}_3 \diagup \mathcal{Z}_3$ which is conjugate to a $\sigma_j\diagup \mathcal{Z}_3$, but they may be powers of a single periodic element of $\mathcal{B}_3 \diagup \mathcal{Z}_3$ (equivalently, the isotopy class may contain a locally holomorphically trivial holomorphic bundle).

Consider an irreducible special holomorphic $(0,4)$-bundle on $X$
which is not isotopic
to a locally holomorphically trivial bundle.
Let $F(x),\, x \in X,$ be the set of finite distinguished points in the fiber over $x$.
By the Holomorphic Transversality Theorem \cite{KZ}
the mapping $F:X\to C_3(\mathbb{C})\diagup \mathcal{S}_3$ can be approximated on relatively compact subsets of $X$ by holomorphic mappings
%to the symmetrized configuration space
that are transverse to $\mathcal{H}$.  Similarly as in the proof of Theorem \ref{thm1}
we will therefore assume in the following (after slightly shrinking $X$ to a deformation retract of $X$ and approximating $F$) that $F$ is transverse to $\mathcal{H}$.

By Proposition \ref{prop4}
there exists a complex affine mapping $M$ and a point $q\in X$ such that $M\circ
F(q)$ is contained in $C_3(\mathbb{R})\diagup \mathcal{S}_3$, and for an arc $\alpha$ in $X$ with initial point $q_0$ and terminating point $q$ and each element $e_j\in {\rm Is}_{\alpha}(\mathcal{E})$ the monodromy $(M\circ F)_*(e_j)\diagup \mathcal{Z}_3$ of the bundle can be written as product of at most $6$ elements  $\textsf{b}_{j,k},\, k=1,2,3,4,5,6,$ of $\mathcal{B}_3 \diagup \mathcal{Z}_3$ with
\begin{equation}\label{eq7}
\mathcal{L}_-(\vartheta(\textsf{b}_{j,k})) \leq 2\pi \lambda_{10}(X).
\end{equation}
\noindent The mappings $F$ and $M\circ F$ from $X$ into the symmetrized configuration space are free homotopic.

Consider an isotopy class of special $(0,4)$-bundles that corresponds to a conjugacy class of homomorphisms $\pi_1(X,q_0)\to \mathcal{B}_3\diagup \mathcal{Z}_3$ whose image is generated by a single periodic element of $\mathcal{B}_3\diagup \mathcal{Z}_3$.
Up to conjugacy we may assume that this element is one of the following:
${\rm Id},\, \Delta_3\diagup \mathcal{Z}_3,\, (\sigma_1 \sigma_2 )\diagup \mathcal{Z}_3 , \,(\sigma_1 \sigma_2)^{-1} \diagup \mathcal{Z}_3$. For each of these elements $\textsf{b}$ the equality $\mathcal{L}_-(\vartheta(\textsf{b}))=0$ holds. Hence, in this case the isotopy class contains a smooth mapping $\tilde F$ such that for each $e_{j,0}\in \mathcal{E}$ the monodromy $(M\circ F)_*(e_{j,0})\diagup \mathcal{Z}_3$ of the bundle can be written as product of at most $6$ elements  $\textsf{b}_{j,k},\, k=1,2,3,4,5,6,$ of $\mathcal{B}_3 \diagup \mathcal{Z}_3$ satisfying inequality \eqref{eq7}.

The same argument as in the proof of Theorem \ref{thm1} shows the following fact. Each irreducible free homotopy class of mappings $X\to C_3(\mathbb{C})\diagup \mathcal{S}_3$ that contains a holomorphic mapping contains a smooth mapping $\tilde F$ such that for each $e_{j,0}\in \mathcal{E}$ the monodromy $\tilde{F}_*(e_{j,0})\diagup \mathcal{Z}_3$ of the bundle can be written as product of at most $6$ elements  $\textsf{b}_{j,k},\, k=1,2,3,4,5,6,$ of $\mathcal{B}_3 \diagup \mathcal{Z}_3$ satisfying inequality \eqref{eq7}.

Using Lemma 1 of \cite{Jo3} the number of elements of $  \textsf{b}\in \mathcal{B}_3 \diagup \mathcal{Z}_3$ (including the identity),
for which $\mathcal{L}_-(\vartheta(\textsf{b}))\leq 2\pi \lambda_{10}(X)$,
is estimated as follows.
The element ${\sf w}\stackrel{def}=\vartheta({\sf b}) \in \mathcal{PB}_3 \diagup \mathcal{Z}_3$ can be considered as a reduced word in the free group generated by $a_1=\sigma_1^2\diagup \mathcal{Z}_3$ and $a_2=\sigma_2^2\diagup \mathcal{Z}_3$.
By Lemma 1 of \cite{Jo3}
there are no more than $ \frac{1}{2}\exp(6 \pi \lambda_{10}(X))+1\leq \frac{3}{2}\exp(6 \pi \lambda_{10}(X))$ reduced words $\sf w$ in $a_1$ and $a_2$
(including the identity) satisfying the inequality $\mathcal{L}_-(\textsf{w})\leq 2\pi \lambda_{10}(X)$.

For a given element ${\sf w}\in \mathcal{PB}_3 \diagup \mathcal{Z}_3$ (including the identity) we describe now all elements ${\textsf{b}}$
of $\mathcal{B}_3 \diagup \mathcal{Z}_3$ with
$\vartheta(\textsf{b})=\textsf{w}$.
If ${\textsf{w}}\neq \mbox{Id}$ these are
the following elements.
If the first term of ${\sf w}$ equals
$a_j^k$ with $k\neq 0$, then the possibilities are $\textsf{b}={\sf w} \cdot (\Delta_3^{\ell}\diagup \mathcal{Z}_3)$ with $\ell=0$ or $1$,  $\textsf{b}=(\sigma_j^{\begin{tiny}{\mbox{sgn}}\end{tiny} k}\diagup\mathcal{Z}_3)\cdot  {\sf w}  \cdot (\Delta_3^{\ell}\diagup \mathcal{Z}_3)$ with $\ell=0$ or $1$, or  $\textsf{b}=(\sigma_{j'}^{\pm 1}\diagup \mathcal{Z}_3) \cdot {\sf w} \cdot (\Delta_3^{\ell}\diagup \mathcal{Z}_3)$ with $\ell=0$ or $1$ and $\sigma_{j'}\neq
\sigma_{j}$. Hence, for $\textsf{w}\neq \mbox{Id}$ there are $8$ possible choices of elements ${\textsf{b}}\in \mathcal{B}_3 \diagup \mathcal{Z}_3$ with $\vartheta({\textsf{b}})= {\sf w} $.

If $\textsf{b}= \mbox{Id}$ then the choices are $\Delta^{\ell}\diagup \mathcal{Z}_3$ and $(\sigma_j^{\pm 1}\Delta^{\ell})\diagup \mathcal{Z}_3$ for $j=1,2,$ and $\ell=0$ or $\ell=1$. These are $10$ choices.  Hence, there are no more than  $15 \exp( 6 \pi \lambda_{10}(X))$ different
elements $\textsf{b}\in \mathcal{B}_3 \diagup \mathcal{Z}_3$  with $\mathcal{L}_-(\vartheta(\textsf{b}))\leq 2 \pi \lambda_{10}(X)$.

Each monodromy is the product of at most six elements ${\sf{b}}_j$ of
$\mathcal{B}_3 \diagup \mathcal{Z}_3$ with $\mathcal{L}_-(\vartheta(\textsf{b}_j))\leq 2 \pi \lambda_{10}(X)$.
Hence, for each monodromy there are no more than $(15 \exp( 6 \pi \lambda_{10}(X)))^{6}$
possible choices. We proved that there are up to isotopy no more than $(15 \exp( 6 \pi \lambda_{10}(X)))^{6(2g+m)}$
%$= (3^6\cdot 5^6\cdot \exp(36 \pi \lambda_8(X)))^{2g+m}$
irreducible holomorphic $(0,3)$-bundles with a holomorphic section over $X$.
Theorem \ref{thm3} is proved. \hfill $\Box$
\smallskip

Notice that we proved a slightly stronger statement, namely, over a Riemann surface of genus $g$ with $m+1\geq 1$ holes there are no more than
$(15 \exp( 6 \pi \lambda_{10}(X)))^{6(2g+m)}$ isotopy classes of smooth $(0,3)$-bundles with a
smooth section that contain a holomorphic bundle with a holomorphic section that is either irreducible or isotopic to the trivial bundle.

\medskip

\noindent {\bf Proof of Theorem \ref{thm2}.} Proposition \ref{prop1} and Theorem \ref{thm3} imply Theorem \ref{thm2} as follows. Suppose an isotopy class of smooth $(1,1)$-bundles over a finite open Riemann surface $X$ contains a holomorphic bundle. By Proposition \ref{prop1} the class contains a holomorphic bundle which is the double branched covering of a holomorphic special $(0,4)$-bundle.
If the $(1,1)$-bundle is irreducible then also the $(0,4)$-bundle is irreducible.
There are up to isotopy no more than $\big(15(\exp(6 \pi \lambda_{10}(X)))\big)^{6(2g+m)}$
holomorphic special $(0,4)$-bundles over $X$ that are either irreducible or isotopic to the trivial bundle.

By Theorem G and Theorem \ref{thm3} there are no more than $\big(15\big(\exp(6 \pi \lambda_{10}(X)))\big)^{6(2g+m)}$
conjugacy classes of monodromy homomorphisms that correspond to
a special holomorphic $(0,4)$-bundle over $X$
that is either irreducible or isotopic to the trivial bundle.
Each monodromy homomorphism of the holomorphic double branched covering is a lift of the respective monodromy homomorphism of the holomorphic special $(0,4)$-bundle.
Different lifts of a monodromy mapping class of a special $(0,4)$-bundle differ by involution, and the fundamental group of $X$ has $2g+m$ generators.
Using Theorem G for $(1,1)$-bundles, we see that
there are no more than  $2^{2g+m}\big(15(\exp(6 \pi \lambda_{10}(X)))\big)^{6(2g+m)}=\big(2 \cdot 15^6\cdot\exp(36 \pi \lambda_{10}(X))\big)^{2g+m}$
isotopy classes of $(1,1)$-bundles that contain a holomorphic bundle that is either irreducible or isotopic to the trivial bundle.
Theorem \ref{thm2} is proved.
\hfill $\Box$

\medskip

\noindent For convenience of the reader we give the short proofs of the Corollaries \ref{cor1a} and \ref{cor1b}.
Such statements
are known in principle, but the case  considered here is especially simple.

\medskip

\noindent {\bf Proof of Corollary \ref{cor1a}.}
We will prove that on a punctured Riemann surface there are no non-constant reducible holomorphic mappings to the twice punctured complex plane and that any homotopy class of mappings from a punctured Riemann surface to the twice  punctured complex plane contains at most one holomorphic mapping. This implies the corollary.

Recall that a holomorphic mapping $f$ from any punctured Riemann surface $X$ to the twice punctured complex plane extends by Picard's Theorem to a meromorphic function $f^c$ on the closed  Riemann surface $X^c$.
Suppose now that $X$ is a punctured Riemann surface and that the mapping $f:X \to \mathbb{C}\setminus \{-1,1\}$ is reducible, i.e. it is homotopic to a mapping into a punctured disc contained in $\mathbb{C}\setminus \{-1,1\}$. Perhaps after composing $f$ with a M\"obius transformation we may suppose that this puncture equals $-1$.  Then the meromorphic extension $f^c$ omits the value $1$. Indeed, if $f^c$ was equal to $1$ at some puncture of $X$, then $f$ would map
the boundary of a sufficiently small disc on $X^c$
that contains the puncture to a loop in $\mathbb{C}\setminus \{-1,1\}$ with non-zero winding number around $1$ ,
which contradicts the fact that $f$ is homotopic to a mapping into a disc punctured at $-1$ and contained in $\mathbb{C}\setminus \{-1,1\}$. Hence, $f^c$ is a meromorphic function on a compact Riemann surface that omits a value, and, hence $f$ is constant. Hence, on a punctured Riemann surface there are no non-constant reducible holomorphic mappings to $\mathbb{C}\setminus \{-1,1\}$.

Suppose $f_1$ and $f_2$ are non-constant homotopic holomorphic mappings from the punctured Riemann surface $X$ to the twice punctured complex plane. Then for their meromorphic extensions $f^c_1$ and $f^c_2$ the functions  $f^c_1-1$ and $f^c_2-1$ have the same divisor on the closed Riemann surface $X^c$. Indeed, suppose, for instance, that $f^c_1-1$ has a zero of order $k >0$
at a puncture $p$. Then for the boundary $\gamma$ of a small disc in $X^c$
around $p$
the curve $(f_1-1) \circ \gamma$ in $\mathbb{C}\setminus \{-2,0\}$ has index $k$
with respect to the origin.
Since $f_2-1$ is homotopic to $f_1-1$ as mapping to $\mathbb{C}\setminus \{-2,0\}$, the curve  $(f_2-1) \circ \gamma$ is free homotopic to $(f_1-1) \circ \gamma$ . Hence, $f_2-1$ has a zero of order $k$
at $p$. Applying the same arguments with $0$ replaced by $\infty$, we obtain that $f^c_1-1$ and $f^c_2-1$ have the same divisor. Hence, $f^c_1-1$ and $f^c_2-1$ differ by a non-zero multiplicative constant. %$\alpha$.
Since the functions are non-constant they must take the value $-2$. By the same reasoning as above the functions are equal to $-2$ simultaneously. Hence, the multiplicative constant is equal to $1$.
We proved that non-constant homotopic holomorphic maps from punctured Riemann surfaces
to $\mathbb{C}\setminus \{-1,1\}$ are equal. \hfill $\Box$

\medskip
\noindent {\bf Proof of Corollary \ref{cor1b}.}
We need the following fact. For each special $(0,4)$-bundle $\mathfrak{F}=(X\times \mathbb{P}^1, {\rm pr}_1, {\mathbold E}, X)$ there is a finite unramified covering $\hat{{\sf P}}:\hat X\to X$ of $X$, such that $\mathfrak{F}$ lifts to a  special $(0,4)$-bundle
$(\hat{X}\times\mathbb{P}^1,{\rm pr}_1,\hat{\mathbold{E}},X)$,
for which the complex curve  $\hat{\mathbold{E}}$ is the union of four disjoint complex curves $\hat{\mathbold{E}}^k,\, k=1,2,3,4,$ each intersecting each fiber $\{\hat{x}\}\times\mathbb{P}^1$ along a single point $(\hat{x},\,\hat{g}^k(\hat{x}))$. %Indeed,
This can be seen as follows. Let $q_0$ be the base point of $X$.
The monodromy mapping class along each element $e$ of $\pi_1(X,q_0)$ takes the set of distinguished points ${\mathbold{E}}\cap (\{q_0\}\times\mathbb{P}^1)$
onto itself, permuting them by a permutation $\sigma(e)$. Consider the set $ N$ of
elements $e\in \pi_1(X,q_0)$ for which $\sigma(e)$ is the identity. The set $N$ is a normal subgroup of $\pi_1(X,q_0)$. Its index is finite, since two left cosets $e_1\,N$
and $e_2 \, N$ are equal if $ \sigma(e_2\,e_1^{-1}) =    \sigma(e_2)\sigma(e_1)^{-1}={\rm Id}$, and there are only finitely many distinct permutations of points of ${\mathbold{E}}\cap (\{q_0\}\times\mathbb{P}^1)$. The quotient $\hat{X}\stackrel{def}=\tilde {X}\diagup {\rm Is}^{\tilde{q}_0}(N)$ of the universal covering of $X$ by the group of covering transformations corresponding to $N$ and the canonical projection $\hat{X}\to X$ define the required covering.

To prove the corollary, we have to show, that any reducible holomorphic $(1,1)$-bundle over a punctured Riemann surface $X$ is locally holomorphically trivial, and that two isotopic (equivalently, smoothly isomorphic) holomorphically non-trivial holomorphic $(1,1)$-bundles over $X$ are holomorphically isomorphic.

The second fact is obtained as follows.
Suppose the holomorphically non-trivial holomorphic $(1,1)$-bundles $\mathfrak{F}_j,\,j=1,2,$ have conjugate monodromy homomorphisms. By Proposition \ref{prop1}
each $\mathfrak{F}_j$ is holomorphically isomorphic to a double branched covering of a special holomorphic $(0,4)$-bundle
$(X\times\mathbb{P}^1,{\rm pr}_1,\mathbold{E}_j,X)\stackrel{def}= {P}(\mathfrak{F}_j)$. The bundles ${P}(\mathfrak{F}_j)$ are isotopic, since they have conjugate monodromy homomorphisms. There is a finite unramified covering $\hat{{\sf P}}:\hat X\to X$ of $X$, such that the bundles  $ P(\mathfrak{F}_j)$ have isotopic lifts
$(\hat{X}\times\mathbb{P}^1,{\rm pr}_1,\hat{\mathbold{E}}_j,X)$ to $\hat X$, and for each $j$ the complex curve  $\hat{\mathbold{E}}_j$ is the union of four disjoint complex curves $\hat{\mathbold{E}}_l^k,\, k=1,2,3,4,$ each intersecting each fiber $\{\hat{x}\}\times\mathbb{P}^1$ along a single point $(\hat{x},\,\hat{g}_j^k(\hat{x}))$.
The lifted bundles are not isotopic to the trivial bundle. The mappings $\hat{X}\ni \hat{x}\to \hat{g}_j^k(\hat{x})$ are holomorphic.
We may assume that $\hat{g}_j^4(\hat{x})=\infty$ for each $\hat{x}$.
Define for $j=1,2,$ a holomorphic isomorphism of the bundle $(\hat{X}\times\mathbb{P}^1,{\rm pr}_1,\hat{\mathbold{E}}_j,X)$ by
$$
\{\hat{x}\}\times \mathbb{P}^1\ni (\hat{x},\,\zeta)\to \Big(\hat{x},\,-1 + 2\frac{\hat{g}_j^1(\hat{x})-\zeta}{\hat{g}_j^1(\hat{x})-\hat{g}_j^2(\hat{x}) }\Big)\,.
$$

The image $\hat{\mathbold{E}}'_j$ of $\hat{\mathbold{E}}_j$ under the $j$-th isomorphism
intersects the fiber over each $\hat{x}\in \hat X$ along the four points $(\hat{x},-1),\,(\hat{x},1),\,(\hat{x},\infty),$ and $(\hat{x},\mathring{g}_j(\hat{x}))$ for a holomorphic function $\mathring{g}_j$ on $\hat X$
that avoids $-1,\,1$ and $\infty$. The functions $\;\mathring{g}_j\;$, $j=1,2,\;$ are homotopic, since the bundles
are isotopic. They are not homotopic to a constant function since the bundles are not isotopic to the trivial bundle.
By Corollary \ref{cor1a} the functions $\;\mathring{g}_1\;$ and $\;\mathring{g}_2\;$ coincide.
Hence, the bundles $\;\;(\hat{X}\times\mathbb{P}^1,{\rm pr}_1,\hat{\mathbold{E}}_j,X)\;\;$
are holomorphically isomorphic. This means that there is a nowhere vanishing holomorphic function $\hat{\alpha}$ on $\hat X$, such that for each $\hat{x}\in \hat{X}$ the equality
$\{\hat{x}\}\times\hat{E}_2(\hat{x})=\{\hat{x}\}\times\hat{\alpha}(\hat{x}) \hat{E}_1(\hat{x})$ holds. Here $\hat{E}_j(\hat{x})$ is defined by the equality $\hat{\mathbold{E}}_j\cap (\{\hat{x}\}\times \mathbb{P}^1)=\{\hat{x}\}\times \hat{E}_j(\hat{x})$. Define also $E_j(x)$ by the equality $\{x\}\times E_j(x)=\mathbold{E}_j\cap (\{{x}\}\times \mathbb{P}^1)$.
For a point $x\in X$ and $\hat{x}_1,\hat{x}_2\in {\hat{{\sf P}}}^{-1}(x)$
the equalities $\hat{{E}}_j(\hat{x}_1)= \hat{{E}}_j(\hat{x}_2)=E_j(x),\,j=1,2,$ hold. Hence,
$E_2(x)=\hat{\alpha}(\hat{x}_1){E}_1({x})=\hat{\alpha}(\hat{x}_2) {E}_1({x})$.
For a set $E\subset C_3(\mathbb{C})\diagup \mathcal{S}_3$ and a complex number $\alpha$ the equality $E=\alpha E$ is possible only if $\alpha=1$, or $\alpha= -1$ and $E$ is obtained from $\{-1,0,1\}$ by multiplication with a non-zero complex number, or $\alpha=e^{\pm \frac{2\pi i}{3}}$ and $E$ is obtained from the set of vertices of an equilateral triangle with barycenter $0$ by multiplication with a non-zero complex number.

For $x$ in a small open disc on $X$ and $x\to\hat{x}_j(x),\,j=1,2,$ being two local inverses of $\hat{{\sf P}}$ the functions $x\to \hat{\alpha}(\hat{x}_j(x))$ are two analytic functions whose ratio is contained in a finite set, hence the ratio is equal to a constant. If the constant was different from one, then all fibers of ${ P}(\mathfrak{F}_1)$ would be conformally equivalent to each other and, hence, ${ P}(\mathfrak{F}_1)$ would be locally holomorphically trivial.
Since the bundles $\mathfrak{F}_j$, and, hence, also the ${P}(\mathfrak{F}_j)$, are locally holomorphically non-trivial, the ratio  of the two functions equals $1$.
We saw that for each pair of points $\hat{x}_1,\hat{x}_2\in \hat{X}$, that project to the same point $x\in X$, $\hat{\alpha}(\hat{x}_1)=\hat{\alpha}(\hat{x}_2)$. Put $\alpha(x)=\hat{\alpha}(\hat{x}_j)$ for any point $\hat{x}_j\in (\hat{{\sf P}})^{-1}(x)$.
We obtain
${E}_2({x})={\alpha}({x}){E}_1({x})$, that means,
the bundles ${P}(\mathfrak{F}_j)$ are holomorphically isomorphic. Since the bundles $\mathfrak{F}_j$, $j=1,2,$ are double branched coverings of the ${P}(\mathfrak{F}_j)$ and have conjugate monodromy homomorphism, they are holomorphically isomorphic.

The first fact is obtained as follows. After a holomorphic isomorphism we may assume %from the beginning
that the reducible holomorphic $(1,1)$-bundle is a double branched covering of a reducible special $(0,4)$-bundle $P(\mathfrak{F})=(X\times \mathbb{P}^1, {\rm pr}_1, \mathring{\mathbold{E}}\cup \mathbold{s}^{\infty}, X)$. After a further isomorphism the bundle $P(\mathfrak{F})$ lifts to a holomorphic bundle $\reallywidehat{P(\mathfrak{F})}= (\hat{X} \times \mathbb{P}^1, {\rm pr}_1, \hat{\mathring{\mathbold{E}}}\cup \widehat{\mathbold{s}^{\infty}},\hat{ X}) $, such that
$\hat{\mathring{\mathbold{E}}}$ intersects each fiber $\{x\}\times\mathbb{P}^1$ along a set of the form $\{\hat{x}\}\times \{-1,1,\mathring{g}(\hat{x}))$. Since $\mathfrak{F}$ is reducible, hence $P(\mathfrak{F})$ and also $\reallywidehat{P(\mathfrak{F})}$ are reducible, the mapping $\mathring{g}$ is constant by Corollary \ref{cor1a}. Hence, all fibers of $\reallywidehat{P(\mathfrak{F})}$ are conformally equivalent, and, hence, all fibers of $P(\mathfrak{F})$ are conformally equivalent. Since $\mathfrak{F}$ is the double branched covering of $P(\mathfrak{F})$, all fibers of $\mathfrak{F}$ are conformally equivalent.
The first fact is proved.
\hfill $\Box$

\medskip

\noindent {\bf Proof of Proposition \ref{prop1a}.}
Denote by $S^{\alpha}$ a skeleton of $T^{\alpha,\sigma} \subset T^{\alpha}$ which is the union of two circles each of which lifts
under the covering ${\sf P}:\mathbb{C}\to T^{\alpha}$ to
a straight line segment which is parallel to an axis in the complex plane. Denote the intersection point of the two circles by $q_0$.
Note that $S_{\alpha}$ is a standard bouquet of circles for $T^{\alpha,\sigma}$ with base point $q_0$, and
${\sf P}^{-1}(T^{\alpha,\sigma})$ is the $\frac{\sigma}{2}$-neighbourhood of
${\sf P}^{-1}(S^{\alpha})$.

Denote by $e$ the generator of $\pi_1(T^{\alpha,\sigma},q_0)$, that lifts to a vertical line segment and $e'$  the generator of $\pi_1(T^{\alpha,\sigma},q_0)$, that lifts to a horizontal line segment. Put $\mathcal{E}=\{e,e'\}$. We show first the inequality
\begin{equation}\label{eq20}
\lambda_3(T^{\alpha,\sigma}) \leq \frac{4(2\alpha+1)}{\sigma}\;.
\end{equation}
For this purpose we take any primitive element $e''$ of the fundamental group $\pi_1(T^{\alpha,\sigma},q_0)$ which is the product of at most three factors, each of the factors being an element of $\mathcal{E}$ or the inverse of an element of $\mathcal{E}$.
We represent
the element $e''$ by a piecewise $C^1$ mapping $f_1$ from an interval $[0,l_1]$
to the skeleton $S^{\alpha}$.
We may consider $f_1$ as a piecewise $C^1$ mapping from the circle
$\mathbb{R}\diagup (x \sim x+l_1)$ to the skeleton, and assume that for all points $t'$ of the circle where $f_1$ is not smooth, $f_1(t')=q_0$.
Let $t_0\in [0,l_1]$ be a point for which $f_1(t_0)\neq q_0$.
Let $\tilde{f}_1$ be a piecewise smooth mapping from $[t_0,t_0+l_1]$ to the universal covering $\mathbb{C}$ of $T^{\alpha}\subset T^{\alpha,\sigma}$
which projects to $f_1$.
We may take $f_1$ so that the equality $|\tilde{f}_1'|=1$ holds. The mapping may be chosen so that $l_1\leq 2 \alpha  + 1$.
(Recall that $\alpha\geq 1$ and the element $e$ is primitive.)

Take any $t'$ for which $f_1$ is not smooth. We may assume that $f_1$ is chosen so that
the direction of $\tilde{f}_1'$ changes by the angle $\pm \frac{\pi}{2}$ at each such point. Hence,
there exists a neighbourhood $I(t')$ of $t'$ on $(t_0,t_0+\ell_1)$, such that the restriction $\tilde{f}_1'| I(t')$
covers two sides of a square of side length $\frac{\sigma}{2}$. Denote $\tilde{q}'_0$  the common vertex  $\tilde{f}_1'(t')$ of these sides, and
by  $\tilde{q}''_0$ the vertex of the square that is not a vertex of one of the two sides.
Replace the union of the two sides of the square that contain $\tilde{q}'_0$ by a quarter-circle of radius $\frac{\sigma}{2}$ with center at the vertex $\tilde{q}''_0$, and
parameterize the latter by $t \to \frac{\sigma}{2} e^{\pm i\frac{2}{\sigma}t}$ so that the absolute value of the derivative equals $1$.
Notice that the quarter-circle is shorter than the union of the two sides.

Proceed in this way with all such points $t'$. After a reparameterization
we obtain a $C^1$ mapping $\tilde f$ of the interval $[0,l]$ of length $l$ not exceeding $2 \alpha +1$ whose image is contained in the union of
${\sf P}^{-1}(S^{\alpha})$
with some quarter-circles, such that $|\tilde{f}'|=1$.
The distance of each point of the image of $\tilde f$ to the boundary of ${\sf P}^{-1}(T^{\alpha,\sigma})$ %$\widetilde{T^{\alpha,\sigma}}$
is not smaller than $\frac{\sigma}{2}$. The mapping $\tilde f$ is piecewise of class $C^2$.
The normalization condition  $|\tilde{ f}'|=1$ implies $|\tilde{f}''|\leq \frac{2}{\sigma}$.

The projection $f={\sf P}\circ\tilde{f} $ can be considered as a mapping from the circle $\mathbb{R}\diagup (x\sim x +l)$ of length $l$ not exceeding $2 \alpha + 1$  to $T^{\alpha,\sigma}$,  that represents
the free homotopy class $\reallywidehat{e''}$ of the chosen element of the fundamental group.

Consider the mapping $ x+iy \to \tilde{F}(x+iy) \stackrel{def}=\tilde{f}(x) + i \tilde{f}'(x) y \in \mathbb{C}$, where
$x+iy$ runs along the rectangle $R_l= \{x+iy  \in \mathbb{C}:
x\in [0,l], |y| \leq \frac{\sigma}{4}\}$.
The image of this mapping is contained in the closure of ${\sf P}^{-1}(T^{\alpha,\sigma})$.
Since $2\frac{\partial}{\partial z} \tilde{F} (x+iy)= 2 \tilde{f}'(x) + i \tilde{f}''(x) y $ and $2\frac{\partial}{ \partial \bar z} \tilde{F }(x+iy)=  i \tilde{f}''(x) y $
the Beltrami coefficient $\mu
_{\tilde{F}}(x+iy)= \frac{\frac{\partial}{\bar z} \tilde{F} (x+iy)}{\frac{\partial}{z} \tilde{F} (x+iy)}$ of $\tilde{F}$  satisfies the inequality $|\mu
_{\tilde{F}}(x+iy) | \leq \frac{1}{3}$.  Hence, for $K=\frac{1+\frac{1}{3}}{1-\frac{1}{3}}=2$ the mapping $\tilde{F}$ descends to a $K$-quasiconformal mapping $F$ from the annulus $A_l$ to $T^{\alpha,\sigma}$ of extremal length $\lambda(A_l)=\frac{l}{\frac{\sigma}{2}}\leq 2 \frac{(2 \alpha +1)}{\sigma}$
that represents the free homotopy class of the element $e''$ of the fundamental group $\pi_1(T^{\alpha,\sigma} ,q_0)$. Realize $A_l$ as an annulus in the complex plane. Let $\varphi$ be the solution of the Beltrami equation on $\mathbb{C}$ with Beltrami coefficient $\mu_{\tilde{F}}$ on $A_l$ and zero else.  Then the mapping $g=F \circ \varphi^{-1}$ is a holomorphic mapping of the annulus $\varphi(A_l)$ of extremal length not exceeding $K \lambda(A_l) \leq  \frac{4(2\alpha +1 )}{\sigma} $ into $T^{\alpha,\sigma}$ that represents the chosen element of the fundamental group $\pi_1(T^{\alpha,\sigma},q_0)$. Inequality \eqref{eq20} is proved.

By Theorem \ref{thm1}
for tori with a hole there are up to homotopy no more than $3(\frac{3}{2}e^{24 \pi \lambda_3(T^{\alpha,\sigma})})^2\leq \frac{27}{4}e^{3\cdot 2^4 \pi \frac{2\alpha +1}{\sigma}}< 7e^{3\cdot 2^4 \pi \frac{2\alpha +1}{\sigma}} $
non-constant irreducible holomorphic mappings from $T^{\alpha,\sigma}$ to the twice punctured complex plane.

We give now the proof of the lower bound. Let $\delta=\frac{1}{10}$.
We consider the annulus $A^{\alpha,\delta}\stackrel{def}=\{z\in \mathbb{C}:|\mbox{Re}z|< \frac{5\delta}{2}\}\diagup (z\sim z+\alpha i)$.
The extremal length of the annulus equals $\frac{\alpha}{5\delta}=2\alpha$.

For any natural number $j$ we consider all elements of $\pi_1(\mathbb{C}\setminus \{-1,1\},0)$ of the form
\begin{equation}\label{eq21}
a_1^{\pm 2}a_2^{\pm 2} \ldots a_1^{\pm 2}a_2^{\pm 2}
\end{equation}
containing $2j$ terms, each of the form $a_j^{\pm 2}$. The choice of the sign in the exponent of each term is arbitrary. There are $2^{2j}$ elements of this kind.  By \cite{Jo2} there is a relatively compact domain $G$ in the twice punctured complex plane $\mathbb{C}\setminus \{-1,1\}$ and a positive constant $C$ such that the following holds. For each $j$, each element of the fundamental group of the form \eqref{eq21}, and for each annulus of extremal length at least $2Cj$ there exists a base point $q$ in the annulus, and a holomorphic mapping from the annulus to $G$ that maps $q$ to $0$ and represents
the element.
Put $j=[\frac{\alpha}{10C\delta}]$, where $[x]$ is the largest integer not exceeding a positive number $x$. Then each element of the form \eqref{eq21}
with this number $j$ can be represented by a holomorphic map $\mathring{g}$ from the annulus $A^{\alpha,\delta}$ to $G$. There is a constant $C_1$ that depends only on $G$ such that
the mapping $\mathring{g}$ satisfies the inequality $|\mathring{g}|<C_1$. Let $g$ be the lift of $\mathring{g}$ to a mapping from the strip $\{z\in \mathbb{C}:|\mbox{Re}z|< \frac{5\delta}{2}\}$ to $G$.
On the thinner strip
$\{|\mbox{Re}z|< \frac{3\delta}{2}\}$ the derivative of $g$
satisfies the inequality $|g'|\leq \frac{C_1}{\delta}$.

We will associate to the holomorphic mapping $\mathring{g}$ on the annulus a smooth mapping $g_1$ from $T^{\alpha,\delta}\subset T^{\alpha}$ to $G$, such that (with $\sf P$ being the projection ${\sf P}:\mathbb{C}\to T^{\alpha}$) the monodromy
along the circle ${\sf P}(\{\mbox{Re}z=0\})$
with base point ${\sf P}(0)$
is equal to \eqref{eq21}, and the monodromy along  ${\sf P}(\{\mbox{Im}z=0\}$
with the same base point equals the identity. This is done as follows.
Let $F_{\alpha}=[-\frac{1}{2},\frac{1}{2})\times [-\frac{\alpha}{2},\frac{\alpha}{2})\subset \mathbb{C}$ be a fundamental domain for the projection ${\sf P}:\mathbb{C}\to T^{\alpha}$. Put $\Delta^{\alpha,\delta}=F_{\alpha}\cap {\sf P}^{-1}(T^{\alpha,\delta})$.
Let $\chi_0:[0,1]\to \mathbb{R}$ be a non-decreasing function of class $C^2$ with $\chi_0(0)=0,\, \chi_0(1)=1$, $\chi'_0(0)=\chi'_0(1)=0$ and $|\chi_0'(t)|\leq \frac{3}{2}$. Define $\chi:[\frac{-3\delta}{2},\frac{+3\delta}{2}]\to [0,1]$ by
\begin{equation}\label{eq22}
\chi(t) =
\begin{cases} \chi_0(\frac{1}{\delta}t+ \frac{3}{2})& \; t \in [\frac{-3\delta}{2},\frac{-\delta}{2}] \\
1 & \; t \in [\frac{-\delta}{2},\frac{+\delta}{2}]\\
\chi_0(-\frac{1}{\delta}t+ \frac{3}{2})& \; t \in [\frac{\delta}{2},\frac{3\delta}{2}]\,.
\end{cases}
\end{equation}
Notice that $\chi$ is a $C^2$-function that vanishes at the endpoints of the interval  $[\frac{-3\delta}{2},\frac{+3\delta}{2}]$ together with its first derivative, is non-decreasing on $[\frac{-3\delta}{2},\frac{-\delta}{2}]$, and non-increasing on $[\frac{\delta}{2},\frac{3\delta}{2}]$.
Put
$g_1(z)= \chi (\mbox{Re}z)\; g(z) + (1-\chi (\mbox{Re}z))\; g(0)$ for $z$ in the intersection of
$\Delta^{\alpha,\delta}$ with $\{|\mbox{Re}z|<\frac{3\delta}{2}\}$,
and $g_1(z)=g(0)$ for $z$ in the rest of $\Delta^{\alpha,\delta}$.

Put $\varphi(z)= \frac{\partial}{\partial \bar z} g_1(z)$ on $\Delta^{\alpha,\delta}$.
Since $\frac{\partial}{\partial \bar z} \chi (\mbox{Re} z)=0$ for $|\mbox{Re} z|<\frac{\delta}{2}$ and for $|\mbox{Re} z|>\frac{3\delta}{2}$,
the function  $\varphi(z)$ vanishes on $\Delta^{\alpha,\delta}\setminus Q$ with $Q\stackrel{def}=([ -\frac{3\delta}{2},  +\frac{3\delta}{2}] \times [-\frac{\delta}{2},\frac{\delta}{2}])$.
On $Q \cap \Delta^{\alpha,\delta}$ the inequality
\begin{equation}\label{eq23}
|\varphi(z)|\leq\frac{1}{2} |\chi'(\mbox{Re} z)|\, |g(z)-g(0)|\leq \frac{3}{4\delta} \cdot \frac{C_1}{\delta}|z|< \frac{3}{4\delta^2} \cdot  C_1 \cdot 2\delta= \frac{3}{2} \frac{C_1}{\delta}\,
\end{equation}
holds. Notice that the functions $g_1$ and $\varphi$ extend to ${\sf P}^{-1}( T^{\alpha,\delta})$ as continuous doubly periodic functions. Hence, we may consider them as functions on $T^{\alpha,\delta}$.

\begin{figure}[h]
\begin{center}
\includegraphics[width=55mm]{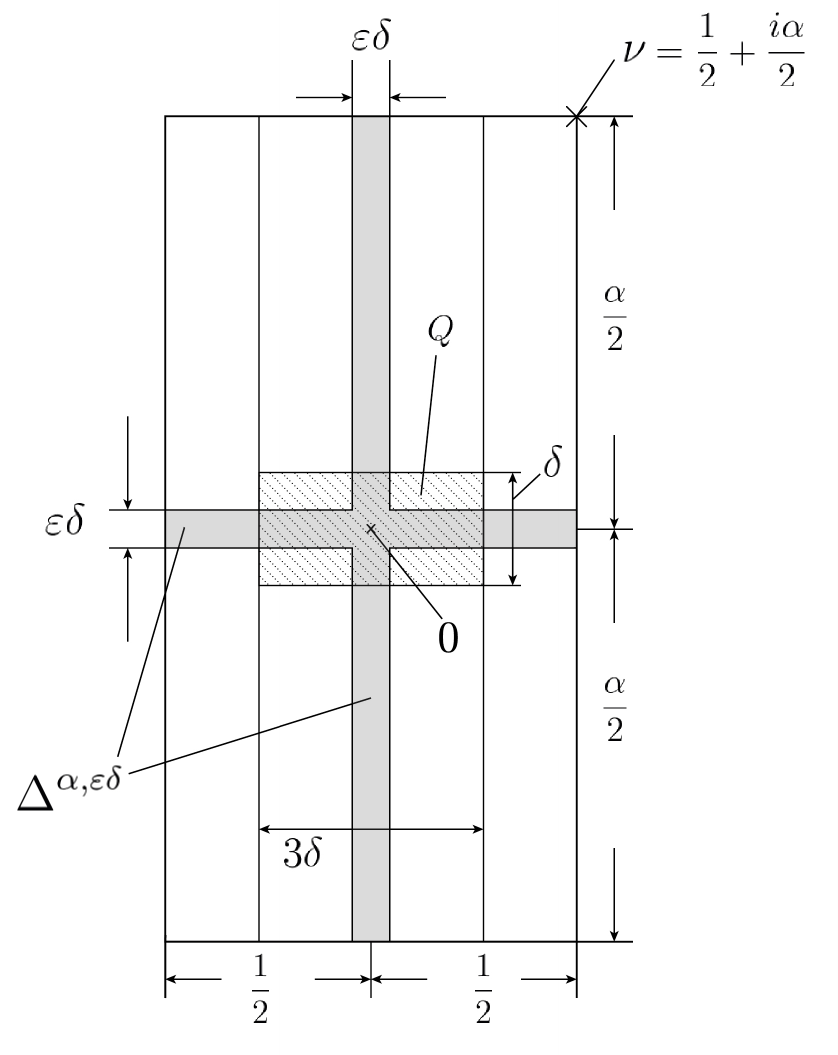}
\end{center}
\caption{A fundamental domain for a torus with a hole and the poles of the kernel for the $\overline{\partial}$-equation}\label{fig8.4}
\end{figure}

We want to find a small positive number $\varepsilon$ that depends on $C$ and $C_1$, but not on $\alpha$, such that the following holds. For $\sigma\stackrel{def}= \varepsilon \delta$ there exists a solution $f$ of the equation $\frac{\partial}{\partial \bar z}f(z)=\varphi(z)$ on $T^{\alpha,\sigma}$ such that for each $z$ the value $|f(z)|$ is smaller than the Euclidean distance in $\mathbb{C}$ of $\pm 1$ to $\overline G$.
Then $g_1-f$ is a holomorphic mapping from $T^{\alpha,\sigma}$ to
$\mathbb{C}\setminus \{-1,1\}$ whose class has monodromies equal to \eqref{eq21}, and to the identity, respectively. %There are at least
This gives
$\,2^{2 [\frac{\alpha}{10C\delta}]}\geq      2^{2 \frac{\alpha}{10C\delta}-2}=\frac{1}{4}e^{\frac{2\varepsilon\log 2}{10C}\frac{\alpha}{\sigma}}$ different homotopy classes of mappings from $T^{\alpha,\sigma}$ to $\mathbb{C}\setminus \{-1,1\}$, and, hence proves the lower bound.

To solve the $\bar{\partial}$-problem on $T^{\alpha,\varepsilon\delta}=  T^{\alpha,\sigma}$, we consider an explicit kernel function which mimics the Weierstraß $\wp$-function. The author is grateful to Bo Berndtsson who suggested to use this kernel function.

Recall that the Weierstraß $\wp$-function related to the torus $T^{\alpha}$ is the doubly periodic meromorphic function
$$
\wp_{\alpha}(\zeta)=\frac{1}{\zeta^2}+ \underset{(n,m) \in \mathbb{Z}^2\setminus (0,0)}{\sum} \Big(\frac{1}{(\zeta-n- i m \alpha)^2 }- \frac{1}{(n+ i m \alpha)^2}\Big)\,
$$
on $\mathbb{C}$. It defines a meromorphic function on $T^{\alpha}$ with a double pole at the projection of the origin and no other pole.

Put $\nu=\frac{1}{2} + \frac{\alpha i}{2}$.
Since for $\zeta\not\in (\mathbb{Z}+i\alpha \mathbb{Z})\cup  (\nu+\mathbb{Z}+i\alpha \mathbb{Z})$
%$n+i m \alpha\neq 0$
the equality
\begin{align*}
\frac{1}{(\zeta-n- i m \alpha) } - \frac{1}{(\zeta-n- i m \alpha-\nu) }+ \frac{\nu}{(\zeta-n- i m \alpha)^2 }
= \frac{-\nu^2}{(\zeta -n- i m \alpha)^2 (\zeta - n- i m \alpha-\nu) }\,
\end{align*}
holds, and the series with these terms converges uniformly on compact sets not containing poles, the expression
$$
\wp_{\alpha}^{\nu}(\zeta)= \frac{1}{\zeta}-\frac{1}{\zeta-\nu} + \,
\underset{(n,m) \in \mathbb{Z}^2\setminus (0,0)}{\sum} \Big(\frac{1}{(\zeta-n- i m \alpha) } - \frac{1}{(\zeta-n- i m \alpha-\nu) }     + \frac{\nu}{(n+ i m \alpha)^2}\Big)\,
$$
defines a doubly periodic meromorphic function on $\mathbb{C}$ with only simple poles. The function descends to a meromorphic function on $T^{\alpha}$ with two simple poles and no other pole.

Recall that the support of $\varphi$ is contained in $Q$. The set $Q$
is contained in the $2\delta$-disc in $\mathbb{C}$ (in the Euclidean metric) around the origin.
If $\zeta$ is contained in the $2\delta$-disc around the origin and $z \in \Delta^{\alpha,\delta}$, then the point $\zeta-z$ is contained in the ${2\delta}$-neighbourhood (in $\mathbb{C}$) of $\Delta^{\alpha,\delta}$.
By the choice of $\delta$ the distance of any such point  $\zeta-z$
to any lattice point $n+i \alpha m$ except $0$ is larger than  $\frac{1}{2}-2\delta>  \frac{1}{4}$.
Further, for $z\in \Delta^{\alpha,\delta}$ and $\zeta$ in the $2\delta$-disc around
the origin the distance of the point  $\zeta-z$
to any point $n+i \alpha m +\nu$ (including the point $\nu$) is not smaller than $ \frac{1}{2}-\frac{5\delta}{2}=  \frac{1}{4}$.
Put $Q_{\varepsilon}\stackrel{def}= Q\cap  \Delta^{\alpha,\varepsilon\delta}=      ([  -\frac{3\delta}{2},  +\frac{3\delta}{2}] \times [-\frac{\varepsilon\delta}{2},+\frac{\varepsilon\delta}{2}]) \bigcup ([  -\frac{\varepsilon\delta}{2},  +\frac{\varepsilon\delta}{2}] \times [-\frac{\delta}{2},+\frac{\delta}{2}])$.
Then the function
\begin{align}\label{eq24}
f(z)= - \frac{1}{\pi}\iint_{Q_{\varepsilon}} \varphi(\zeta) \wp_{\alpha}^{\nu} (\zeta-z) dm_2(\zeta)\;,
\end{align}
for $z$ in $\Delta^{\alpha,\varepsilon\delta}$ is holomorphic outside  $Q_{\varepsilon}$
and satisfies the equation $\frac{\partial}{\partial \bar z}f=\varphi$ on $Q_{\varepsilon}$.
It extends continuously to
a doubly periodic function on ${\sf P}^{-1}(T^{\alpha,\varepsilon\delta})$ and hence descends to a continuous function on $T^{\alpha,\varepsilon\delta}$.
It remains to estimate the supremum norm of the function $f$ on  $\Delta^{\alpha,\sigma}=\Delta^{\alpha,\varepsilon\delta}$.
The following inequality holds for  $z \in \Delta^{\alpha,\sigma}$
\begin{align}\label{eq26}
|\iint_{Q_{\varepsilon}} \varphi(\zeta) \wp_{\alpha}^{\nu} (\zeta-z) dm_2(\zeta)|=
|& \frac{1}{\pi}\iint_{Q_{\varepsilon}} \varphi(\zeta)\Big(\frac{1}{\zeta-z} + (\wp_{\alpha}^{\nu}(\zeta-z) -\frac{1}{\zeta-z})\Big)  dm_2(\zeta)| \nonumber\\
\leq & \frac{1}{\pi}\iint_{Q_{\varepsilon}} \frac{3C_1}{2\delta}\Big(|\frac{1}{\zeta-z}|+ C_2\Big)   dm_2(\zeta)  \,.
\end{align}
We used the upper bound for $\varphi$ and the fact that for $z \in \Delta^{\alpha,\sigma}$ and $\zeta$ in $Q_{\varepsilon}$ the expression $|\wp_{\alpha}^{\nu}(\zeta-z)  -\frac{1}{\zeta-z}|$ is bounded by a universal constant $C_2$.
The integral  of the second term on the right hand side does not exceed $\frac{3 C_1}{2 \delta}\cdot C_2 \cdot 4\varepsilon \delta^2=6 C_1 C_2 \varepsilon \delta$.
The integral $\quad \iint_{Q_{\varepsilon}}|\frac{1}{\zeta-z}|dm_2(\zeta)\; $ does not exceed the sum of the two integrals
$\quad I_1=\iint_{(- \frac{3}{2}  \delta, \frac{3}{2} \delta)\times (- \frac{1}{2}\varepsilon\delta, \frac{1}{2}\varepsilon\delta)}\; \mid \frac{1}{\zeta-z}\mid dm_2(\zeta)\;, \;\; $
and
$\;\; I_2=\iint_{(-\frac{1}{2}\varepsilon\delta, \frac{1}{2}\varepsilon\delta)\times (-\frac{1}{2}\delta,\frac{1}{2}\delta)}\; \mid \frac{1}{\zeta-z}\mid  dm_2(\zeta)\;$.
The first integral $I_1$ is largest when $z= 0 $. Hence, it does not exceed
\begin{align}\label{eq27}
\iint_{|\zeta|< (\sqrt{2})^{-1} \varepsilon\delta} \; \mid \frac{1}{\zeta} \mid dm_2(\zeta) + 2\varepsilon \delta \int_{\frac{1}{2}\varepsilon\delta}^{\frac{3}{2}\delta}\, \frac{1}{\eta}\, d\eta
\leq  \sqrt{2}\pi\varepsilon\delta + 2\varepsilon\delta \log{\frac{3}{\varepsilon}}\,.
\end{align}
The second integral $I_2$ is smaller.
We obtain the estimate
\begin{equation}\label{eq28}
|f(z)|\leq  \frac{6 C_1 C_2 \varepsilon \delta}{\pi} + 3\frac{C_1}{\pi\delta}(\sqrt{2}\pi\varepsilon\delta + 2\varepsilon\delta \log{\frac{3}{\varepsilon}})\,.
\end{equation}
Recall that we have chosen $\delta=\frac{1}{10}$.
We may choose $\varepsilon_0>0$ depending only on $C_1$ (and, hence, only on the domain $G$) so that if $\varepsilon< \varepsilon_0$
%$\sigma<\frac{1}{2}$
the supremum norm of $f$ is less than the distance of $\pm 1$ to $\overline G$. The proposition is proved.
\hfill  $\Box$

\medskip

\noindent {\bf Proof of Proposition \ref{prop1b}.} Let $\ell_0$ be the length in the K\"ahler metric of the longest circle in the bouquet. For each natural number $k$ and each positive $\sigma<\sigma_0$ the value $\lambda_k(S_{\sigma})$ satisfies the inequalities
\begin{align}
C_1' \frac{\ell_0}{\sigma}\leq \lambda_k(S_{\sigma})\leq C_1'' \frac{\ell_0}{\sigma}
\end{align}
for constants $C_1'$ and $C_1''$ depending on $k$, $X$, $S$ and the K\"ahler metric. This can be seen by the argument used in the proof of Proposition \ref{prop1a}.

The upper bound in inequalities \eqref{eqabc} follows from Theorem \ref{thm1}.

The proof of the lower bound in \eqref{eqabc} follows along the same lines as the proof of Proposition \ref{prop1a}. It leads to a $\overline{\partial}$-problem on an open Riemann surface, for which H\"ormander's $L^2$-method can be used. The case of open Riemann surfaces is easier to treat as the general case of pseudo-convex domains. The needed results for Riemann surfaces are explicitly formulated in \cite{Na}.

To obtain the lower bound we
consider for each positive number $\delta<\sigma_0$ the $\delta$-neighbourhood of the longest
%(in the K\"ahler metric)
circle $\gamma_0$ of the bouquet. Consider a curvilinear rectangle $R^X_{\delta}$,
%of extremal length at least $\frac{3}{4}\lambda_0$,
that is contained in the $\delta$-neighbourhood of
the largest circle $\gamma_0$,
%of the bouquet,
whose ''vertical curvilinear sides'' are contained in the boundary of $S_{\delta}$ and whose open ''horizontal curvilinear sides'' are contained in $S_{\delta}$. Choosing $\sigma_0$ small enough, we may choose
$R^X_{\delta}$  so that for its extremal length the inequality $\lambda(R^X_{\delta})>c \frac{\ell_0}{\delta} +4$ holds for a number $c>0$ that depends only on $X$, $S$ and the K\"ahler metric.
For any positive $\delta<\sigma_0$ we denote by $R_{\delta}$ the true rectangle
$R_{\delta}\stackrel{def}=\{x+iy: x\in (-\delta,\delta), y \in (-c\ell_0-2\delta,c\ell_0+2\delta)\} $
in the complex plane.
Shrinking perhaps $R^X_{\delta}$, we may assume that $R^X_{\delta}$
is conformally equivalent to $R_{\delta}$. Denote by $\omega$ the conformal mapping $R^X_{\delta}\to R_{\delta}$ for which the orientation of the curve $\gamma_0$  corresponds to the
positive orientation of the imaginary axis.

Let $G\subset \mathbb{C}\setminus \{-1,1\}$  be the same relatively compact domain as in the proof of Proposition \ref{prop1a}, and let
$\mathring{R}_{\delta}\subset R_{\delta} $ be the rectangle in the complex plane with the same center and horizontal side length as $R_{\delta} $, and with
vertical side length $2c\ell_0$.
There is an absolute constant $C>0$ such that for  $j=[\frac{c}{C} \frac{\ell_0}{\delta}]$ and any word of the form $a_1^{\pm 1}\, a_2^{\pm}\ldots a_2^{\pm 1}$ in the relative fundamental group  $\pi_1(\mathbb{C}\setminus\{-1,1\}, (-1,1))$ with $2j$  terms
there exists a holomorphic mapping $g:\mathring{R}_{\delta}\to G\subset \mathbb{C}\setminus \{-1,1\}$
%whose restriction to $\mathring{R}_{\delta}$
that represents this word and vanishes at $ \pm i c\ell_0$ (see Theorem 1 of \cite{Jo2}).
The function $g$ extends by reflection through the horizontal sides of $\mathring{R}_{\delta}$ to a holomorphic function on $R_{\delta}$, that we also denote by $g$.
%Let $C_3$ be the distance in the Euclidean metric of the
%domain $G$ from the tuple of points $\{-1,1,\infty\}$. Then
Since $|g|\leq C_1$ on $\mathring{R}_{\delta}$ and, hence, $|g|\leq C_1$ also on $R_{\delta}$, for any positive $\alpha<1$ the inequality $|g'|\leq  \frac{C_1}{\delta(1-\alpha)}$ holds for the derivative of the mapping $g$ on the smaller rectangle $R_{\delta\alpha}$ (defined as $R_\delta$ with $\delta$ replaced by $\alpha\delta$). This fact implies that $|g|\leq \frac{\sqrt{2} C_1\alpha}{1-\alpha}$ on $Q_{\alpha\delta}^{\pm}\stackrel{def}=
\{x+iy: x\in (-{\alpha\delta},{\alpha\delta}), \pm y\in
(c\ell_0, c\ell_0 + \alpha \delta)\}$. We took into account that
%the function $g$ on this set is obtained by reflection from
%the restriction of $g$ to a subset of $R_{\alpha \delta}$ and
$g(\pm i c \ell_0)=0$.
We take $\alpha$ so that  $\frac{\sqrt{2} C_1\alpha}{1-\alpha}
=\frac{1}{2}$.
%$\frac{2C_3\alpha}{\delta(1-\alpha)}$

With the same function $\chi_0$ as in the proof of Proposition \ref{prop1a} we define
%$\chi(t)=\chi_0(\pm \frac{t-c\ell_0}{\alpha\delta})$
\begin{equation}\label{eq22+}
\chi(t) =
\begin{cases}
1 & \; \;t \;\in [-c\ell_0,c\ell_0]\\
\chi_0( \frac{c\ell_0+ \alpha\delta -|t|}{\alpha\delta}) &\; |t| \in (c\ell_0, c\ell_0+ \alpha \delta)\,,\\
0&\; \;t\; \in \mathbb{R}\setminus[-c\ell_0-\alpha\delta,c\ell_0+ \alpha\delta]\, .
\end{cases}
\end{equation}

Consider the function $\tilde{g}(z)=g(z)\cdot\chi({\rm Im}(z))$
%$g_1(z)=g(z)\cdot \chi_0(\pm \frac{{\rm %Im}(z)-c\ell_0}{\alpha\delta})$
and the
continuous $(0,1)$-form $\varphi\stackrel{def}= \bar{\partial}{\tilde{g}}$ on  $R_{\alpha\delta}$. The form $\varphi$ vanishes outside $Q_{\alpha\delta}^{\pm}$.
Let $\varepsilon$
%on $R_{\delta}$
be a small positive number that will be chosen later.
Consider the measurable $(0,1)$-form $\varphi_{\varepsilon}$ on $R_{\delta}$ that equals $\varphi$ on $Q_{\alpha\delta,\varepsilon}^{\pm}\stackrel{def}=
\{x+iy: x\in (-{\alpha\delta\varepsilon},{\alpha\delta\varepsilon}), \pm y\in
(c\ell_0, c\ell_0 + \alpha \delta)\}$
%$\cap R_{\varepsilon\alpha\delta}$
and vanishes outside this set.
Extend its pullback under the conformal mapping $\omega: R^X_{\delta}\to R_{\delta}$
%$\varphi_{\varepsilon}$
%and extend the obtained form
to a measurable $(0,1)$-form on $X$ by putting it equal to zero outside $R^X_{\delta}$. Denote the obtained form by $\varphi_{\varepsilon}^X$.

By Corollary 2.14.2
%p. 82 existence of strictly subharmonic exhaustion function
of \cite{Na} there exists a strictly subharmonic exhaustion function $\psi$ on $X$.
The $L^2$-norm of $\varphi_{\varepsilon}$ with respect to the Euclidean metric on the complex plane does not exceed $C_2 \sqrt{\varepsilon}$ for an absolute constant $C_2$.
Hence, the weighted $L^2$-norm on $X$ of $\varphi^X_\varepsilon$ with respect to the K\"ahler metric and the weight $e^{-\psi}$ (see Definition 2.6.1
%p.54
of \cite{Na}) does not exceed $C_3 C_2 \sqrt{\varepsilon}$ for a constant $C_3$ that depends on $\psi$ and
on the K\"ahler metric on a the compact subset
$R^X_{\delta}$ of $X$. By Corollary 2.12.6 of \cite{Na}
%p.74, existence of $solution in weighted $L^2$-space
there exists a function $f^X$ with $\bar{\partial} f^X=\varphi_X$ in the weighted $L^2$-space on $X$ with respect to the K\"ahler metric and the weight $e^{-\psi}$ (see Definition 2.6.1
%p.54, definition of weighted $L^"$-spaces for functions, definition 2.8.1 p.64 $\theta_{\varphi}$, p.73 %$\theta_{\omega}$,
of \cite{Na}), whose norm in this space does not exceed $C_4 C_3 C_2 \sqrt{\varepsilon}$ for a constant $C_4$ depending only on $X$, $\psi$, and the K\"ahler metric. Let $(Q_{\alpha\delta,\varepsilon}^{\pm})^X$ be the preimages of $Q_{\alpha\delta,\varepsilon}^{\pm}$ under $\omega$.
The function $f^X$ is holomorphic on
$X\setminus \big((Q_{\alpha\delta,\varepsilon}^{+})^X\cup (Q_{\alpha\delta,\varepsilon}^{-})^X\big)$.
Put $\tilde{Q}_{\delta}^{\pm}\stackrel{def}=\{x+iy\in R_{\delta}: \pm y\in (c\ell_0 -\delta, c\ell_0 + 2\delta)\}$, and $(\tilde{Q}_{\delta}^{\pm})^X=\omega^{-1}(\tilde{Q}_{\delta}^{\pm})$. Then $(Q_{\alpha\delta,\varepsilon}^{\pm})^X$ is relatively compact in $(\tilde{Q}_{\delta}^{\pm})^X$.
On a relatively compact open subset of $X$, containing the closed subset $\overline{S_{\delta_0}} \setminus  ((\tilde{Q}_{\delta}^+)^X \cup  (\tilde{Q}_{\delta}^-)^X)$ of $X$, the supremum norm of $|f^X|$ is estimated by its weighted $L^2$-norm: $|f^X|< C_5 \sqrt{\varepsilon}$ in a neighbourhood of  $\overline{S_{\delta_0}} \setminus  ((\tilde{Q}_{\delta}^+)^X \cup  (\tilde{Q}_{\delta}^-)^X)$ for a constant $C_5$ that depends on the K\"ahler metric, on $\psi$ and on the constants chosen before (see Theorem 2.6.4
%p. 59 sup norm estimate by $L^2$-norm
of \cite{Na}).

On the other hand the classical Cauchy-Green formula on the complex plane provides a solution $\tilde f$ of the equation $\overline{\partial}\tilde{f}=\varphi_{\varepsilon}$ on the set $\tilde{Q}_{\delta}^+ \cup  \tilde{Q}_{\delta}^-$. The supremum norm of the function $\tilde{f}$ is estimated by $C_6\sqrt{\varepsilon}$ for an absolute  constant $C_6$.
Let $\tilde{f}^X$ be the pullback of $\tilde f$ to  $(\tilde{Q}_{\delta}^+)^X \cup  (\tilde{Q}_{\delta}^-)^X$. The function $f^X - \tilde{f}^X$ is holomorphic on $(\tilde{Q}_{\delta}^+)^X \cup  (\tilde{Q}_{\delta}^-)^X$ and satisfies the inequality $|f^X - \tilde{f}^X|<(C_5+C_6)\sqrt{\varepsilon}$ at all points of the set $(\tilde{Q}_{\delta}^+)^X \cup  (\tilde{Q}_{\delta}^-)^X$, that are close to its boundary. Hence, the inequality is satisfied on $(\tilde{Q}_{\delta}^+)^X \cup  (\tilde{Q}_{\delta}^-)^X$. As a consequence, $|f^X|<(C_5+2C_6)\sqrt{\varepsilon}$ on $(\tilde{Q}_{\delta}^+)^X \cup  (\tilde{Q}_{\delta}^-)^X$.

Choose $\varepsilon$ depending on $C_5$ and $C_6$, so that
\begin{align}\label{eqabcd'}
|f^X|<\min\Big\{{\rm dist}(G,\{-1,1\}),\,\frac{1}{2}\Big\} \;\mbox{ on} \; S_{\sigma_0}\,.
\end{align}
Put $\sigma=\varepsilon\alpha \delta$. Consider the smooth function $g_{\sigma}^X$ on  $S_{\alpha\delta\varepsilon }=S_{\sigma}$ which equals $\tilde{g}^X$
on $R_{\delta}^X\cap S_{\sigma }$ and vanishes on the rest of $S_{\sigma}$. Hence, it vanishes on all circles of the bouquet except $\gamma_0$, and therefore, the monodromy of its homotopy class along each
such circle is the identity. The restriction of $g_{\sigma}^X$ to $\mathring{R}_{\delta}^X\cap S_{\sigma}$ represents the element $a_1^{\pm} a_2^{\pm}\ldots a_2^{\pm}\in\pi_1(\mathbb{C}\setminus\{-1,1\}, (-1,1))$. Moreover, on $\big({R}_{\delta}^X\setminus \mathring{R}_{\delta}^X\big)\cap S_{\sigma} $  the inequality $|g_\sigma^X|<\frac{1}{2} $ holds, and on $S_{\sigma}\setminus {R}_{\delta}^X$ the mapping $g_{\sigma}^X$
vanishes. Hence, the monodromy of the homotopy class of $g_{\sigma}^X$ along $\gamma_0$ equals $a_1^{\pm} a_2^{\pm}\ldots a_2^{\pm}$.  By the inequality \eqref{eqabcd'} the monodromies of the homotopy class of $g_{\sigma}^X   -f^X$
along all circles of the bouquet coincide with those of the homotopy class of $g_{\sigma}^X$.
The function
$g_{\sigma}^X   -f^X$ is holomorphic on $S_{\sigma }$.
We put
$C_7= \frac{2 \log 2 c\alpha\varepsilon}{C}$.
For each positive $\sigma< \epsilon\alpha\sigma_0$ we found no less than
$\frac{1}{4} 2^{\frac{C_7 \ell_0}{\sigma}}$ irreducible non-homotopic holomorphic mappings from $S_{\sigma}$ to $\mathbb{C}\setminus\{-1,1\}$.
The proposition is proved. \hfill $\Box$

\end{document}